\documentclass[12pt]{article}
\usepackage{amsfonts}
\usepackage{epsfig}
\title{The Plaid Model and Outer Billiards on Kites}
\author{Richard Evan Schwartz \thanks{\hskip 5 pt Supported by 
N.S.F. Research Grant DMS-1204471}}

\newtheorem{theorem}{Theorem}[section]

\newtheorem{lemma}[theorem]{Lemma}

\def\startproof{{\bf {\medskip}{\noindent}Proof: }}

\def\endproof{$\spadesuit$  \newline}

\def\Q{\mbox{\boldmath{$Q$}}}%
\def\R{\mbox{\boldmath{$R$}}}%
\def\Z{\mbox{\boldmath{$Z$}}}%

\begin{document}
\maketitle

\begin{abstract}
This paper establishes an affine
quasi-isomorphism between the
arithmetic graph associated to
outer billiards on kites, and the
plaid model, a combinatorial 
construction for producing embedded
polygons in the plane.
\end{abstract}

\section{Introduction}

The purpose of this paper,
which continues the study of the plaid model,
is to prove the Quasi-Isomorphism Conjecture
from [{\bf S1\/}].

\begin{theorem}[Quasi-Isomorphism]
Let $p/q \in (0,1)$ be any rational such
that $pq$ is even.  Let $\Pi$ denote the
plaid model at the parameter $p/q$.
Let $\Gamma$ denote the arithmetic graph
associated to the special orbits
of outer billiards on the kite $K_{p/q}$.
Then $\Pi$ and $\Gamma$ are 
$2$-affine-quasi-isomorphic.
\end{theorem}

\noindent
{\bf Quasi-Isomorphisms:\/}
First let me explain what the
term {\it quasi-isomorphic\/} means.
We say that two embedded polygons are
$C$-{\it quasi-isomorphic\/} if they can be
parametrized so that corresponding
points are within $C$ units of each other.
We say that two unions $\Gamma$ and $\Pi$ of polygons
are $C$-{\it quasi-isomorphic\/} if there is a bijection
between the members of $A$ and the members
of $\Gamma$ which pairs up $C$-quasi-isomorphic polygons.
Finally, we say that $A$ and $B$ are
$C$-{\it affine-quasi-isomorphic\/}
if there is 
an affine transformation $T: \R^2 \to \R^2$
such that $\Pi$ and $T(\Gamma)$ are $C$-quasi-isomorphic.
\newline
\newline
{\bf Example:\/}
Figure 1 shows 
the Quasi-Isomorphism Theorem in action
for the parameter
$p/q=3/8$.  The black polygons are part of the plaid model
and the 
grey polygons are part of the affine image of the
arithmetic graph. 

\begin{center}
\resizebox{!}{5.3in}{\includegraphics{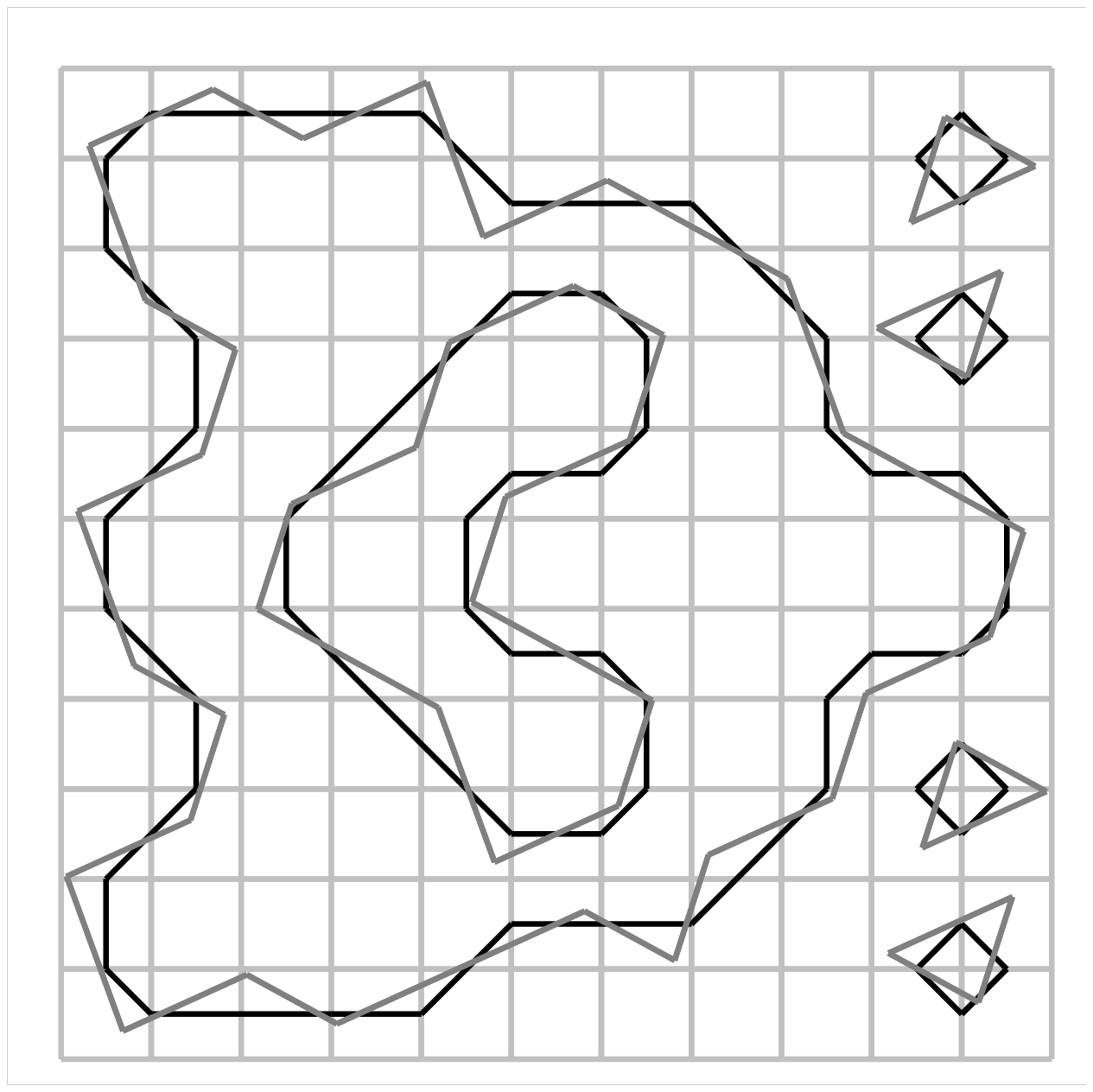}}
\newline
{\bf Figure 1:\/} The Quasi-Isomorphism Theorem
in action.
\end{center}

In general, the vertices of the black curves lie
in the lattice $\frac{1}{2}\Z^2$, which has co-area $1/4$,
whereas the vertices
of the grey curves lie in a lattice having co-area $1+p/q$.
Even through the polygons line up in a global
sense, their vertices lie on quite different lattices.
This fact points to the highly nontrivial
nature of quasi-isomorphism being established.
\newline
\newline
{\bf The Plaid Model:\/}
The plaid model, introduced in [{\bf S1\/}], has a
combinatorial flavor that
is similar in spirit to corner percolation and
P. Hooper's Truchet tile system [{\bf H\/}]. I
will give a precise definition in \S 2.  The plaid
model has a rich and inherently hierarchical
structure which I explored in [{\bf S1\/}] and
[{\bf S2\/}].
In [{\bf S2\/}] I showed that the plaid model (which
can be defined even for irrational parameters) has
unbounded orbits for every irrational parameter.
The significance of the Quasi-Isomorphism Theorem,
aside from its intrinsic beauty, is that it allows
one to use the plaid model as a tool for studying
the dynamics of outer billiards on kites.  
\newline
\newline
{\bf Outer Billiards:\/}
Outer billiards was introduced by B. H. Neumann in the
late $1950$s. See [{\bf N\/}].
One of the central questions about outer billiards has
been: Does there exist
a compact convex set $K$ (other than a line segment) and a point
$p_0$ so that the orbit $\{p_n\}$ exits every compact
set?  This question, dating from around $1960$, is called
the Moser-Neumann question.  See  [{\bf M1\/}] and [{\bf M2\/}].
  My monograph [{\bf S0\/}] has a
survey of the known results about this problem. 
See also [{\bf T1\/}] and [{\bf T2\/}]. The main
work on this Moser-Neumann problem is contained in 
[{\bf B\/}], [{\bf D\/}], [{\bf VS\/}], [{\bf K\/}], [{\bf GS\/}], [{\bf G\/}],
[{\bf S0\/}], [{\bf DF\/}].

In [{\bf S0\/}] I showed that
there are unbounded orbits with respect to any irrational kite.
Any kite is affinely equivalent to the kite
$K_A$ with vertices
\begin{equation}
\label{KITE}
(-1,0), \hskip 30 pt
(0,1), \hskip 30 pt (0,-1), \hskip 30 pt (A,0),
\hskip 30 pt A \in (0,1).
\end{equation}
$K_A$ is rational if and only if $A \in \Q$.
Outer billiards on
$K_A$ preserves the infinite family of
horizontal lines $\R \times \Z_1$.
Here $\Z_1$ is the set of odd integers.
A {\it special orbit\/} is an orbit which
lies in this set of horizontal lines.
\newline
\newline
{\bf Polytope Classifying Spaces:\/}
Say that a {\it grid\/} is an affine image
of $\Z^2$ in the plane.
Suppose that $G$ is a grid and
$\Gamma$ is a collection of polygons and polygonal
arcs having vertices in $G$.  Suppose also
that there are only finitely many local pictures
for $\Gamma$, up to translation.

A {\it polytope classifying space\/} is a
triple $(X,\Phi,{\cal P\/})$ where
$X$ is some Euclidean manifold,
$\Phi: G \to X$ is an affine map, and
$\cal P$ is a partition of $X$ into
polytopes, each one labeled by one of
the local pictures for $\Gamma$.
We say that the triple $(X,\Phi,{\cal P\/})$ is
a classifying space {\it for\/} $\Gamma$ if
the following is true for each $c \in G$.
The point $\Phi(c)$ lies in the interior
of one of the polytopes of $\cal P$ and
the local picture of $\Gamma$ at $c$
coincides with the label of the
polytope containing $\Phi(c)$.
\newline
\newline
{\bf The Plaid Master Picture Theorem:\/}
Our Isomorphism Theorem from
[{\bf S1\/}] says that, for any even
rational parameter $A \in (0,1)$, the
plaid model has a $3$ dimensional
classifying space
$(X_A,\Phi_A,{\cal P\/}_A)$. 
The grid in question here is the
set of centers of integer unit squares.
We call this grid the {\it plaid grid\/}.

In hindsight, it would have been better to
call the Isomorphism Theorem the Plaid
Master Picture Theorem, 
and so I will rechristen it here.

The plaid classifying spaces fit together.
There is a $3$ dimensional group 
$\Lambda$ of affine transformations
acting on $$\widehat X=\R^3 \times [0,1],$$ together
with a $\Lambda$ invariant partition
partition of $\widehat X$ into
integer convex polytopes such that
the intersection of these objects
with $\R^3 \times \{P\}$ gives the
pair $(X_A,{\cal P\/}_A)$.  Here
$P=2A/(1+A)$.  

I will also explain how to interpret the
quotient $\widehat X/\Lambda$, together
with the quotient partition, as a
$4$ dimensional affine polytope
exchange transformation.  This
is the same as the
{\it curve following dynamics\/} 
defined in  [{\bf H\/}].
\newline
\newline
{\bf The Graph Master Picture Theorem:\/}
Our (Graph) Master Picure Theorem from
[{\bf S0\/}] does 
for the arithmetic graph what the
Plaid Master Picture Theorem does for the
plaid model.

There are several minor differences
between the plaid case and the graph case.
First, the way the partition determines
the local picture is different; in fact
there are really two partitions in
this case.  Another difference is the
way that the slices fit together. 
The pair $(X_A,{\cal Q\/}_A)$ associated
to the arithmetic graph at $A$ is
the slice at $A$, rather than the
slice at $P$.  This seemingly
minor point turns out to be 
important and helpful.

Another difference
is that we will compose the classifying
map, $\Phi'_A$ with the affine
map from the Quasi-Isomorphism Theorem,
so that the Master Picture Theorem
pertains directly to the set of
polygons we wish to compare with the
plaid model. We define the
{\it graph grid\/} to be the affine
image of $\Z^2$ under this affine map.

The classifying space picture for the
arithmetic graph can also be interpreted
in terms of affine polytope exchange 
transformations, but we will not dwell
on this point.  In the plaid case
we really need the interpretation for
some of our constructions, but in the
graph case we can do without it.
\newline
\newline
{\bf The Projective Intertwiner:\/}
Let $G_{\Pi}$ and $G_{\Gamma}$ respectively denote the
plaid grid and the graph grid.  Let
$X_{\Pi}$ and $X_{\Gamma}$ respectively denote the
plaid PET and the graph PET.  Let
$\Phi_{\Pi}: G_{\Pi} \to X_{\Pi}$ and
$\Phi_{\Pi}: G_{\Gamma} \to X_{\Gamma}$ be the
two intertwining maps.  Technically, we should
write $\Phi_{\Pi}=\Phi_{\Pi,A}$, and $G_{\Gamma,A}$, etc.,
to denote dependence on the parameter.

The grid $G_{\Gamma}$ has the property that no point
of $G_{\Gamma}$ lies on the boundary of an integer unit
square.  This means that to each $b \in G_{\Gamma}$,
there is a unique associated point $c(b) \in G_{\Pi}$.
Here $c(b)$ is the center of the unit integer
square containing $b$. We call $b \to c(b)$ the
{\it centering map\/}.  We will find a piecewise
projective map $\Psi: X_{\Pi} \to X_{\Gamma}$ with the following
property:
\begin{equation}
\label{intertwine}
\Phi_{\Gamma}(b)=\Psi \circ \Phi_{\Pi}(c(b)), \hskip 30 pt
\forall b \in G_{\Gamma}.
\end{equation}
This equation is meant to hold true for each
even rational parameter, but the map
$\Psi$ is globally defined and piecewise projective.
That is, the restriction of $\Psi$ to the interiors of
suitable polytopes is a projective transformation.
We call $\Psi$ the {\it projective intertwiner\/}.
\newline
\newline
{\bf Main Idea:\/}
A plaid $3$-{\it arc\/}
is the portion of a connected component of the plaid
model which occupies
$3$ consecutive integer unit squares.
We will see that the plaid classifying space
$X_{\Pi}$ has a partitioned into finitely many polytopes,
each one labeled by a plaid $3$-arc. (This is
where we use the dynamical interpretation.)
For each point $a \in G_{\Pi}$, the plaid $3$-arc
centered at $a$ is a translate of the label
of the partition piece containing
$\Phi_A(a)$.   We call this auxiliary partition the
{\it plaid triple partition\/}.

At the same time the graph classifying space
$X_{\Gamma}$ has a partition into finitely
many polytopes, with the following properties:  Each of
these polytopes is labeled by a (combinatorial) length $2$ polygonal
arc with integer vertices.  
For each lattice point
$b \in G_{\Gamma}$, the length $2$ arc of the arithmetic
graph centered at $b$ is determined
by the polytope label.
We call this partition the {\it graph double partition\/}.

The way we prove the Quasi-Isomorphism Theorem is that
we consider the image of the plaid triple partition under
the projective intertwining map. We study how it sits
with respect to the graph double partition. This allows
us to line up the plaid model with the arithmetic
graph piece by piece.  We will see also that there is
some amount of fussing when it comes to globally 
fitting the different pieces together.
\newline
\newline
{\bf Self-Containment:\/}
This paper is the culmination of several years of work,
and about $8$ years of thinking, about the structure of
outer billiards on kites.  The reader might wonder whether
it is possible to read this paper without an enormous
buildup of machinery and background.  In one sense this
paper is very far from self-contained, and in another sense,
it is self-contained.  

The proof of the Quasi-Isomorphism Theorem, as stated,
is far from self-contained.  It relies on the
Plaid and Graph Master Picture Theorems.
In this paper, I will state these results precisely,
but I will not re-do the proofs.
On the other hand, I will give a complete account of the
relationship between the graph PET and the plaid PET. 
So, if one view the Quasi-Isomorphism Theorem as a statement
about some kind of quasi-isomorphism between $4$ dimensional
affine PETs, then the paper is self-contained.  The
projective intertwining map $\Psi$ is not a conjugacy
between these two PETs, but for even rational parameters,
it does induce a bijection between the orbits of the two
systems which preserves periods up to a factor of $2$.
This seems to be a striking situation. I haven't
seen anything like it in the study of PETs.
\newline
\newline
{\bf Paper Overview:\/}
Here is an overview of the paper.
\begin{itemize}
\item In \S 2, I will recall the definition of
the plaid model.
\item In \S 3, I will recall the definition of the
arithmetic graph.  I will also prove some geometric
facts about the {\it graph grid\/}, which is
the image of $\Z^2$ under the affine transformation
implicit in the Quasi-Isomorphism Theorem.
\item In \S 4 I will deduce the Quasi-Isomorphism
from two other results, the Pixellation Theorem
and the Bound Chain Lemma.  The Pixellation Theorem,
a $6$-statement result,
is the main thrust of the paper. It makes a statement
to the effect that locally the polygons in the
two collections follow each other.
The Bound Chain Lemma provides the glue
needed to match the local pictures together.
\item In \S 5 I will deduce the Bound Chain
Lemma from the Pixellation Theorem.
The deduction is a tedious case by case
analysis.
Following \S 4-5, the rest of the paper
is devoted to proving the Pixellation Theorem.
\item In \S 6 I will give an account of the
Plaid Master Picture Theorem, and I will
also define the Plaid Triple Partition.
\item In \S 7, I will give an account of
the Graph Master Picture Theorem.
\item In \S 8 I will introduce a formula,
called the Reconstruction Formula, which gives
information about the graph grid in terms of
the graph classifying map. I will use this
formula to prove Statement 1 of the
Pixellation Theorem modulo 
some integer computer calculations.
\item In \S 9 I will prove the existence
of the projective intertwiner, the map
which relates the Plaid and Graph Master Picture
Theorems.  The result is summarized in the
Intertwining Theorem.
\item In \S 10, I will combine the
Plaid Master Picture Theorem,
the Graph Master Picture Theorem, and
the Intertwining Theorem to prove
Statment 2 of the Pixellation
Theorem, modulo some integer
computer calculations. I will also set up the
problems which one must
solve in order to Prove the
remaining statements of
the Pixellation Theorem.
\item In \S 11, I will introduce some auxiliary
sets which augment the power of the machinery
developed above.  I will also
give some sample arguments which illustrate
most of the features of the general proof
given in \S 12.
\item In \S 12 I will combine all the
machinery to prove Statements 3,4, and 5 of the
Pixellation Theorem, modulo some
integer computer calculations.
\item In \S 13 I will explain the integer
computer calculations used in 
\S 8,10,11.  In \S \ref{COMP}, I
will explain how to access the computer
program which runs the proof.
\end{itemize}

\noindent
{\bf Companion Program:\/}
This paper has a companion computer program.
I checked essentially all the details of
the paper on the computer program.  I think
that the reader would find this paper much easier
to read if they also used the program.  Most
of the technical details of the paper are simply
encoded into the computer program, where they
are illustrated as colorful plots and pictures.
In particular, the polytope partitions which form the
backbone of our proof are much easier to understand
when viewed with the program.
Again, in \S \ref{COMP}, I explain
how to access the program. 
\newline
\newline
{\bf Acknowledgements:\/}
I would like to thank Pat Hooper,
Sergei Tabachnikov, and Barak Weiss
for helpful and interesting conversations
related to this work.

\newpage

\section{The Plaid Model}

\subsection{The Mass and Capacity Labels}
\label{def1}

In [{\bf S1\/}] we gave a general definition
of the plaid model and explored its symmetries.
Rather than repeat the definition from
[{\bf S1\/}] we will explain the plaid
model for the parameter $p/q=2/9$.  The
definition generalizes in a straightforward
way to the other parameters.  The parameter
$2/9$ is small enough to allow us to
present the whole definintion concretely
but large enough so that the definition isn't
muddled by too many coincidences between 
small numbers.  In any case, the planar definition
of the plaid model is just given for the sake
of completeness; we will use the PET definition
given in \S 3.

So, we fix
$p/q=2/9$.  We have the following auxiliary
parameters:
\begin{equation}
\omega=p+q=11, \hskip 20 pt
P=\frac{2p}{p+q}=\frac{4}{11}, \hskip 20 pt
Q=\frac{2q}{p+q}=\frac{18}{11}, \hskip 20 pt
\tau=3.
\end{equation}
The parameter $\tau$ is the solution in
$(0,\omega)$ to the equation $2p\tau \equiv 1$ mod $\omega$.
In this case, we have $2(2)(3)=12 \equiv 1$ mod $11$.

The plaid model is invariant under the action of the
lattice generated by vectors $(0,11)$ and $(11^2,0)$.
Moreover, the model does not intersect the boundaries
of the $11$ {\it blocks\/}
$$[0,11]^2 + (11k,0), \hskip 30 pt k=0,...,10.$$
So, it suffices to describe the plaid model in
each of these $11$ blocks.  

We think of $[0,11]^2$ as tiled by the $12 \times 12$
grid of unit integer squares.  We introduce the
{\it capacity labeling\/} on the integer points
along the coordinate axes as follows.
\begin{itemize}
\item We give the labels $0,2,4,6,8,10,-10,-8,-6,-4,-2,0$ to
the points $(0,3)$,$(0,6)$,$(0,9),...,(0,33)$.
\item We extend the labels so that they are periodic
with respect to translation by $(0,11)$.
\item We give the same labels to integer points along the
$x$-axis, simply interchanging the roles of the $x$ and
$y$ axes.
\end{itemize}

We now describe the {\it mass labeling\/} of the integer
points along the $y$-axis. We let $[x]_{11}$ denote
the point $x+11k$ which lies in $(-11,11)$.  Here
$k \in \Z$ is a suitably chosen integer. When
$x \in 11\Z$ we set $[x]_{11}=\pm 11$. (This
situation does not actually arise.)
 If $\mu \not = 0$ is a capacity label of some point,
then exactly one of the two quantities
\begin{equation}
[\mu/2]_{11}, \hskip 30 pt
[(\mu/2)+11]_{11}
\end{equation}
is odd.  We define the mass label of the point
to be this odd quantity.  When
$\mu=0$ we define the mass label to be $\pm 11$,
without specifying the sign.  We denote this
unsigned quantity by $[11]$.

Figure 2.1 shows the capacity and mass labelings of
the first block.

\begin{center}
\resizebox{!}{4.5in}{\includegraphics{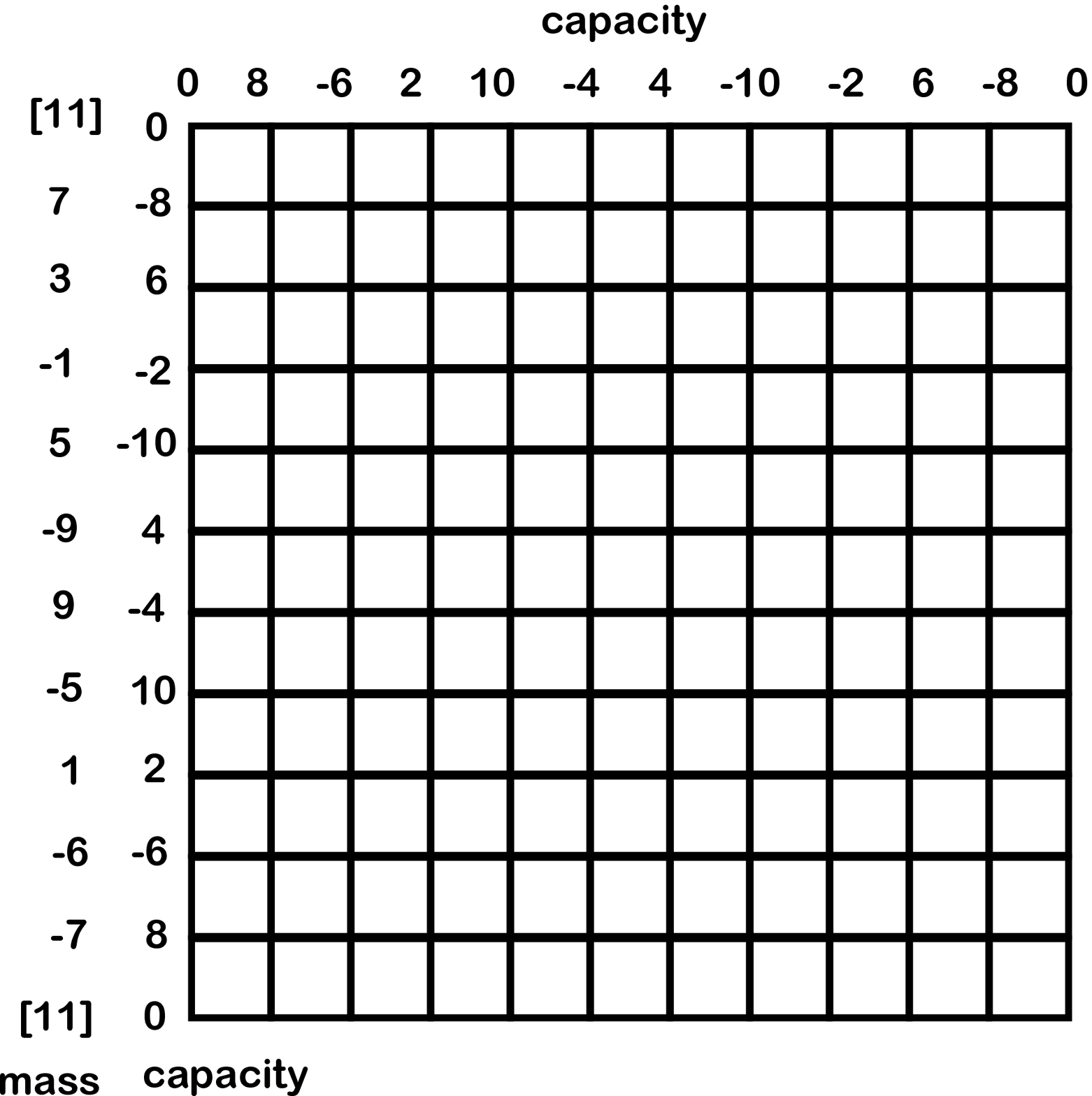}}
\newline
{\bf Figure 2.1:\/} The mass and capacity labelings for the
first block.
\end{center}

\subsection{Definition of the Model}

We define $6$ families of lines
\begin{itemize}

\item $\cal H$: Horizontal lines
intersecting the $y$-axis at integer points.

\item $\cal V$:  Vertical lines
intersecting the $x$-axis at integer points.

\item ${\cal P\/}_-$:  Lines of slope
$-4/11$ intersecting the $y$-axis at integer points.

\item ${\cal Q\/}_-$: Lines of slope
$-18/11$ intersecting the $y$-axis at integer points.

\item ${\cal P\/}_+$: Lines of slope
$+4/11$ intersecting the $y$-axis at integer points.

\item ${\cal Q\/}_+$:  Lines of slope
$+18/11$ intersecting the $y$-axis at integer points.
\end{itemize}
We call the $\cal H$ and $\cal V$ lines
{\it grid lines\/}.
We call the lines from the latter
$4$ families the {\it slanting lines\/}.
The grid lines divide
$\R^2$ into the usual grid of unit squares.
We call the edges of these equares
{\it unit grid segments\/}.  When we
want to be more specific, we will
speak of {\it unit horizontal segments\/}
and {\it unit vertical segments\/}.

Say that an {\it intersection point\/} is
a point on a unit grid segment which is
also contained on a slanting line.
Every intersection point is actually
contained on two slanting lines, one
of positive slope and one of negative slope.
With the correct multiplicity convention.
every unit grid segment contains $2$
intersection points. The multiplicity
convention is that the horiozontal unit segments
which cross the lines $x=(11/2)+11k$ for
$k \in \Z$ are declared to
contain $2$ intersection points  What is going on
is that these segments contain a single point
which lies on $4$ slanting lines. 
we call these intersection points
{\it light\/} and {\it dark\/}
according to the following rules:
\newline
\newline
{\bf Vertical Rule\/}
Suppose that $z$ is an intersection point
contained in a vertical unit line segment.
We call $z$ a {\it light point\/} if
and only if the $3$ lines through $z$
all have the same sign, and the common
sign is usual.  
\newline
\newline
{\bf Horizontal Rule:\/}
If $z$ happens to lie on the special segments
mentioned in the section on the hexagrid -- i.e.
if $z$ is an intersection point of multiplicity two -- then
$z$ is not a light point.  Otherwise, $z$ is a
light point if and only if the signs of the three
lines are usual, and the sign of the horizontal
line matches the sign of the line of negative slope
and is opposite from the sign of the line of
positive slope.
\newline
\newline
{\bf Remark:\/}
It might happen that the point $z$ is a vertex of one
of the unit squares.  In this case, we apply each
rule separately.  It turns out that this situation
happens only when $z$ lies on a vertical line of
the form $x=n\omega$ for $n \in \Z$.  In this case,
only the horizontal rule can make $z$ a light point
because the vertical line just mentioned has an unusual
sign.
\newline

Now we define the plaid model.
We say that a unit grid segment is {\it on\/} if it
contains exactly $1$ light point.  We proved in
[{\bf S1\/}] that each unit square has either $0$ or
$2$ edges which are on. (This result takes up the
bulk of [{\bf S1\/}].)  We form the {\it plaid model\/}
as follows:  In each square which has $2$ on edges, we
connect the centers of these edges by a segment.
By construction, the union of all these connector
segments is a countable union of embedded polygons.
This union of polygons is the plaid model for the
given parameter.

\begin{center}
\resizebox{!}{4.2in}{\includegraphics{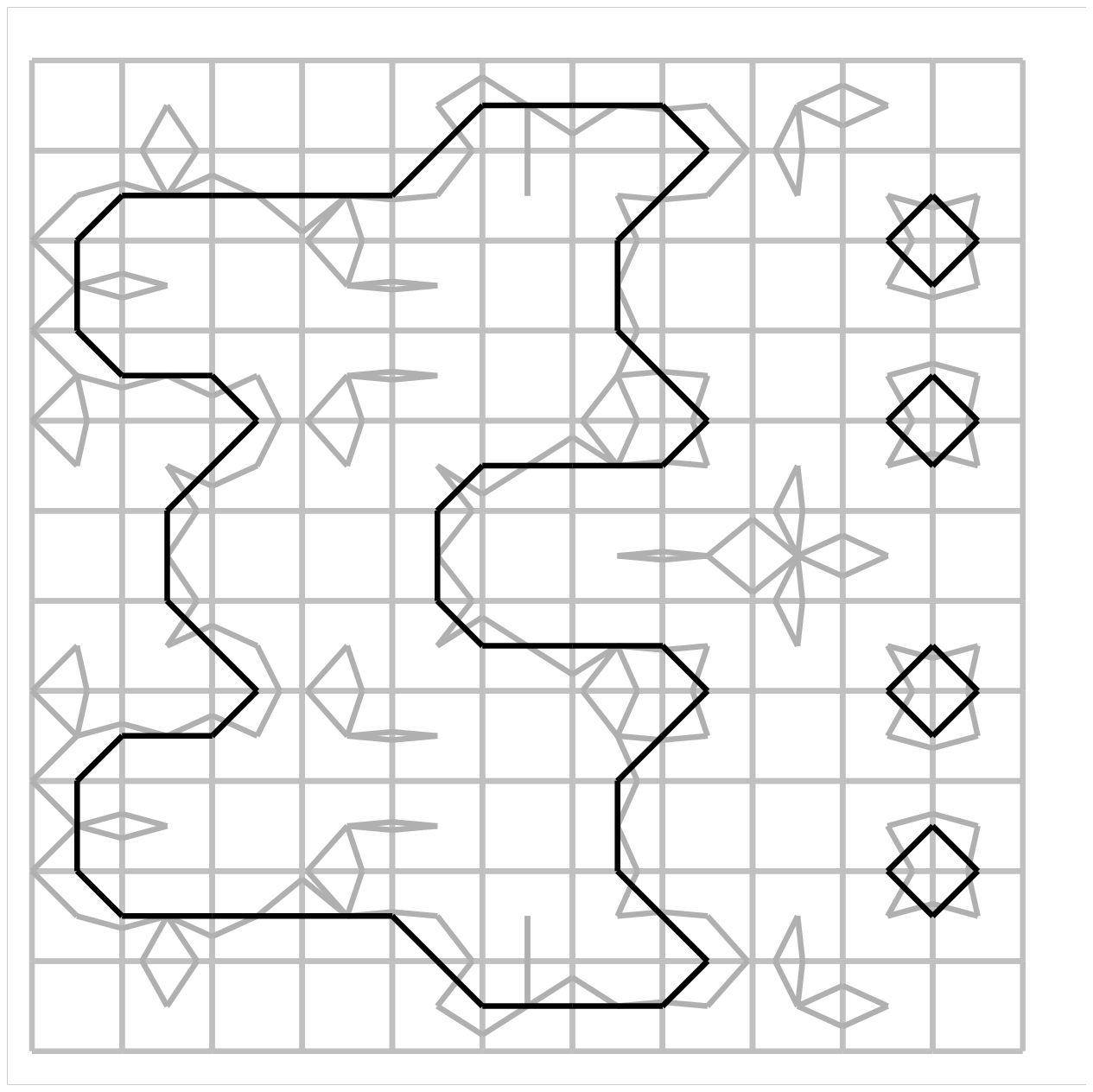}}
\newline
{\bf Figure 2.2:\/} The first block of the plaid model
for $p/q=2/9$.
first block.
\end{center}

Figure 2.2 needs some explanation.  The black polygons
are the plaid polygons.  We have illustrated the
light points in each unit square by drawing a grey
arc from the center of the square to the corresponding
light points.  Notice that none of the black arcs 
touch edges which have two grey arcs incident to them.

It is worth pointing out that the picture looks
different in different blocks.  Figure 2.3 shows
the picture in the block
$$[0,11]^2+(55,0).$$

\begin{center}
\resizebox{!}{4.5in}{\includegraphics{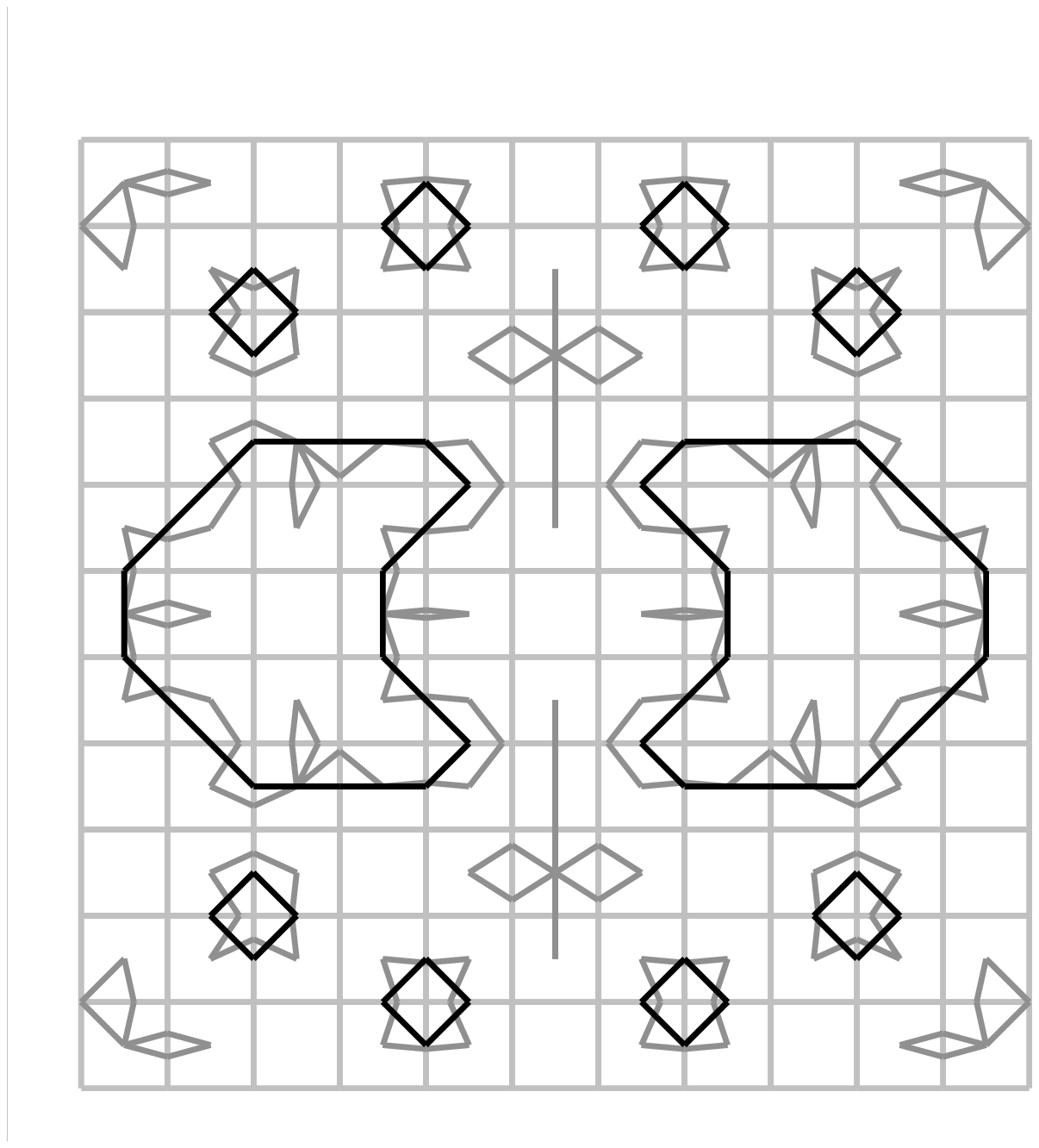}}
\newline
{\bf Figure 2.2:\/} The plaid model for
for $p/q=2/9$ in $[55,66]\times[0,11]$.
\end{center}
The center of rotational symmetry is
$(11^2/2,11/2)$.

The plaid model for this parameter is
invariant under the reflections in the
lines $y=11/2$ and $x=11^2/2$.  In
particular, the point $(11^2/2,11/2)$ is
a center of rotational symmetry. All
this is proved, for the general parameter,
in [{\bf S1\/},\S 2].

\subsection{Alternate Definition}

There is a second definition which just uses the
lines in the families
families $\cal V$, $\cal H$, ${\cal P\/}_-$
and ${\cal Q\/}_-$. 
In this definition, an intersection point is light
if and only if it is the intersection point of a
grid line a negative slanting line such that
the mass of the slanting line is less than the
capacity of the grid line and the signs of
the two lines are the same.

This second definition points out the hierarchical
nature of the plaid model.  The plaid model
has $4$ grid lines of capacity $2k$, two horizontal
and two vertical, and each of them carries at most
$2k$ light points.  This hierarchical structure is
exploited in [{\bf S1\/}] and [{\bf S2\/}] to get
control over the geometry of the components of
the plaid model.

It is not at all clear from either definition that the
number of on sides of an integer unit square is either
$0$ or $2$.  However, we proved this in
[{\bf S1\/}].  It is a consequence of the Isomorphism
Theorem.

\newpage

\section{The Arithmetic Graph}
\label{ob}

\subsection{Special Orbits and the First Return Map}

We will describe the picture when
$p/q \in (0,1)$ is rational.  We need not
take $pq$ even.  Let $A=p/q$. 
We consider outer billiards on the kite $K_A$,
described in Equation \ref{KITE} and shown
in Figure 3.1.
Recall from the introduction that the {\it special orbits\/}
are the orbits which lie on the set $S$ of
horizontal lines having
odd integer $y$-intercepts. 

Let $\psi_A$ denote the
second iterate of the outer billiards map.
This map is a piecewise translation.
Let $\Psi_A$ denote the
first return map of $\psi_A$ to the union 
$\Xi$ of $2$ rays shown in Figure 3.1 below
Here
\begin{equation}
\Xi=\R_+ \times \{-1,1\} \subset S
\end{equation}

\begin{lemma} 
Every special orbit is combinatorially identical
to the orbit of a point of the form
\begin{equation}
\label{zoop}
\bigg(2mA+2n+\frac{1}{q},\pm 1\bigg).
\end{equation}
Here $m,n \in \Z$ and the left hand side of the
equation is meant to be positive.
\end{lemma}

\startproof
This is proved in [{\bf S1\/}]. Here is a sketch of
the proof.
What is going on is that $\psi_A$ permutes the
intervals of a certain partition of $S$ into
intervals of length $2/q$, namely
$$S=\bigcup_{{m,n} \in \Z^2} \bigg(\frac{2m}{q},n\bigg)$$
and so every orbit is combinatorially identical
to the orbit of a center point of one of these
intervals.  Moreover, all such orbits intersect
the set $\Xi$.  When the orbit of a center
point of an interval in the partition intersects
$\Xi$, it does so in a point of the form
given in Equation \ref{zoop}.  The left hand side
of Equation \ref{zoop} is an alternate way of
writing the set 
$$\bigg( \bigcup_{m \in \Z} \frac{2m}{q}\bigg) \cap \R_+.$$
\endproof

\noindent
{\bf Remark:\/} The same result holds if $1/q$ is in
Equation \ref{zoop} is replaced
by any $\iota \in (0,2/q)$.  We will take
$\iota=1/q$ when we state the Master
Picture Theorem.

\subsection{The Arithmetic Graph}

Suppose that $m_0,n_0$ are integers such that
$2m_0A+2n_0 \geq 0$.  Then, by definition,
there are integers $m_1,n_1$, with
$2m_1A+2n_1 \geq 0$ such that

\begin{equation}
\label{ag1}
\Psi_A\Big(2m_0A+2n_0+\frac{1}{q},1\bigg)=
\Big(2m_1A+2n_1+\frac{1}{q},\epsilon \bigg),
\end{equation}
where $\epsilon \in \{-1,1\}$ is given by
\begin{equation}
\epsilon=(-1)^{m_0+m_1+n_0+n_1}.
\end{equation}
Reflection in the $x$-axis conjugates $\Psi_A$
to $\Psi_A^{-1}$, so we also have
\begin{equation}
\Psi_A^{-1}\Big(2m_0A+2n_0+\frac{1}{q},-1\bigg)=
\Big(2m_1A+2n_1+\frac{1}{q},-\epsilon \bigg),
\end{equation}
Actually, we won't end up caring about
$\epsilon$.

Figure 3.1 shows an cartoon
of what Equation \ref{ag1} looks like geometrically when
$\epsilon=-1$.

\begin{center}
\resizebox{!}{2.5in}{\includegraphics{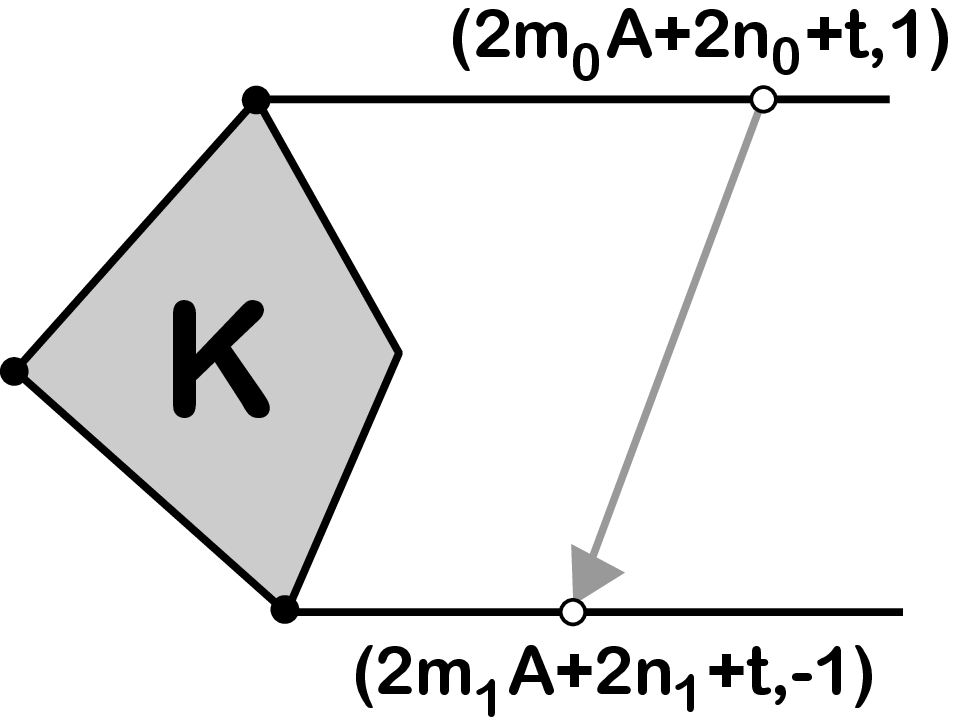}}
\newline
{\bf Figure 3.1:\/} The arithmetic graph construction.
\end{center}

We form a graph $\Gamma_A$ whose vertices are $\Z^2$ by
joining $(m_0,n_0)$ to $(m_1,n_1)$ by an edge
if and only if these points are related as in
Equation \ref{ag1}.  We proved in [{\bf S1\/}] that all the
edges of $\Gamma_A$ have length at most $\sqrt 2$.
That is, there are just $8$ kinds of edges.
We also proved that $\Gamma_A$ is a union of
pairwise disjoint embedded polygonal paths.
All these paths are closed when $pq$ is even.

The nontrivial components of $\Gamma$ all lie in the
half plane above the line $L$ of slope $-A$ through the
origin.  The map
\begin{equation}
f(m,n)=\bigg(2Am+2n+\frac{1}{q},(-1)^{m+n}\bigg)
\end{equation}
carries each component of $\Gamma_A$ to a
different special orbit.  The image of
this map is ``half'' the special special orbits,
in the sense that every special orbit, or its
mirror, is represented.  The {\it mirror\/}
of a special orbit is its reflection in the $x$-axis.
If we wanted to get the other half of the special
orbits, we would use the map
$\rho \circ f$, where $\rho$ is reflection in
the $x$-axis.

It is possible to extend $\Gamma_A$ in a canonical
way so that it fills the entire plane, and not just
the half plane.  We want to do this so that the
Quasi-Isomorphism Theorem is true as stated.
There are two ways to do this, and
they give the same answer.
\begin{enumerate}
\item Dynamically, we can consider the
first return map to the negative-pointing rays
$-\Xi$ and make the same construction.
\item Using the classifying space picture
described below, we can simply take the
domain to be all of $\Z^2$ rather than just
the portion of $\Z^2$ above $L$.
\end{enumerate}
We will take the second approach below.
We call this extended version of $\Gamma_A$ the
{\it arithmetic graph\/} at the parameter $A$.

\subsection{The Canonical Affine Transformation}
\label{CAT}

The Quasi-Isomorphism Theorem compares
the plaid model with a certain affine
image of the arithmetic graph.  In this
section we describe the affine map.
Recall that $A=p/q$ and $\omega=p+q$.  The affine map from
the Quasi-Isomorphism Theorem is given by
$$
T\left(\matrix{x\cr y}\right)=
\frac{1}{A+1}\left(\matrix{A^2+A&A+1 \cr -A^2+2A+1&-2A}\right)
\left(\matrix{x\cr y}\right)+$$
\begin{equation}
\label{k2p}
\left(\matrix{\frac{1}{2q} \cr-\frac{-1}{2q}+\frac{1}{p+q}+\tau}\right).
\end{equation}
Here $\tau$ is the solution in $(0,\omega)$ to the equation
$2p\tau \equiv 1$ mod $\omega$.
The linear part of the map $T$ is defined for irrational
parameters as well as rational parameters, but the map itself
is only defined when $p/q$ is an even rational.
\newline
\newline
{\bf Normalized Arithmetic Graph:\/}
The canonical affine transformation $T$ is
the implied affine map in the Quasi-Isomorphism
Theorem.  We define the {\it normalized arithmetic graph\/}
to be the $T$-image of the arithmetic graph.
The Quasi-Isomorphism Theorem says that the
plaid model and the normalized arithmetic graph are
$2$-quasi-isometric at each parameter.

The inverse of the canonical affine transformation is
$$
T^{-1}\left(\matrix{x\cr y}\right)=
\frac{1}{(A+1)^2}\left(\matrix{2A&A+1\cr
-A^2+2A+1 & -A^2-A}\right) \left(\matrix{x\cr y}\right)+$$
\begin{equation}
\label{p2k}
\frac{1}{2p+2q}\left(\matrix{-1-2q \tau \cr -1+2p\tau}\right).
\end{equation}
Again, the linear part is defined even for irrational
parameters.  However, the entire map is only defined
for even rationals.
\newline
\newline
{\bf The Graph Grid:\/}
We define the {\it graph grid\/} to be the
grid $T(\Z^2)$.  The vertices  of the
normalized arithmetic graph lie in the graph grid.
A calculation shows that $\det(dT)=1+A$, so
the graph grid has co-area $1+A$, as mentioned
in the introduction.  What makes $T$ canonical
is that $T(\Z^2)$ has rotational symmetry
about the origin.
\newline
\newline
{\bf The Anchor Point:\/}
Let
\begin{equation}
\label{anchor}
\zeta=\bigg(\frac{1+A}{2},\frac{1-A}{2}\bigg).
\end{equation}
We compute
$$
T^{-1}(\zeta)=\bigg(\frac{1}{2},\frac{1}{2}\bigg)+
\bigg(\frac{-1-2q\tau}{2\omega},\frac{-1+2p\tau}{2\omega}\bigg).
$$
The second summand is a vector having half-integer
coordinates.  Hence, $T^{-1}(\zeta) \in \Z^2$.
Hence $\zeta \in T(\Z^2)$.
We call $\zeta$ the {\it anchor point\/}.
Given the existence of $\zeta$, we
cen redefine the graph grid 
to be the translate of $dT(\Z^2)$ which
contains $\zeta$. This allows us to define
the graph grid even at irrational parameters,
even though the canonical affine transformation
is only defined for even rational parameters.
As $A$ varies from $0$ to $1$,
the grid $T(\Z^2)$ interpolates between the
grid of half integers and the grid of
integers whose coordinates have odd sum,
and each individual point travels along a
hyperbola or straight line.
\newline
\newline
{\bf Distinguished Edges and Line:\/}
Say that a {\it distinguished edge\/} in the
grid graph is one connecting distinct points of the form
\begin{equation}
T(\zeta), \hskip 40 pt
T(\zeta+(i,j)), \hskip 20 pt
i,j \in \{-1,0,1\}.
\end{equation}
Let ${\cal F\/}(i,j)$ denote the family
of distinguished edges corresponding to
the pair $(i,j)$.  In [{\bf S0\/}] we proved
that the arithmetic graph is embedded, and each
edge is one of the $8$ shortest vectors
in $\Z^2$.  Thus, the edges of the
normalized arithmetic graph are all distinguished.

Say that a {\it distinguished line\/} is a line that contains
a distinguished edge of the grid graph.  If these
lines are parallel to edges in ${\cal F\/}(i,j)$ we
say that the lines have {\it type\/} $(i,j)$.

\subsection{Geometry of the Graph Grid}
\label{gg0}

In this section we prove a number of
statements about the geometry of the graph grid.
The word {\it square\/} always
means integer unit square, as in the plaid
model. 

\begin{lemma}[Grid Geometry]
The following is true at each parameter.
\begin{enumerate}
\item No point of the grid graph lies on the
boundary of a square.
\item Two points of the grid graph cannot lie
in the same square.
\item A union of $3$ horizontally consecutive
squares intersects the grid graph.
\item A union of $2$ vertically consecutive
squares intersects the grid graph.
\item Two graph grid points in consecutive squares
are always connected
by a distinguished edge.
\item If two parallel distinguished lines intersect
the interior of the same edge $e$ of a square, then
the lines have type $(-1,1)$ and $e$ is horizontal.
\item The slopes of the distinguished edges
are never in $\{-1,1,0,\infty\}$.
\end{enumerate}
\end{lemma}

\noindent
{\bf Proof of Statement 1:\/}
We just need to prove, for each
$(x,y) \in \Z^2$, that neither coordinate of
$T(x,y)$ is an integer.  The two coordinates of
$T(x,y)$ are
$$\frac{1}{2q}+\frac{px}{q}+y,
\hskip 20 pt 
\tau-\frac{q-p}{2q(p+q)}+\frac{-2p^2x+4pqx+2q^2x-4pqy}{2q(p+q)}$$
The first number is always $k/(2q)$ where $K$ is an odd integer.
The second number is always $n/(2q(p+q))$, where $n$ is odd.
\endproof

\begin{lemma}
\label{l2norm}
Suppose $v$ is a vector with $\|v\|<\sqrt 2$.
Then $\|dT^{-1}(v)\|<2$.
\end{lemma}

\startproof
For any matrix $M$, we have
$\|M(v)\| \leq \|M_2\|\|v\|.$
Here $\|M\|_2$ is the $L_2$ norm of $M$.
Setting $M=dT^{-1}$ we get
$$\|M\|_2^2=\sum_{ij} M_{ij}^2=
\frac{2(1+3A+4A^2-A^3+A^4)}{(1+A)^4} \leq 2$$
The inequality on the right, which
holds for all $A \in [0,1]$, is
an exercise in calculus:
We check that the derivative of
this expression
does not vanish in $[0,1]$, and then we evaluate at $A=0$ and
$A=1$ to get the bound on the right. 
The final result is that
$\|M\|_2 \leq \sqrt 2.$.  
\endproof

\noindent
{\bf Proof of Statement 2:\/}
Call a vector
in $dT(\Z^2)$ {\it bad\/} if both its coordinates
are less than $1$ in absolute value.
Thanks to the previous result, if this lemma is false
then there is some bad vector in $dT(\Z^2)$.
A bad vector has length less than $\sqrt 2$.
Hence, by the previous result, a bad vector
must have the form 
$dT(\zeta)$, where $\zeta$ is
one of the $8$ shortest nonzero vectors in
$\Z^2$.  
By symmetry we just have
to check $4$ out of the $8$, namely:
\begin{itemize}
\item $dT(1,0)=\bigg(A,\frac{1+2A-A^2}{1+A}\bigg)$.
\item $dT(0,1)=\bigg(1,\frac{-2A}{1+A}\bigg)$.
\item $dT(1,1)=(1+A,1-A)$.
\item $dT(1,-1)=\bigg(-1+A,\frac{1+4A-A^2}{1+A}\bigg).$
\end{itemize}
Here we have set $A=p/q \in (0,1)$.  In all cases, we
can see that least one coordinate is at least $1$
in absolute value.
\endproof

\noindent
{\bf Proof of Statement 3:\/}
It suffices to prove that
every translate of $T^{-1}(R)$ intersects
$\Z^2$.
The vertices of $T^{-1}(R)$ are
$$R_0=(0,0), \hskip 6 pt
R_1=\bigg(1,-\frac{2A}{1+A}\bigg), \hskip 6 pt
R_2=\bigg(3A,\frac{3+6A-3A^2}{1+A}\bigg), \hskip 6 pt
R_3=R_1+R_2.$$
Clearly $R_1$ lies below the line $y=0$. A
bit of calculus shows that $R_2$ and $R_3$ both
lie above the line $y=2$.  Furthermore, any
line of the form $y=h$ with $h \in [0,2]$ intersects
$T^{-1}(R)$ in a segment of width
$$\frac{1+2A+A^2}{1+2A-A^2}>1.$$
This is sufficient to see that every
translate of $T^{-1}(\R)$ intersects $\Z^2$.
\endproof

\noindent
{\bf Proof of Statement 4:\/}
Let $R$ be a union of two vertically consecutive
squares.  Let $C_R$ denote the infinite
column of integer unit squares containing $R$.
Note that $T(0,1)=(1,-P)$.  Say that a
{\it distinguished array\/} is an infinite
set of points of the form $T(m_0,n)$, where
$m_0$ is held fixed and $n \in \Z$.  Because
the first coordinate of $T(0,1)$ is $1$,
we see that every row intersects $C_R$ in one
point.

Let $\rho_{m}$ denote the $m$th row of points.
Let $y_m$ denote the $y$-coordinate of the
point $v_m$ where $\rho_m$ intersects $C_R$.  
Looking at the formulas above, observe that
the difference in the $x$ coordinates
of $dT(1,0)$ and $dT(1,-1)$ is $1$ and that
both $y$ coordinates are in $[1,2]$.
For this reason, one of the two points
$v_m+dT(1,0)$ or
$v_m+dT(1,-1)$ lies on a row above $v_m$
and moreover the vertical distance between
these points is less than $2$. Hence
$|h_{m+1}-h_m|<2$ for all $m$.  But then
$R$ contains some $v_m$.
\endproof

\noindent
{\bf Proof of Statement 5:\/}
Let $$f=v_1-v_2=(f_1,f_2)=dT(i,j).$$  When the
squares are stacked on top of each other,
we have the constraints
$|f_1|<1$ and $|f_2|=2$. When the
squares are stacked on top of each other,
we have the constraints
$|f_1|<2$ and $|f_2|=1$.  Both
cases give $\|f\|<\sqrt 5.$  This
combines with
$\|dT^{-1}\|_2<\sqrt 2$ to show that
$i^2+j^2<10$.
An explicit case-by-case rules everything
out but \begin{itemize}
\item $(i,j) \in \{\pm(1,0),\pm(1,-1)\}$ in
the first case.
\item $(i,j) \in \{\pm(0,1),\pm(1,1)\}$ in
second case.
\end{itemize}
These cases do actually occur.
\endproof

\noindent
{\bf Proof of Statement 6:\/}
Let us first
consider in detail the case when $L_1$ and $L_2$ 
have type $(1,0)$ and the edge $e$ is vertical.
To rule out this case,
we just need to intersect two
adjacent lines with the $y$-axis and see that the distance
between the intersection points is at least $1$ unit.
Let $dT$ be the linear part of $T$.
Two consecutive lines are given by
$$(1-s)dT(0,0)+sdT(1,0), \hskip 30 pt
(1-s)dT(0,1)+sdT(1,1).$$
The first line contains $(0,0)$.  To see where the
second line intersects the $y$-axis, we set the
first coordinate equal to $0$ and solve for $t$.
This gives $t=-1/A$.  Plugging this into the equation,
we see that the $y$ intercept is $-1-1/A$.
Hence, the vertical distance is
$1+1/A$, a quantity that always exceeds $1$.
We record this information by writing
$d(1,0,V)=1+1/A$.  Following the same method, we
do the other $7$ cases. Here is the result:
$$d(1,0,V)=1+\frac{1}{A}, \hskip 30 pt
d(1,0,H)=\frac{1+2A+A^2}{1+2A-A^2}.$$
$$d(0,1,V)=1+A, \hskip 30 pt
d(0,1,H)=\frac{1+2A+A^2}{2A}.$$
$$d(1,1,V)=1, \hskip 30 pt d(1,1,H)=\frac{1+A}{1-A}.$$
$$d(-1,1,V)=\frac{1+A}{1-A}, \hskip 30 pt
d(-1,1,H)=\frac{1+2A+A^2}{1+4A-A^2}.$$
Only the last quantity can drop below $1$.
\endproof

\noindent
{\bf Proof of Statement 7:\/}
We compute that the possible
slopes for the distinguished edges are
$$\frac{1+2A-A^2}{A+A^2}, \hskip 25 pt
\frac{-2A}{1+A}, \hskip 25 pt
\frac{1-A}{1+A}, \hskip 25 pt
\frac{1+4A-A^2}{A^2-1}.$$
It is an easy exercise in algebra
to show that these quantities
avoid the set $\{-1,0,1,\infty\}$ for
all $A \in (0,1)$.
\endproof

\newpage

\section{The Proof in Broad Strokes}

\subsection{Pixellated Squares}

We fix some even rational parameter
$A=p/q$ for the entire chapter.
All the definitions are made with respect
to this parameter.  Also, to save words,
{\it square\/} will always denote a unit
integer square.  Finally, by {\it arithmetic graph\/},
we mean the image of the arithmetic graph under
the canonical affine transformation.
\newline
\newline
{\bf Plaid Nontriviality:\/}
 We say that a square $\Sigma$ is
{\it plaid nontrivial\/} if the plaid
tile at $\Sigma$ is nontrivial. Otherwise we
call $\Sigma$ {\it plaid trivial\/}.
When $\Sigma$ is plaid nontrivial, we define the
{\it plaid edge set\/} of $\Sigma$ to be the two
edges of $\Sigma$ crossed by the
plaid polygon that enters $\Sigma$.
\newline
\newline
{\bf Grid Fullness:\/}
We say that $\Sigma$ is
{\it grid full\/} if $\Sigma$ contains
a point of the graph grid. Otherwise we
call $\Sigma$ {\it grid empty\/}.
When $\Sigma$ is grid full, there is
a unique point of the graph grid
contained in $\Sigma$, and this
point lies in the interior of $\Sigma$.
See the Grid Geometry Lemma.
\newline
\newline
{\bf Graph Nontriviality:\/}
This definition only applies when
$\Sigma$ is grid full.
We say that $\Sigma$ is
{\it graph trivial\/} if 
the point $\sigma$ of the graph grid
contained in $\Sigma$ is isolated in
the arithmetic graph.  Otherwise we
call $\Sigma$ {\it graph nontrivial\/}.
When $\Sigma$ is graph nontrivial, we call
the two edges of the arithmetic graph
incident to $\sigma$
the {\it graph edges associated to\/} $\Sigma$.
\newline
\newline
{\bf Pixellation:\/}
This definition only applies to
squares $\Sigma$ which are graph full.
We call $\Sigma$ {\it pixellated\/} if 
the following is true:
\begin{itemize}
\item $\Sigma$ is graph trivial
if and only if $\Sigma$ is plaid trivial.
\item If $\Sigma$ is plaid nontrivial then
the graph edges associated to $\Sigma$
cross $\Sigma$ in the interiors of the
edges in the plaid edge set of $\Sigma$.
\end{itemize}
When $\Sigma$ is pixellated, the plaid
model at $\Sigma$ determines the local
picture of the graph model in $\Sigma$
in the cleanest possible way. 
\newline

We can see that all the grid full 
squares in Figure 1 are pixellated.
Indeed, for the parameter
$3/8$ and many others, every grid full square
is pixellated.   However, this perfect
situation fails for some parameters.

\subsection{Bad Squares}

\noindent
{\bf Bad Squares\/}
We keep the notation and terminology
from the previous section.
We call the square
$\Sigma$ {\it bad\/} if $\Sigma$
is graph full but not pixellated. 
It turns out that this only happens when
$\Sigma$ is both plaid nontrivial and
graph nontrivial. (See the Pixellation
Theorem below.
\newline
\newline
{\bf Offending Edges:\/}
We say that
an {\it offending edge\/} is a graph
edge associated to $\Sigma$ which does not
cross the boundary of $\Sigma$ in the
interior of one of the edges of the edge set.
\newline
\newline
{\bf Unused Sides:\/}
We say that an
{\it unused side\/} is an edge in the
plaid edge set of $\Sigma$ which is not
crossed by a graph edge associated to
$\Sigma$.  The existence of an
unused side implies the existence of an
offending edge, and {\it vice versa\/}.
\newline

\begin{center}
\resizebox{!}{3.5in}{\includegraphics{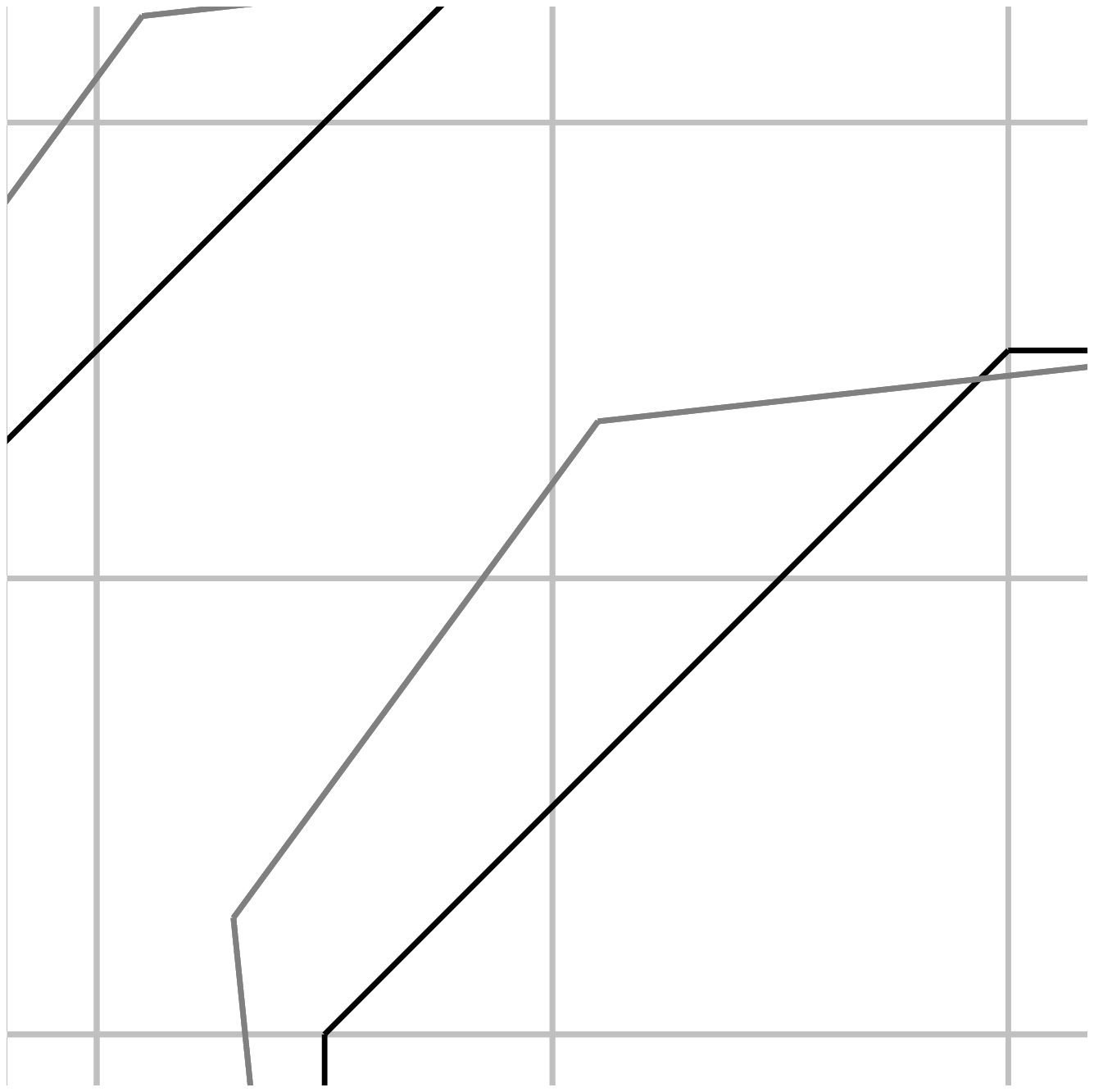}}
\newline
{\bf Figure 4.1:\/} Some bad squares for the parameter $p/q=4/5$.
\end{center}

Figure 4.1 shows a closeup of the picture when $p/q=4/5$.
The plaid polygons are in black and the arithmetic
graph polygons are in grey.
The top right and bottom left squares are bad.  The top left
and bottom right squares are grid empty.  

\noindent
{\bf Catches:\/}
Figure 4.2 shows a picture of two possible patterns of
squares which involve a bad square, an offending
edge and an unused side.  In both cases,
the bad square is meant to be the top right one.
The plaid segment in the top right square either 
connects the south edge to the north edge or to the
east edge.
The shaded squares are meant to be grid empty.
The offending edge is drawn as a curved grey segment.
We mean to consider now just these patterns, but also
the ones which arise by applying a symmetry of the
square grid to these.  In other words, we don't want
to fix the orientation of the pattern in the plane.

In both cases, the portion of the plaid model in the
bottom two squares
makes a straight diagonal segment,
and the sign of the slope of the
offending edge is the same as the sign of the portion
of the diagonal segment.   In these pictures, we say
that the unused side and the offending edge are
{\it associated\/}.
We call these two patterns (and their isometric images)
{\it catches\/} for the offending edge.   

\begin{center}
\resizebox{!}{3in}{\includegraphics{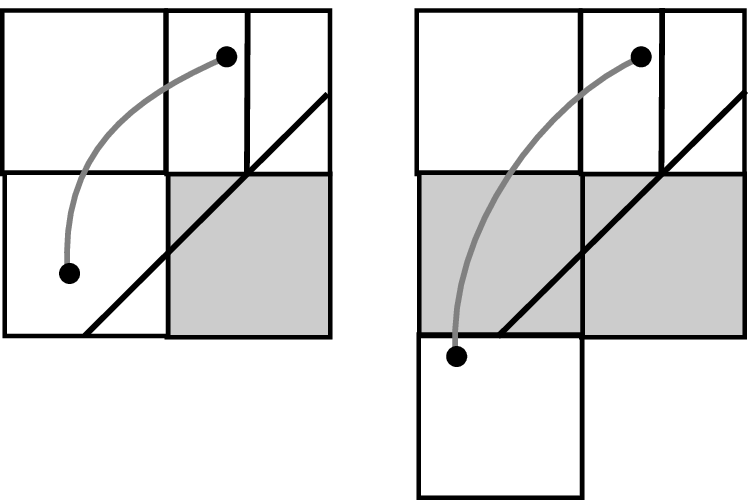}}
\newline
{\bf Figure 4.2:\/} The two catches for an offending edge in a bad square.
\end{center}

The Pixellation Theorem, stated below, says that the pictures
in Figure 4.2 and their rotated/reflected images are always
present when there is an unused side and an offending edge.
We defer the statement of the Pixellation Theorem for a
while, because we want to bundle some other minor results
into it.

\subsection{Errand Edges and Double Crossings}
\label{doublecross}

\noindent
{\bf Errand Edges:\/}
Say that an {\it errant edge\/} is an edge of
the arithmetic graph with sits with respect
to the plaid model as in Figure 4.3.
The grey arc is the arithmetic graph edge.
We have drawn it curved because here and below
a curved picture will look nicer.  Of course,
we are representing a straight line segment.
The numbers $1$ and $2$ denote squares
$\Sigma_1$ and $\Sigma_2$.  The edge
$\Sigma_1 \cap \Sigma_2$ is in the plaid edge set
of $\Sigma_1$ (and of $\Sigma_2$) and hence the
arithmetic graph edge is not an offending
edge with respect to $\Sigma_1$.  However, it rises
up at least one unit above the top of $\Sigma_1$ and $\Sigma_2$
while the plaid polygon in $\Sigma_2$ either goes straight
across $\Sigma_2$ or else moves downward.  As usual,
we mean to consider all possible orientations of
these configurations.  

\begin{center}
\resizebox{!}{1.7in}{\includegraphics{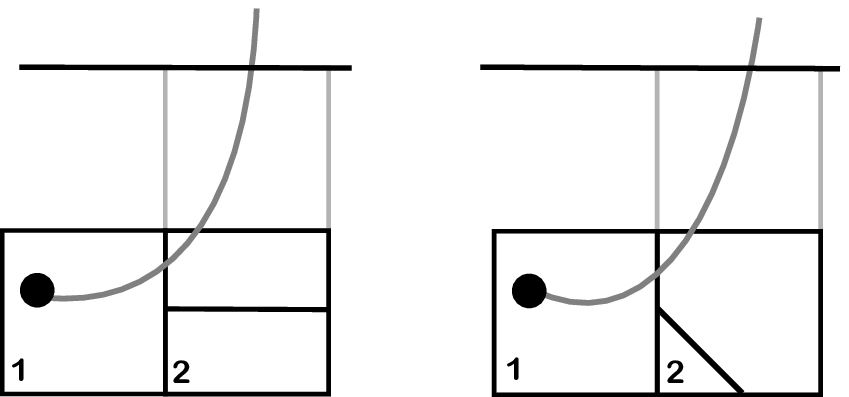}}
\newline
{\bf Figure 4.3:\/} errant edges
\end{center}

\noindent
{\bf Double Crossings:\/}
We say that a {\it double crossing\/} is a union of
two disjoint distinguished edges $e_1$ and 
$e_2$ which have
endpoints $v_1$ and $v_2$ in adjacent 
unit integer squares $Sigma_1$ and $Sigma_2$
and both cross $\Sigma_1 \cap \Sigma_2$ at interior points.  The two
squares may either be stacked on top of each other,
as in Figure 4.4, or stacked side by side.

\begin{center}
\resizebox{!}{1.6in}{\includegraphics{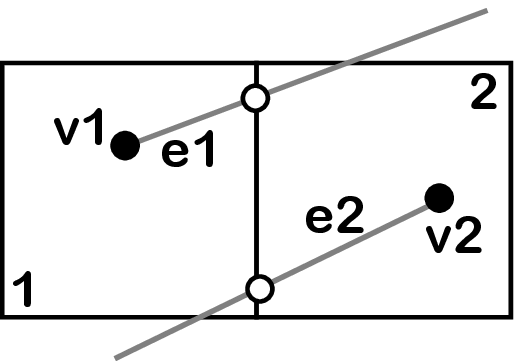}}
\newline
{\bf Figure 4.4:\/} A double crossing
\end{center}

\subsection{The Pixellation Theorem}

Here is our main result.

\begin{theorem}[Pixellation]
The following is true for any even rational parameter.
\begin{enumerate}
\item There are no double crosses in the arithmetic graph.
\item A square is plaid nontrivial if and only if it
is grid nontrivial.
\item There are no errant edges in the arithmetic graph.
\item When a square is graph nontrivial,
the two associated graph edges must
cross distinct sides of the square.
\item In a bad square, there is a bijection
between unused sides and offending edges,
and each matched pair of objects is involved
in a catch.
\end{enumerate}
\end{theorem}

The Pixellation Theorem says that the vast majority of
grid full squares are pixellated, and it precisely
charactarizes the local picture around the ones which
are not.  To deduce the Quasi-Isomorphism Theorem
from the Pixellation Theorem, we need to see that
the various local pictures implied by the Pixellation
Theorem piece together correctly.  The next
section contains the result we need, the Bound Chain Lemma.  Here are
some remarks on the wording in the Pixellation Theorem.

Statement 5 has a complicated phrasing which we want to explain.
Logically, the existence of an unused side implies the
existence of an offending edge, but we want to make sure
that there is an offending egde that specifically is
associated to the unused side of interest to us.  Likewise,
the existence of an offending edge implies the existence
of an unused side, but we want to make sure that there
is an unused side that specifically is associated to the
offending edge of interest to us. 

Statements 1 and 3 will not be used directly in our
deduction of the Quasi-Isomorphism Theorem.  However,
they will be used, in the next chapter, for the deduction
of the Bound Chain Lemma.

\subsection{The Bound Chain Lemma}

The main weakness in the Pixellation Theorem
is that it only deals with grid full squares.
Here we deal with the grid empty squares.

We say that a finite union of squares
$\Sigma_1,...,\Sigma_m$ is {\it linked\/} if
\begin{itemize}
\item $\Sigma_1$ and $\Sigma_m$ are grid fill, and
the remaining squares are grid empty.
\item $\Sigma_k \cap \Sigma_{k+1}$ is an edge, for
each $k=1,...,(m-1)$.
\item A single plaid
polygon intersects $\Sigma_k$ for all $k$.
\end{itemize}
To avoid trivialities, we take $m \geq 2$.
If follows from Statements 3 and 4 of
the Grid Geometry Lemma
that a linked chain has
length at most $4$.

We say that our linked sequence
is {\it bound\/} if a single edge in the
arithmetic graph joins the graph grid point
in $\Sigma_1$ to the graph grid point in
$\Sigma_m$.

Figure 4.5 shows a linked and bound chain of
length $4$.  We have drawn the arithmetic
graph edge as curved, to get a nicer picture.

\begin{center}
\resizebox{!}{3.4in}{\includegraphics{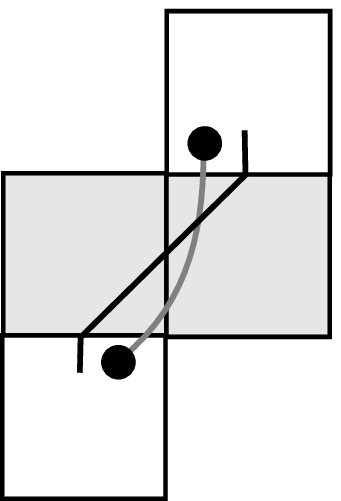}}
\newline
{\bf Figure 4.5:\/} A Linked and Bound Chain of Length $4$
\end{center}

In the next chapter, we will prove the
following result, essentially using the
Pixellation Theorem and a case-by-case
analysis.  

\begin{lemma}[Bound Chain]
Every linked chain is bound.
\end{lemma}

\subsection{Proof of the Quasi-Isomorphism Theorem}

The Quasi-Isomorphism Theorem is a fairly immediate
consequence of the Pixellation Theorem and the
Bound Chain Lemma.  Here we
produce an explicit homeomorphism between the
union of polygons in the plaid model and the union
of polygons in the arithmetic graph.  The
homeomorphism moves no point more than $2$ units.
The homeomorphism is defined in several pieces,
depending on the type of square we have.

Let $\Sigma$ be a pixellated square.  Let $v$ be the
graph grid point in $\Sigma$ and let 
$e_1$ and $e_2$ be the two edges incident to 
$v$.  Let $v'$ be the midpoint of the plaid
segment contained in $\Sigma$ and let
$e_1'$ and $e_2'$ be the halves of the
plaid segment on either side of $v'$.
We map $e_k \cap \Sigma$ linearly to $e_k'$
for $k=1,2$.  We choose the labeling so that
$e_k$ and $e_k'$ both intersect the same
side of $\Sigma$. Figure 4.6 shows this simple
map in action.

\begin{center}
\resizebox{!}{2.5in}{\includegraphics{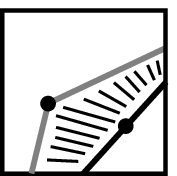}}
\newline
{\bf Figure 4.6:\/} The Map on Pixellated Squares
\end{center}

We use the same notation for discussing bad
squares.  If $\Sigma$ is a bad square and $e_1$ is
not an offending edge of $\Sigma$ then we do the
same map as for pixellated squares.  If $e_1$ is
an offending edge, we map $e_1$ to the union
of squares in the catch as shown in Figure 4.7
for each of the two kinds of catches.
The white dots in Figure 4.7 are just extra
guides for the map.  On the right hand side of
Figure 4.7, we have drawn one particular way that
the plaid segment could look.  The black dot
the bottom square is meant to be the midpoint
of the plaid segment, whatever it looks like.

\begin{center}
\resizebox{!}{3.5in}{\includegraphics{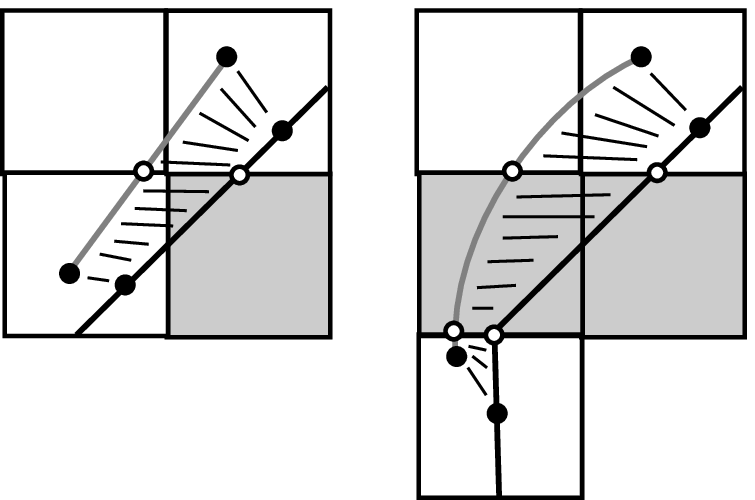}}
\newline
{\bf Figure 4.7:\/} The Map on Bad Squares
\end{center}

Say that a {\it clean linked chain\/} is a chain
$\Sigma_1,...,\Sigma_m$ such that the arithmetic
graph edge $e$ connecting $v_1$ to $v_m$ crosses
$\Sigma_1 \cap \Sigma_2$ and $\Sigma_{m-1} \cap \Sigma_m$.
We say that the {\it graph core\/} associated
to the clean linked chain is the segment 
\begin{equation}
\widehat e=e-(\Sigma_1 \cap \Sigma_m).
\end{equation}
The segment $\widehat e$ has endpoints on
$\Sigma_1 \cap \Sigma_2$ and
$\Sigma_{m-1} \cap \Sigma_m$.  
Corresponding to the $\widehat e$
is the portion of the plaid model connecting
these same two edges.  We call this the
{\it plaid core\/} associated 
the clean linked chain.

Our map is defined everywhere on the
arithmetic graph polygons except on the
graph cores.  These comprise a disjoint
set of segments.  The image of our map
so far is exactly the complement of the
plaid cores. We finish the proof by mapping
the graph cores to the plaid cores according
to the scheme in Figure 4.8.

Figure 4.8 doesn't show every possibility,
but these examples should be sufficient to 
show what we do in every case.
In every case, the Pixellation Theorem and the
Bound Chain Theorem imply that the map is 
well defined and only moves points by at most
$2$ units.  We call our map $\Theta$.
By construction,
$\Theta$ maps each arithmetic
graph polygon homeomorphically to some plaid polygon.

\begin{center}
\resizebox{!}{3in}{\includegraphics{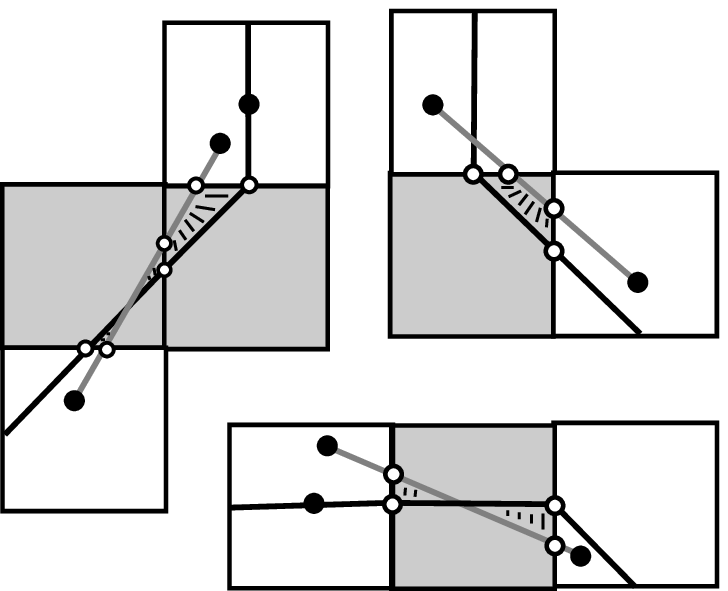}}
\newline
{\bf Figure 4.8:\/} The Map on Bad Squares
\end{center}

\begin{lemma}
$\Theta$ is injective.
\end{lemma}

\startproof
Suppose $\gamma_1$ and $\gamma_2$ are
graph polygons which both map to the
plaid polygon $\pi$.   Suppose 
$\pi$ never enters a square with a 
non-offending edge.  Then, according
to the Pixellation Lemma, $\pi$ travels
in a straight diagonal line.  Since
$\pi$ is closed, this situation is
impossible. Hence,
$\pi$ must enter at least one square $\Sigma$
which has at least one non-offending edge $e$.
The only way $\gamma_1$ and $\gamma_2$ are
both mapped on to $e \cap \Sigma$ is if both
these polygons contain the grid point in
$\Sigma$.  But then $\gamma_1=\gamma_2$.
\endproof

\begin{lemma}
$\Theta$ is surjective.
\end{lemma}

\startproof
Each plaid polygon $\pi$ has length at least $4$.
Hence, by Statements 3 and 4 of the
Grid Geometry Lemma
the polygon $\pi$ intersects at least one
grid full square $\Sigma$.  But then $\Phi$ maps
 the arithmetic
graph polygon which contains the graph grid
point in $\Sigma$ to $\pi$.
\endproof

In short, $\Theta$ is a bijection between the
components and a homeomorphism on each one,
and $\Theta$ moves points by at most $2$ units.
This completes the proof of the Quasi-Isomorphism
Theorem.

\newpage

\section{Proof of the Bound Chain Lemma}

\subsection{Length Two Chains}
\label{chain2}

Let $\Sigma_1,\Sigma_2$ be a linked chain
of length $2$.  This means that
both $\Sigma_1$ and $\Sigma_2$ are plaid
nontrivial and grid full, and a single plaid
polygon runs through both. 
By the Pixellation Theorem, $\Sigma_1$ and $\Sigma_2$ are
graph full as well.

\begin{center}
\resizebox{!}{3in}{\includegraphics{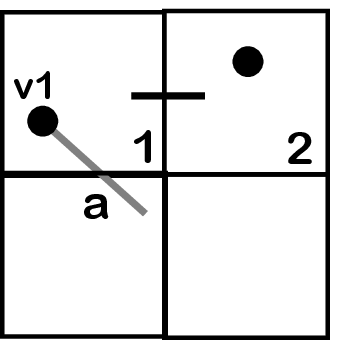}}
\newline
{\bf Figure 5.1:\/} $\Sigma_1$ and $\Sigma_2$, and an uncaught
offending edge.
\end{center}

Let $v_j$ be the graph grid point in $\Sigma_j$.
Suppose that no graph edge incident to $v_1$
crosses $\Sigma_1 \cap \Sigma_2$. Then there
must be an offending edge $a$ incident to $v_1$
and associated to $\Sigma_1 \cap \Sigma_2$.
This edge must cross either the bottom or the
top of $\Sigma_1$, and we have shown it
crossing the bottom.  But the picture
contradicts Statement 5 of
the Pixellation Theorem, because the
catch for $a$ would involve $\Sigma_2$, and
$\Sigma_2$ would be grid empty in the catch.

The argument in the preceding paragraph shows that
there is some graph edge $e_1$ incident to
$v_1$ which crosses $\Sigma_1 \cap \Sigma_2$.
Likewise, there is some graph edge $e_2$ incident
to $v_2 \in \Sigma_2$ which crosses $\Sigma_1 \cap \Sigma_2$.
If $e_1$ and $e_2$ are disjoint, then we have a
bad configuration.  The Bad Configuration Lemma
rules this out.  Since the arithmetic graph is embedded,
one of our two edges connects $v_1$ to $v_2$, and
in fact these two edges must coincide. (Otherwise
we contradict the Pixellation Theorem.)  In short,
$\Sigma_1,\Sigma_2$ is a bound chain.

\subsection{Length Three Chains: Case A}

Here we consider a length $3$ chain
in which the three squares are the same
horizontal row.

\begin{center}
\resizebox{!}{1.8in}{\includegraphics{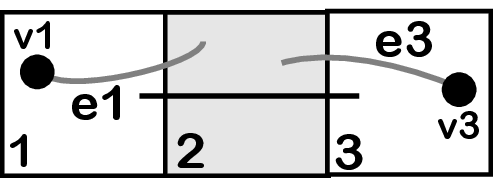}}
\newline
{\bf Figure 5.2\/} A horizontal chain
\end{center}

\begin{lemma}
\label{horizontal}
There is an edge of the arithmetic graph connecting
$v_1$ and $v_3$.
\end{lemma}

Because the plaid polygon in $\Sigma_2$ is not a diagonal
segment, there is no way to involve $\Sigma_2$ in a catch
for an offending edge incident to $v_1$ and
associated to $\Sigma_1 \cap \Sigma_2$. Hence,
$\Sigma_1 \cap \Sigma_2$ is not an unused edge of $\Sigma_1$.
Hence there is some graph edge 
$e_1$ incident to $v_1$ and crossing $\Sigma_1 \cap \Sigma_2$.
The same argument shows that there is an arithmetic graph
edge $e_3$ associated to $\Sigma_3$ which crosses $\Sigma_2 \cap \Sigma_3$.

If the other endpoint $e_1$ of $e_1$ lies in $\Sigma_3$ then
this other endpoint must be $v_3$, because there is at most
one graph grid point in $\Sigma_3$.  
So, to finish the proof, we just have to rule out the other
possible squares which could potentially contain the
endpoint of $e_1$.  

The edge $e_1$ has length less than $2$ because $\|dT\|_2<\sqrt 2$
and $e_1$ is the image of a vector in $\Z^2$ of length at
most $\sqrt 2$. This bound on the length cuts down on the
possibilities where $e_1$ could end up.
For one thing, $e_1$ cannot end up in any square
to the right of $\Sigma_3$.  Also, $e_1$ cannot cross the line
one unit above the top of $\Sigma_2$, by Statement 3 of
the Pixellation Theorem: no errant edges.
EHence, the other endpoint of
$e_1$ lies in the row of squares above $\Sigma_2$.
We will consider the
$3$ cases when $e_1$ ends in a square in the row above
$\Sigma_3$.  The ``below'' case has the same treatment.
\newpage

\noindent
{\bf Case 1:\/}
Suppose first $e_1$ ends in $\Sigma_4$, as
shown on the left in Figure 5.3.
In this case, $e_1$ connects $v_1$ with the
graph grid point $v_4$ in $\Sigma_4$.   The plaid
edge set of $\Sigma_4$ cannot contain
$\Sigma_2 \cap \Sigma_4$ because the $\Sigma_2 \cap \Sigma_4$
is not in the plaid edge set of $\Sigma_2$.
Hence $e_1$ is an offending edge for $\Sigma_4$.  But
then, by the Pixellation Theorem, the portion of
the plaid model inside $\Sigma_1$ must look as drawn.
This is impossible, because it forces the plaid
edge set of $\Sigma_1$ to have $3$ edges in it.

\begin{center}
\resizebox{!}{1.6in}{\includegraphics{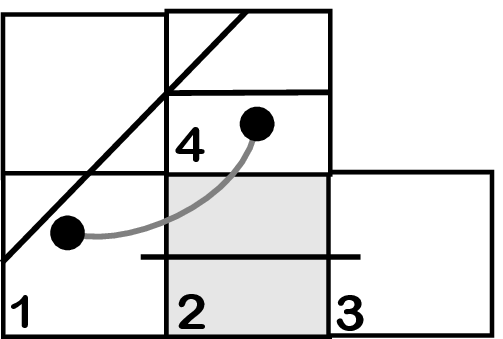}}
\newline
{\bf Figure 5.3:\/} Case 1.
\end{center}

\noindent
{\bf Case 2:\/}
Suppose $e_1$ crosses $\Sigma_3 \cap \Sigma_5$ and
ends in $\Sigma_5$.
Note that $e_1$ cannot be an offending edge for
$\Sigma_5$ because $\Sigma_2$ does not have the
right form to be part of a catch for $e_1$. Hence
$\Sigma_3 \cap \Sigma_5$ is in the plaid edge set
for $\Sigma_5$.  But then the result for length $2$
chains says that some some graph edge joins
$v_3$ to $v_5$.  But then both graph edges incident
to $v_5$ cross the same edge of $\Sigma_5$.
This contradicts Statement 4 ofthe Pixellation Lemma.

\begin{center}
\resizebox{!}{1.6in}{\includegraphics{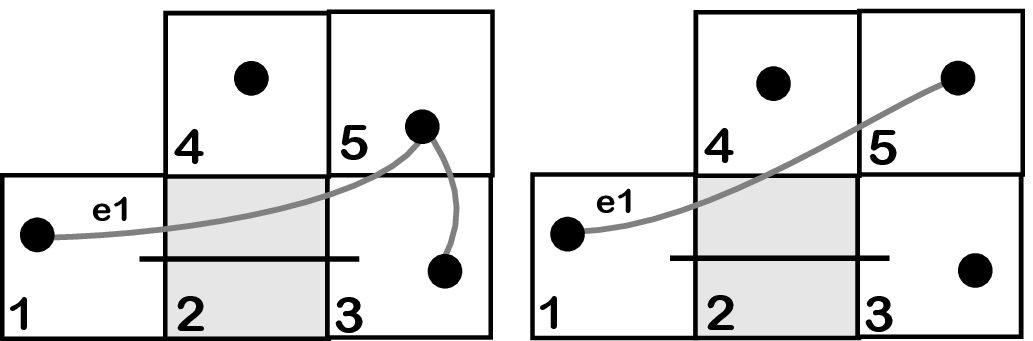}}
\newline
{\bf Figure 5.4:\/} Cases 2 and 3
\end{center}

\noindent
{\bf Case 3:\/}
Suppose that $w_1 \in \Sigma_5$ and $e_1$ crosses
$\Sigma_4 \cap \Sigma_5$.  By Statement 4 of
the Grid Geometry Lemma, the square $\Sigma_4$ is
grid full, just like $\Sigma_3$ is.
The same argument as in Case 2, with
$\Sigma_4$ replacing $\Sigma_3$, takes
care of this case.
\endproof

\subsection{Length Three Chains: Case B}

Here we consider a length $3$ chain
in which the three squares are the same
vertical column.  Our proof refers to
Figure 5.5

This situation
is not quite the same as in Case A
because we can have $2$ horizontally consecutive
grid empty squares.  If neither $\Sigma_4$ nor
$\Sigma_5$, shown in Figure 5.5, is grid
empty, then the same argument as in Case A
applies here.  We just have to
worry about the case when one of
$\Sigma_4$ or $\Sigma_5$ is grid empty.
We will consider the case when $\Sigma_5$
is grid empty.  The other case has the
same treatment.

\begin{center}
\resizebox{!}{3.5in}{\includegraphics{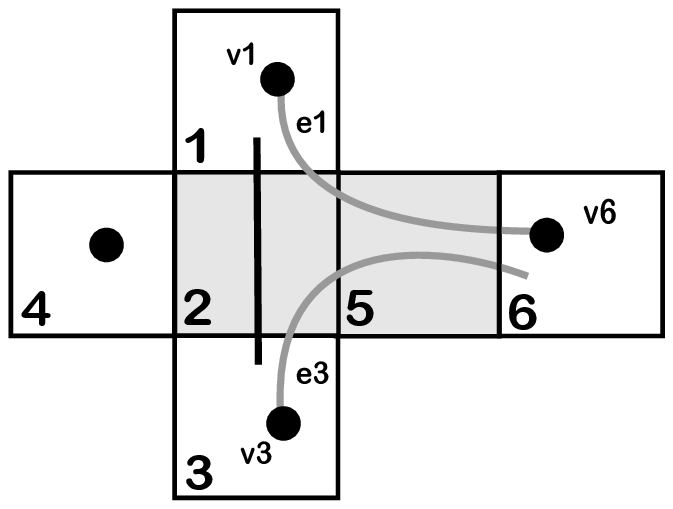}}
\newline
{\bf Figure 5.5:\/} A vertical chain 
\end{center}

In order to avoid finishing the proof as
in Lemma \ref{horizontal} both edges
$e_1$ and $e_3$ must cross
$\Sigma_2 \cap \Sigma_5$.
Neither edge can end in $\Sigma_5$ because
$\Sigma_5$ is grid empty.
If $e_1$ crosses the bottom edge of $\Sigma_5$ then
$e_1$ blocks $e_3$ from exiting $\Sigma_5$, which
is a contradiction.  Hence $e_1$ does not
exit the bottom edge of $\Sigma_5$.
Similarly, $e_3$ cannot cross the top edge of
$\Sigma_5$.  Hence, both $e_1$ and $e_3$ cross
$\Sigma_5 \cap \Sigma_6$.  But then
$e_1$ and $e_3$ are errant edges, and this
contradicts Statement 3 of the Pixellation Theorem.

\subsection{Length Three Chains: Case C}

Now we consider the remaining kind of length $3$ chain,
shown in Figure 5.6.
As usual, let $v_j$ be the graph grid point
in $\Sigma_j$ for $j=1,3$.  We want to prove
that an arithmetic graph edge connects
$v_1$ to $v_3$.  This case is rather painful.

\begin{center}
\resizebox{!}{2.8in}{\includegraphics{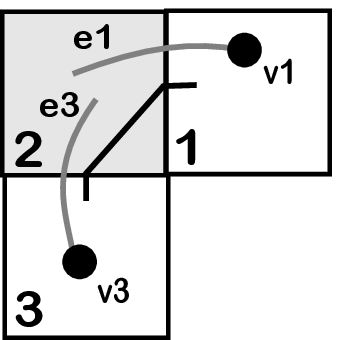}}
\newline
{\bf Figure 5.6\/} 
\end{center}

By Statement 5 of the Pixellation Theorm,
$\Sigma_1 \cap \Sigma_2$ cannot be an
unused edge with respect to $\Sigma_1$,
for the following reasons.
\begin{itemize}
\item If the associated offending edge were
to cross through the top of $\Sigma_1$, then
the plaid segment in $\Sigma_2$ (which is
part of the catch) would have
the wrong position.

\item If the associated offending edge were
to cross through the bottom of $\Sigma_1$, then
the plaid segment in $\Sigma_3$ (which is
part of the catch) would have the wrong position.
\end{itemize}

Hence some graph edge $e_1$ incident to
$v_1$ crosses $\Sigma_1 \cap \Sigma_2$.
The same argument applies to $\Sigma_3$. So, we can
assume that the edges $e_1$ and $e_3$ are as in Figure 5.6.
We will show that $e_1$ has its other endpoint in
$\Sigma_3$.
Since there is only one graph grid point
in $\Sigma_3$, the other endpoint of $e_1$ must be $v_3$.
In short, $e_1$ is the desired arithmetic graph edge
connecting $v_1$ and $v_3$. There are $6$ situations
we must rule out, and we deal with them in turn.
\newline
\newline
{\bf Case 1:\/}
Suppose $e_1$ ends in $\Sigma_4$.
Since $\Sigma_2 \cap \Sigma_4$ is not in the
edge set of $\Sigma_4$, the edge $e_1$ is offending
with respect to $\Sigma_4$.
Given the negative slope of $e_1$, the
associated unused edge must be the right
edge of $\Sigma_4$ and the catch must involve
$\Sigma_1,\Sigma_2,\Sigma_4$.  But then
the plaid segment in $\Sigma_1$ is in the
wrong position.  This is a contradiction.
The long $Y$-shaped graph in Figure 5.7 shows
the shape of the plaid arc implied by the
existence of the catch, and this contradicts the
fact that the plaid arc in $\Sigma_1$ also
crosses into $\Sigma_2$.

\begin{center}
\resizebox{!}{3in}{\includegraphics{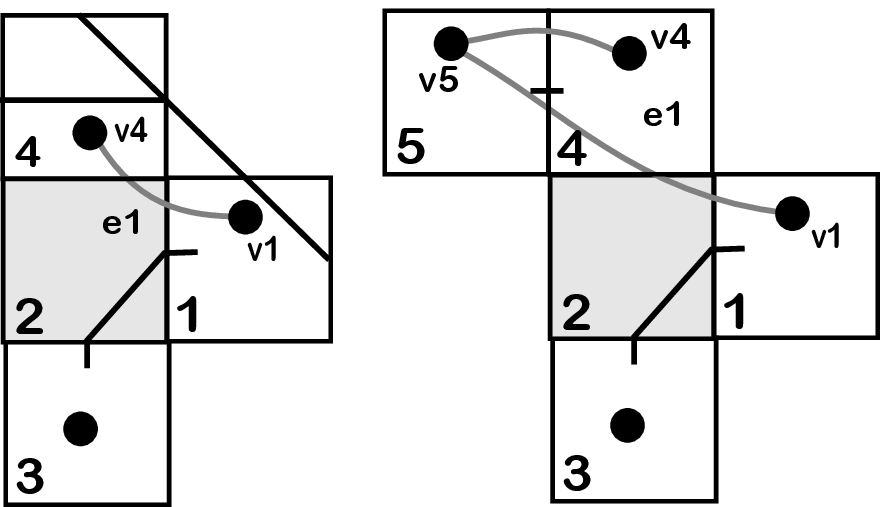}}
\newline
{\bf Figure 5.7:\/}  Cases 1 and 2
\end{center}

{\bf Case 2:\/}
Suppose $e_1$ ends in $\Sigma_5$ and that
$e_1$ crosses $\Sigma_4 \cap \Sigma_5$.
By Statement 4 of the Grid Geometry Lemma, the square
$\Sigma_4$ is grid full.  So, $\Sigma_4$ and
$\Sigma_5$ are both grid full.  If $e_1$ is an
offending edge with respect to $\Sigma_5$, then
the associated unused edge is the bottom edge
of $\Sigma_5$.  But then
the catch for $e_1$ must involve 
$\Sigma_2,\Sigma_4,\Sigma_5$, and the
plaid segment in $\Sigma_2$ is in the wrong position.
Hence $\Sigma_4 \cap \Sigma_5$ is in the
edge set for $\Sigma_5$.  Hence
$\Sigma_4,\Sigma_5$ form a linked chain and
some arithmetic graph edge joins $v_4$ and $v_5$.
This contradicts the Statement 4 of the
Pixellation Theorem.
\newpage

\noindent
{\bf Case 3:\/}
Suppose that $e_1$ ends in $\Sigma_5$ and crosses
$\Sigma_5 \cap \Sigma_6$, as shown on the left in 
Figure 5.8. 
If $\Sigma_6$ is not grid empty, then
the same argument in Case 2 finishes the job.
So, we may assume that $\Sigma_6$ is grid empty.

Now, $e_1$ cannot be an offending edge with
respect to $\Sigma_5$ for the same reason as
in Case 2.  Thus, $\Sigma_5 \cap \Sigma_6$
is in the plaid edge set for $\Sigma_6$.  The
left edge of $\Sigma_6$ cannot be in the 
plaid edge set for $\Sigma_6$ because this
would make $e_1$ an errant edge for $\Sigma_5$.
Hence, the segment of the plaid model in
$\Sigma_6$ connects $\Sigma_5 \cap \Sigma_6$
to $\Sigma_6 \cap \Sigma_7$.  But then,
by Case B in the previous section, some
arithmetic graph edge connects
$\Sigma_5$ to $\Sigma_7$, as shown in Figure 5.10.
But then two arithmetic graph edges cross
$\Sigma_5 \cap \Sigma_6$, contradicting 
Statement 4 of the Pixellation Theorem.

\begin{center}
\resizebox{!}{2.5in}{\includegraphics{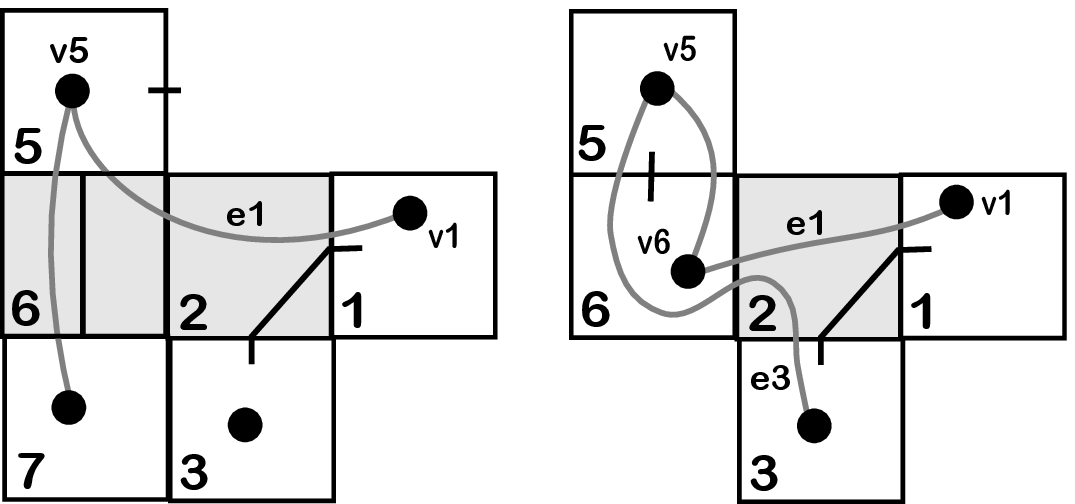}}
\newline
{\bf Figure 5.8:\/}  Cases 3 and 4
\end{center}

\noindent
{\bf Case 4:\/}
Suppose that $e_1$ ends in $\Sigma_6$ and $e_3$ crosses $\Sigma_2 \cap \Sigma_6$.
If $e_3$ ends at
$v_6$, then both $e_1$ and $e_3$ cross the same edge
of $\Sigma_6$ and we contradict Statement 4 of
 the Pixellation Theorem.
Since $e_3$ has length less
than $2$ and is not an errant edge for $\Sigma_3$, the
only possibility is that $e_3$ ends at $v_5$.  But
the same argument as in Case 3 shows that $e_3$ is not an
offending edge for $\Sigma_5$.  Hence $\Sigma_6,\Sigma_5$
is a linked chain. But then some other edge of the
arithmetic graph connects $v_6$ to $v_5$.  This contradicts
Statement 4 of the Pixellation Theorem.

\newpage
\noindent
{\bf Case 5:\/}
Suppose that $e_1$ ends in $\Sigma_6$ and
$e_3$ does not cross
$\Sigma_6 \cap \Sigma_2$.  Note that
$e_1$ blocks $e_3$ from crossing the top
of $\Sigma_2$, so $e_2$ crosses $\Sigma_1 \cap \Sigma_2$.

Note that $e_1$ is an offending edge
for $\Sigma_6$, because $\Sigma_2$ does not
have $\Sigma_6 \cap \Sigma_1$ in its plaid edge set.
The plaid segment in $\Sigma_2$ is not in the
correct position for $\Sigma_2$ to be part of 
a catch for $e_1$, as would happen if $e_1$ had
negative slope.  Hence $e_1$ has positive slope.
Note that $e_3$ also has
positive slope, because it rises up
to cross $\Sigma_3 \cap \Sigma_2$.  So, both
$e_1$ and $e_2$ have the same slope.
(In the reflected case, both would have negative
slope.)

There are two kinds of distingished edges
having positive slope, and only one of them
has a horizontal projection of length greater
than $1$. (The same goes in the negative slope
case, which arises in the reflected case.)  
We know that the horizontal projection
of $e_1$ is greater than $1$. If the horizontal
projection of $e_3$ is greater than $1$, then
$e_1$ and $e_3$ must be parallel.  However, then
we contradict Statement 6 of the Grid Geometry Lemma.
We have found
two parallel distinguished lines which intersect
the same vertical edge of a square.  In short
$e_3$ has horizontal projection at most $1$.
Hence the endpoint of $e_3$ is either in
$\Sigma_1$ or in $\Sigma_7$.

If the endpoint of $e_3$ lies in $\Sigma_1$, we are
done.  So, consider the case when the endpoint
of $e_3$ lies in $\Sigma_7$, as shown in Figure 5.9.
In this case, we get the same contradiction as in Case 4.

\begin{center}
\resizebox{!}{3.3in}{\includegraphics{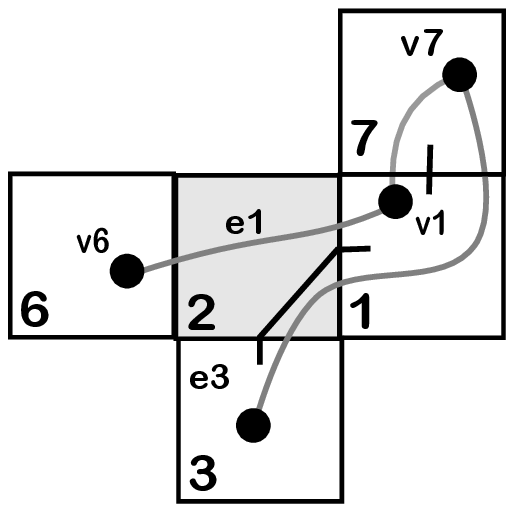}}
\newline
{\bf Figure 5.9:\/}  Case 5
\end{center}

\noindent
{\bf Case 6:\/}
Suppose that $w_1 \in \Sigma_8$ or $w_1 \in \Sigma_9$.
Then $e_1$ has positive slope and horizontal
projection at least $1$.  If $e_3$ does
not cross $\Sigma_2 \cap \Sigma_6$ then we
get the same contradiction as in Case 5.
So, we can assume that $e_3$ crosses
$\Sigma_2 \cap \Sigma_6$, as shown in Figure 5.10.
But then, as Figure 5.10 indicates, $e_3$ cannot
be an offending edge for $\Sigma_6$ because
the plaid segment in $\Sigma_3$ is not in
the right position. Hence, $\Sigma_2 \cap \Sigma_6$
is in the edge set for $\Sigma_6$. But this would
make the edge set of $\Sigma_2$ have $3$ members,
which is a contradiction.

\begin{center}
\resizebox{!}{4.5in}{\includegraphics{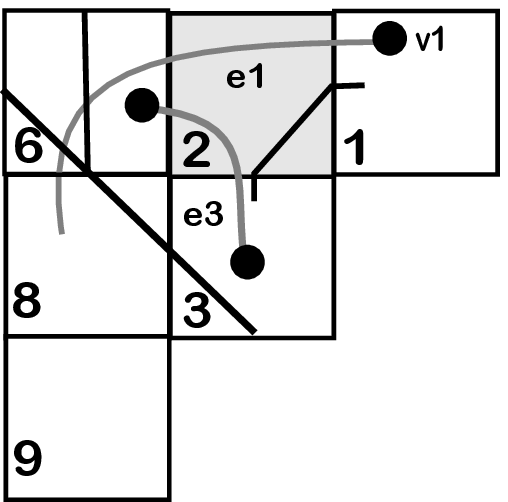}}
\newline
{\bf Figure 5.10:\/}  Case 6
\end{center}

This completes the analysis of Case C.
At this point we have proved the
Bound Chain Lemma for all linked chains, 
except those of length $4$.  
\newpage

\subsection{Length Four Chains: Case A}

Suppose that $\Sigma_1,\Sigma_2,\Sigma_3,\Sigma_4$ is
a length $4$ linked chain, and all four of these

squares are on the same row. 
The same argument as in Case 3A above shows that
there must be some graph edge $e_1$, incident to
$v_1$, which crosses $\Sigma_1 \cap \Sigma_2$.
Figure 5.11 shows a tree of possibilities.
We will show that, actually, this case cannot occur.
Without loss of generality, we will consider
the case when $e_1$ does not end up in a
row of squares below our chain.
Since $e_1$ has length less than $2$, and $e_1$
is not an errant edge, $e_1$ must end
either in $\Sigma_5$ or $\Sigma_6$.

\begin{center}
\resizebox{!}{2.5in}{\includegraphics{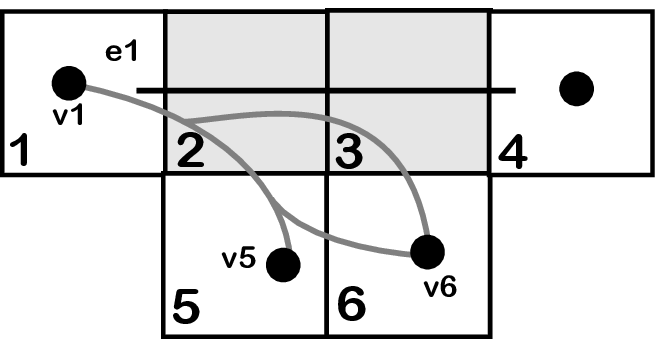}}
\newline
{\bf Figure 5.11:\/}  Four-in-a-row case
\end{center}

\noindent
{\bf Case 1:\/}
Suppose $e_1$ lands in $\Sigma_5$. Since
$\Sigma_2 \cap \Sigma_5$ is an unused edge for
$\Sigma_5$ and $e_1$ crosses this edge, $e_1$
must be an offending edge associated to the
left edg of $\Sigma_5$.  But then
$\Sigma_1,\Sigma_2,\Sigma_5$ are part of the
catch for this edge.  The plaid segment in
$\Sigma_1$ is in the wrong position for this.
\newline
\newline
{\bf Case 2:\/} Suppose
$e_1$ lands in $\Sigma_6$ and crosses
$\Sigma_5 \cap \Sigma_6$. 
Then $e_1$ cannot be an offending edge
for $\Sigma_6$ because the plaid
segment in $\Sigma_3$ is in the wrong position
for it to participate in the required catch.
Hence, some other arithmetic graph edge
joins $v_5$ and $v_6$.  But then
two graph edges are incident to $v_6$
and cross $\Sigma_5 \cap \Sigma_6$. This
contradicts Statement 4 of the Pixellation Theorem.
\newline
\newline
\noindent
{\bf Case 3:\/}
Suppose $e_1$ lands in $\Sigma_6$ and
crosses $\Sigma_3 \cap \Sigma_6$. The same
kind of argument as in Case 1 works here.

\subsection{Length Four Chains: Case B}

Suppose that $\Sigma_1,\Sigma_2,\Sigma_3,\Sigma_4$ is
a length $4$ linked chain, and exactly $3$ of these
squares are on the same row.  Without loss of
generality, we consider the case shown in Figure 5.12.
The same argument as in Case A above shows that
there must be some graph edge $e_1$, incident to
$v_1$, which crosses $\Sigma_1 \cap \Sigma_2$.

Since $e_1$ has length less than $2$ and cannot
be an errant edge, $e_1$ must end in one
of the squares $\Sigma_4,\Sigma_5,\Sigma_6,\Sigma_7$.
Moreover, all these squares are grid full, by
Statement 4 of the Grid Geometry Lemma.  We want to show
that $e_1$ ends in $\Sigma_4$. There are $4$ cases
to rule out. All $4$ cases are handled by
arguments just like those in Cases 4A1 and 4A2 above.

\begin{center}
\resizebox{!}{4.5in}{\includegraphics{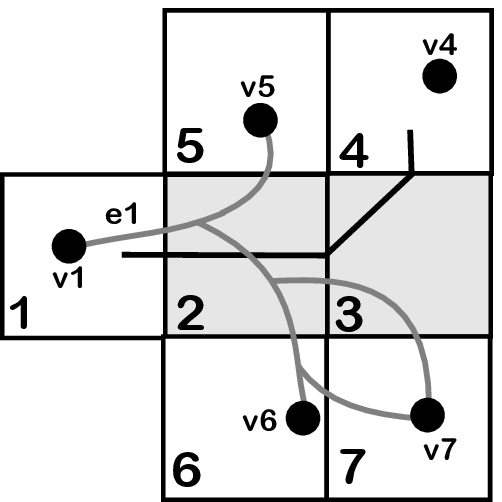}}
\newline
{\bf Figure 5.12:\/}  Three-in-a-row case.
\end{center}

\noindent
{\bf Remark:\/}
We have stopped short of ruling out the existence
of chains like this, but we think that they never
actually occur.

\subsection{Length Four Chains: Case C}

Suppose that $\Sigma_1,\Sigma_2,\Sigma_3,\Sigma_4$ is
a length $4$ linked chain, making a $2 \times 2$ block
as shown in 
Figure 5.13. 

The same argument as in previous
cases shows that there is some graph edge $e_1$
incident to $v_1$ which crosses $\Sigma_1 \cap \Sigma_2$.
We want to show that $e_1$ connects $v_1$ to $v_4$.
Actually, this situation is impossible, given that we
already know that $e_1$ crosses $\Sigma_1 \cap \Sigma_2$.
So, our argument will really show that this kind
of linked chain is impossible.

Before we start our analysis, we make special mention of
the squares $\Sigma_6$ and $\Sigma_9$.  
We do not know in general whether thess squares are
grid full or grid empty.  However, these squares are
irrelevant for all our arguments unless $e_1$ actually
ends in them.  In those cases, the relevant square
is grid full, by definition.  
So, in all relevant cases, $\Sigma_6$ and $\Sigma_9$
are grid full, as drawn.

\begin{center}
\resizebox{!}{3.2in}{\includegraphics{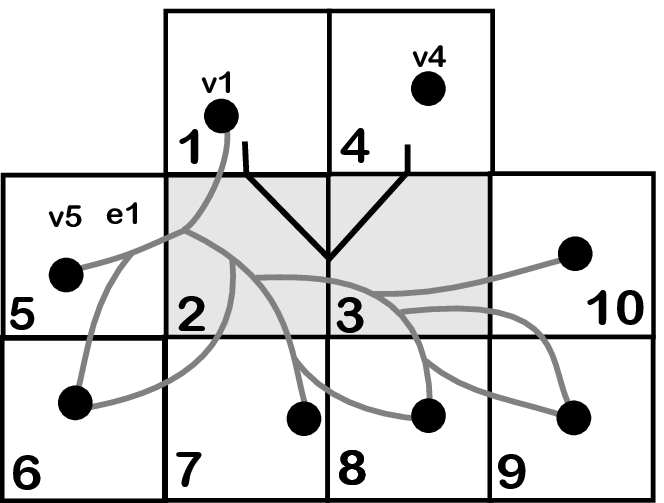}}
\newline
{\bf Figure 5.13:\/}  Block case
\end{center}

\noindent
{\bf Case 1:\/} Suppose $e_1$ ends in $\Sigma_5$
or $\Sigma_10$.
This is the same as Case 4A1.
\newline
\newline
{\bf Case 2:\/} Suppose $e_1$ lands in $\Sigma_6$,
$\Sigma_8$, or $\Sigma_9$.  In all these cases,
the argument is the same as in Case 4A2.
\newline
\newline
{\bf Case 3:\/}
Suppose that $e_1$ lands in $\Sigma_7$.   Here we 
redraw the picture to focus more particularly
on this case. 
Since $\Sigma_7 \cap \Sigma_2$ is not in
the edge set for $\Sigma_2$, the edge
$e_1$ must be an offending edge in $\Sigma_7$
associated to $\Sigma_7 \cap \Sigma_2$.

If $e_1$ has negative slope, then
$\Sigma_2,\Sigma_5,\Sigma_6,\Sigma_7$ form the
catch for $e_1$. But the catch must be of
the first kind, and $e_1$ must connect $v_7$ to
$v_5$. This is a contradiction.
If $e_1$ has positive slope,
then a similar argument shows that actually
$e_1$ connects $v_7$ to $v_4$, another
contradiction.

\begin{center}
\resizebox{!}{3.5in}{\includegraphics{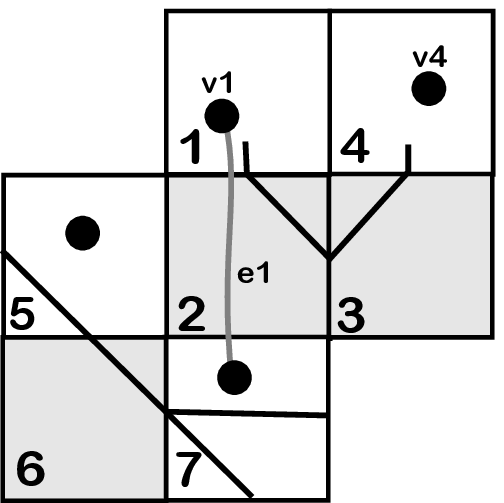}}
\newline
{\bf Figure 5.14:\/}  Case 3
\end{center}

\subsection{Length Four Chains: Case D}

There is only one remaining case, and this
case actually occurs.
Suppose that $\Sigma_1,\Sigma_2,\Sigma_3,\Sigma_4$ is
a length $4$ linked chain, making a zig-zag pattern
as in Figure 5.15 below.

As in Case 3C, it could happen that the
edge $\Sigma_1 \cap \Sigma_2$ is an unused side
for $\Sigma_1$.  In this case, there is an offending
edge $e$ associated to $\Sigma_1 \cap \Sigma_2$.
The catch for $e$ must be of the second
kind because $\Sigma_3$ is grid empty.  The
catch involves $\Sigma_1,\Sigma_2,\Sigma_3,\Sigma_4,\Sigma_{10}$,
and by definition $e$ must end in $\Sigma_4$ as desired.
So, we just have to worry about the
case when some graph edge $e_1$ is incidenct to
$v_1$ and crosses $\Sigma_1 \cap \Sigma_2$.

\begin{center}
\resizebox{!}{3.5in}{\includegraphics{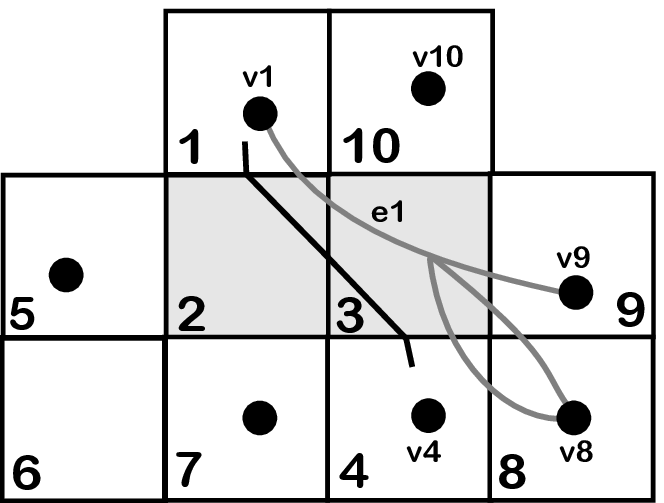}}
\newline
{\bf Figure 5.15:\/}  zig-zag case
\end{center}

The same reasons as in Case 4C rule out the
possibility that $e_1$ ends in
$\Sigma_5$, $\Sigma_6$, or $\Sigma_7$.
We just have to rule out $e_1$
ending in $\Sigma_8$ or $\Sigma_9$.
The same remarks as in Case 4C apply to the drawing
of the square $\Sigma_8$: This square could
be either grid full or grid empty, but in the
case relative to the argument it is grid full.
\newline
\newline
\noindent
{\bf Case 1:\/} Suppose $e_1$ ends in $\Sigma_8$ and
crosses $\Sigma_4 \cap \Sigma_8$.  The argument in
Case 4B2 rules this out: Two graph edges would be incident
to $\Sigma_8$ and would cross $\Sigma_4 \cap \Sigma_8$.
\newline
\newline
{\bf Case 2:\/} Suppose $e_1$ ends in $\Sigma_8$ and
crosses $\Sigma_9 \cap \Sigma_8$.  If $e_1$ is not
an offending edge for $\Sigma_8$ then we get the
same contradiction as in Case 1.  But if $e_1$ is
an offending edge for $\Sigma_8$ then, inspecting
the nature of the catches, we see that $e_1$ must
connect $v_8$ to $v_{10}$.
\newline
\newline
{\bf Case 3:\/}
Suppose $e_1$ ends in $\Sigma_9$.
Looking at the portion of the plaid segment in $\Sigma_3$,
we see that $\Sigma_3 \cap \Sigma_9$ is not in the
plaid edge set for $\Sigma_9$.  Hence $e_1$ is an
offending edge for $\Sigma_9$.  The catch
for $e_1$ must be of the first kind, because
$\Sigma_{10}$ is grid full.  But then
$e_1$ connects $v_9$ to $v_{10}$, a contradiction.

\newpage

\section{The Plaid Master Picture Theorem}
\label{PLAID}

\subsection{The Spaces}
\label{plaidspace}

Define
\begin{equation}
\widehat X=\R^3 \times [0,1].
\end{equation}
The coordinates on $\widehat X$ are
given by $(x,z,y,P)$.  We think of
\begin{equation}
P=\frac{2A}{1+A}, \hskip 30 pt A=p/q,
\end{equation}
but $P$ is allowed to take on any real value in $[0,1]$.
Define the following affine transformations of $\widehat X$.

\begin{itemize}
\item $T_X(x,y,z,P)=(x+2,y+P,z+P,P)$.
\item $T_Y(x,y,z,P)=(x,y+2,z,P)$;
\item $T_Z(x,y,z,P)=(x,y,z+2,P)$;
\end{itemize}
Define two abelian groups of affine transformations:
\begin{equation}
\widehat \Lambda_1=\langle T_X,T_Y,T_Z \rangle,
\hskip 30 pt
\widehat \Lambda_2=\langle T_X^2,T_Y,T_Z \rangle.
\end{equation}
Finally define
\begin{equation}
X_k=\widehat X/\Lambda_k, \hskip 30 pt k=1,2.
\end{equation}

The space $X_2$ is a double cover of $X_1$.
Both spaces should be considered flat affine
manifolds - i.e. manifolds whose overlap functions
are restrictions of affine transformations.
All the affine transformations in sight
preserve the slices $\R^3 \times \{P\}$ and
act as translations on these slices.
This $X_1$ and $X_2$ are fibered by 
$3$-dimensional Euclidean tori.
\newline
\newline
{\bf Remark:\/}
The reason why we keep track of two spaces
is that a suitable partition of $X_1$ into
polytopes determines the unoriented polygons
in the plaid model whereas a suitable partition
of $X_2$ determines the oriented polygons
in the plaid model. 
What we will do is describe a partition of
$X_2$ that is invariant under the action
of $T_X$. To get a partition of $X_1$ we
simply take the quotient and forget some
auxiliary information about the labelings.
The partition of $X_2$ has more information
but we sometimes consider the picture in
$X_1$ because it makes some formulas simpler
and cuts down on the computations.

\subsection{The Plaid Master Picture Theorem}

Here we will describe the Isomorphism
Theorem from [{\bf S1\/}]. (Again, I am
rechristening this result the
Plaid Master Picture Theorem.)
The {\it plaid grid\/} is defined to be the
set $G_{\Pi}$ of centered of integer unit squares.
When we know that we are taking about the
plaid grid, as opposed to the graph grid,
we will set $G=G_{\Pi}$.
For each parameter $A \in (0,1)$, there is a 
linear map $\Phi_A: G \to \widehat X$,
given by
\begin{equation}
\label{classify}
\Phi_{A,k}(x,y)=(2Px+2y,2Px,2Px+2Py,P),
\hskip 30 pt
P=\frac{2A}{1+A}.
\end{equation}

At the same time $\widehat X$ has a $\Lambda_1$-invariant
tiling by $4$ dimensional convex integral polytopes.
Modulo the action of $\Lambda_1$ there are $26$ such
polytopes.   We list the coordinates of $26$ representatives
of the orbit equivalence classes below. We call these
representatives the {\it fundamental polytopes\/}.  Their
union is a fundamental domain for the action of
$\Lambda_1$.  By taking the $\Lambda_1$-orbit of
the fundamental polytopes, we
get a partition of $\widehat X$ into infinitely many
convex integral polytopes.

The polytopes in the partition have an auxiliary labeling.
Either the labeling is by the empty set, or it is
by an unordered pair of letters in
$\{S,W,N,E\}$.
Here is a precise statement of the Plaid Master Picture Theorem.

\begin{theorem}
\label{iso1}
For any even rational parameter $A=p/q$ the following
is true.
At the point $c \in G$ the plaid tile is
determined by the labeling
of the polytope which contains $\Phi_A(c)$ in its interior.
\end{theorem}

In case the labeling is $(0,0)$ the tile centered at $c$
is the empty tile.  Otherwise, the tile is determined
by the correspondence above.  For instance, if
$\Phi_A(c)$ lies in a polytope labeled $\{1,2\}$, then
the tile at $c$ has the connector joining the
South edge to the West edge.  Implicit in the statement
of the theorem is that $\Phi(c)$ does always lie in
the interior of one of the Polytopes.

In [{\bf S1\/}] we proved an enhancement of
Theorem \ref{iso1}.  The labeling of the polytopes
by unordered pairs can be enhanced to a labeling
by {\it ordered pairs\/} provided that we only
require that the labeling has $\Lambda_2$ symmetry
and moreover that the action of the affine map
$T_X$ reverses the ordering of the labels.  
Thus the union of the $26$ polytopes below and
their images under $T_X$, with the labels reversed,
gives a fundamental domain for the action of
$\Lambda_2$ on $\widehat X$.   The enhanced version
of Theorem \ref{iso1} is as follows:

\begin{theorem}
\label{iso2}
For any even rational parameter $A=p/q$ the following
is true.
At the point $c \in G$ the oriented plaid tile is
determined by the enhanced labeling
of the polytope which contains $\Phi_A(c)$ in its interior.
\end{theorem}

We explain what we mean by way of example.  Suppose that
$\Phi_A(c)$ lies in a polytope labeled $\{1,2\}$, then
the tile at $c$ has the connector joining the
South edge to the West edge and pointing from South
to West.

We can think of the classifying map $\Phi_A$ as having its
range in $X_1$.  In this case, the tiling of
$\widehat X$ descends to a tiling of $X_1$ by $26$
convex polytopes, each of which is labeled by an unordered
pair of integers, as above.
Likewise, we can think of the classifying map
$\Phi_A$ as having its range in $X_2$.  In this case,
the tiling of
$\widehat X$ descents to a tiling of $X_2$ by $56$
convex polytopes, each of which is labeled by
an ordered pair of integers.  

Note that
Theorem \ref{iso2} is stronger than
Theorem \ref{iso1}, but for our proof of
the Quasi-Isomorphism Theorem it is useful sometimes
to keep Theorem \ref{iso1} in mind. The calculations
are simpler.
On the other hand, the advantage to working with
$X_2$ is that this space has a natural interpretation
as an affine PET.  We will explain this interpretation
below.  Thus, for the Quasi-Isomorphism Theorem, it is
best to keep track of both spaces and partitions.

\subsection{The Partition}
\label{plaidpet}

In \S \ref{plaid_polytopes}
we list the vertices and
oriented labels of $10$ 
polytopes.  All the polytopes in
the partition of $\widehat X$ can obtained
from these using the group generated
by $\Lambda_1$ and the
following two additional elements:
\begin{itemize}
\item {\bf Negation:\/} The map
$(x,y,z,A) \to (-x,-y,-z,A)$ preserves
the partition and changes the labels as
follows: $N$ and $S$ are swapped and $E$ and $W$
are swapped.
\item {\bf Flipping:\/}
The map $(x,y,z,A) \to (x,z,y,A)$ preserves the
partition and changes the labels as follows:
$N$ and $S$ are swapped, and the order of
the label is reversed.
\end{itemize}
Our listing of polytopes, together with the explicit
formulas above, gives a complete account of
Theorem \ref{iso1} and Theorem \ref{iso2}.
However, this description is not very
satisfying from a geometric point of view.
The partition is actually quite beautiful.

In [{\bf S1\/}] we gave a careful 
geometric description of the partition.
Rather than repeat {\it verbatim\/} what
we said in [{\bf S1\/}], here we describe
enough of the
partition so that the very motivated reader
could reconstruct it from scratch.
In our computer program, we have an interface
which allows the reader to navigate through
the partition and instantly see the
features we discuss here.  We describe
the partition of $X_1$.  The partition of
$X_2$ is obtained using the action of
$\Lambda_1$, as described above.

The rectangular
solid
\begin{equation}
[-1,1]^3 \times [0,1].
\end{equation}
serves as a fundamental domain for the action of
$\Lambda_1$ on $\widehat X$.  We think of this
space as a fiber bundle over the $(x,A)$ plane.
The base space $B$ is the rectangle $[-1,1] \times [0,1]$.
The space $B$ has a partition into $3$ triangles, as
shown in Figure 6.1.

\begin{center}
\resizebox{!}{2.6in}{\includegraphics{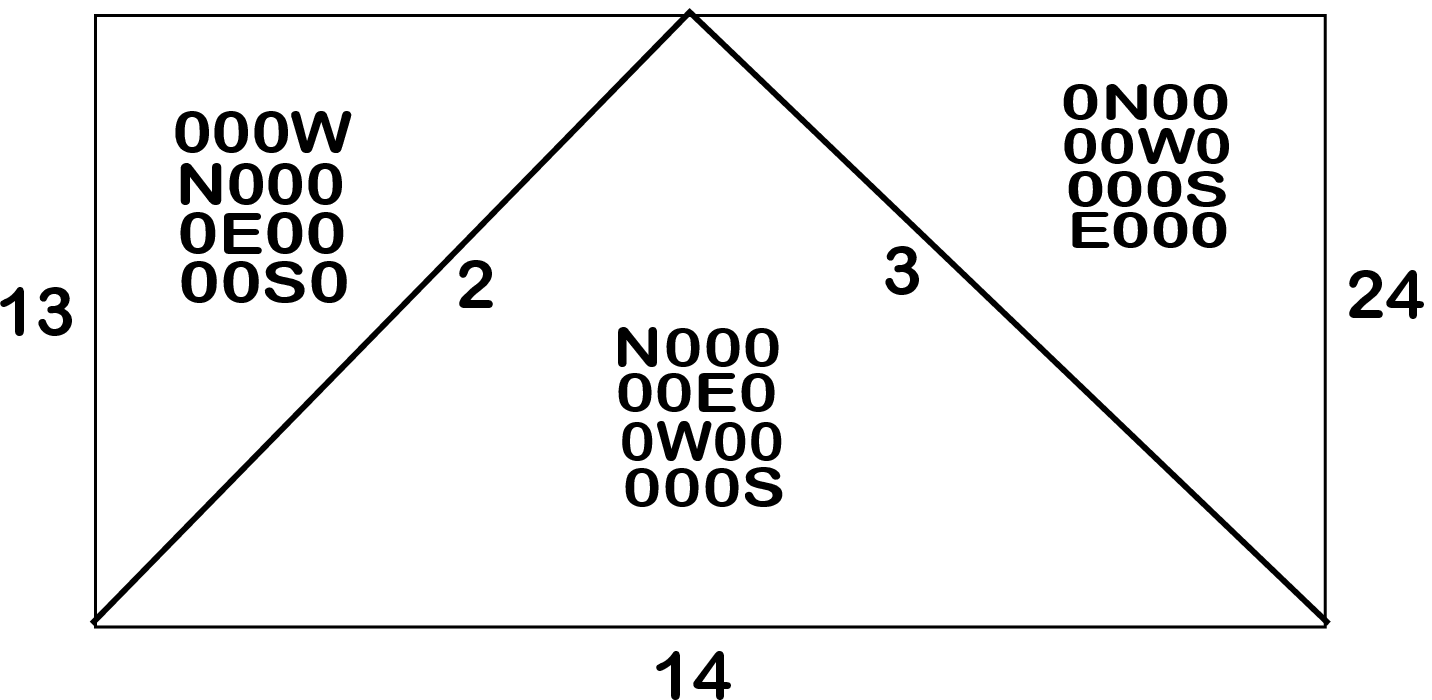}}
\newline
{\bf Figure 6.1:\/} The $(x,A)$ base space.
\end{center}

Each triangular region has been assigned a
$4 \times 4$ matrix which we will explain
momentarily.  The slices above each open triangle
in the partition of $B$ consists of a
$4 \times 4$ grid of rectangles.  For of the
rectangles, corresponding to the polytopes
labeled by the emptyset are squares. (In Figure 6.2
we replace the emptyset by the symbols
$\{S,W,N,E\}$ because this makes the rest of
the labeling pattern more clear.  The pattern of
the squares is indicated by the nonzero entries
in the matrices.  For instance, the grid shown in Figure 6.2 
corresponds to a fiber over the left triangle.
We have also indicated the labels in the picture.

\begin{center}
\resizebox{!}{3.8in}{\includegraphics{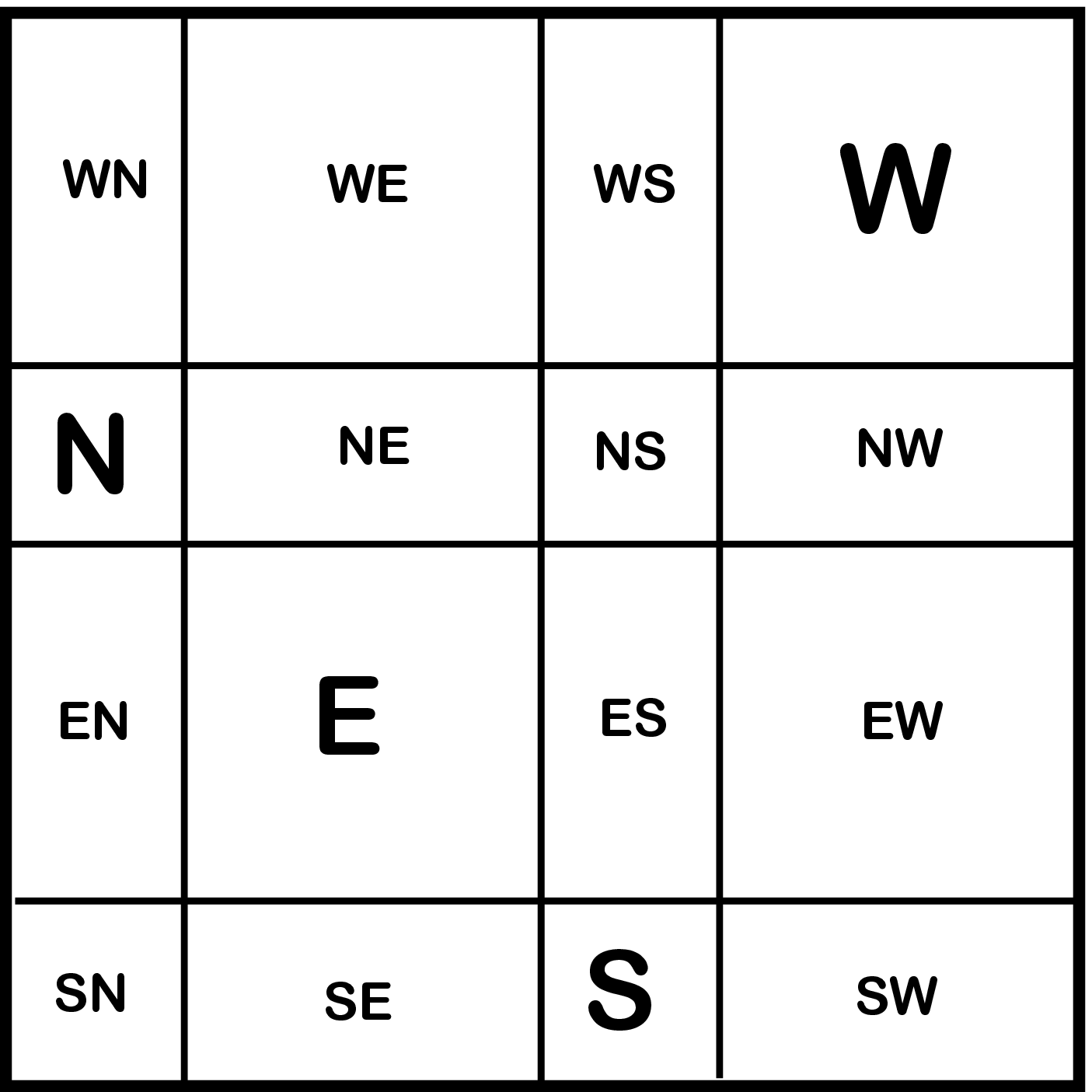}}
\newline
{\bf Figure 6.2:\/} The checkerboard partition
\end{center}

The labels can be derived from the matrices in
Figure 6.1.
Every nonsquare rectangle in a column has the
same first label and every rectangle in a
row has the same second label.

The edge labels in Figure 6.1 encode degenerations
in the fibers.
Away from
the singular edges of the base space, each
polytope is a rectangle bundle over an
open triangle. 
As one approaches the edges of the triangles in
the base, the partition degenerates in that
some of the rectangles shrink to line segments
or points. The numbers on the edges in 
Figure 6.1 indicate which rows and columns
degenerate.  For instance, as we approach
the edge labeled $2$, the thickness of
second column from the left tends to $0$
and at the same time the thickness of the
second row from the bottom tends to $0$.
The whole (unlabeled) picture is symmetric
with respect to reflection in the line $y=x$.

From the description of these degenerations
and from the fact that we know we are looking
at slices of convex integral polytopes, one
can actually reconstruct the entire partition
from the labelings in Figure 6.1.

\subsection{Curve Following Dynamics}
\label{CFD}

Recall that $P=2A/(1+A)$.
Recall also that $\cal C$ is the set of centers
of integer unit squares in the plane.
We fix a parameter $A=p/q$ and consider
the map $f_A: G \to G$
that comes from simply following the arrows
on the tiles.  Given $c_0 \in G$, the
new point $C_1=f_A(c_0)$ is defined to be
the center of the tile into which the
tile at $c_0$ points. For instance,
the tile centered at $c_0$ is NE, then
$c_1=c_0+(1,0)$. In case the tile
centered at $c_0$ is empty, we have
$c_1=c_0$.

In view of Theorem \ref{iso2} there is a
corresponding map $F: \widehat X \to \widehat X$
such that
\begin{equation}
F\circ \Phi_A = \Phi_A \circ f_A.
\end{equation}
Looking at Equation \ref{classify} we see that $F$ has
the following description.
\begin{itemize}
\item $F$ is the identity on tiles labeled by the
empty set.
\item On all polytopes whose labels end in $N$, we
have $$F(x,y,z,P)=(x,y,z,P)+(2,0,2P,0).$$
\item On all polytopes whose labels end in $S$, we
have $$F(x,y,z,P)=(x,y,z,P)-(2,0,2P,0).$$
\item On all polytopes whose labels end in $E$, we
have $$F(x,y,z,P)=(x,y,z,P)+(2P,2P,2P,0).$$
\item On all polytopes whose labels end in $N$, we
have $$F(x,y,z,P)=(x,y,z,P)-(2P,2P,2P,0).$$
\end{itemize}
We can also interpret $F$ as a map 
$F_2: X_2 \to X_2$, because everything in
sight is $\Lambda_2$ invariant.  With
this interpretation, $(X_2,F_2)$ is the
fibered integral affine PET.
What we mean is the $4$ dimensional
system is piecewise affine, and that
the PET preserves
the $P$-slices and acts as an ordinary
piecewise translation in each fiber.
The $4$ dimensional system is
integral in the sense that all the polytopes
in the partition have integer vertices.

\subsection{The Triple Partition}
\label{rtf}

Suppose that ${\cal P\/}_1$ and ${\cal P\/}_2$ are
two partitions of $\widehat X$ into polytopes,
say
\begin{equation}
{\cal P\/}_k=\bigcup_i P_{i,k}.
\end{equation}
We define the {\it common refinement\/} to be
\begin{equation}
{\cal P\/}_1 \bowtie {\cal P\/}_2=\bigcup_{i,j} P_{i,1} \cap P_{j,2}.
\end{equation}
This is the usual definition.
In case the polytopes in the two partitions have rational
vertices, the polytopes, in the common refinement will also
have rational vertices.   This construction may be iterated,
so that we can take the $n$-fold common refinement of
$n$ partitions.

Let ${\cal P\/}_0$ be the partition of $\widehat X$ that
we have described above.  We define
\begin{equation}
{\cal P\/}_k=F^k({\cal P\/}_0).
\end{equation}
In other words, we take the original partition and apply
a power of the PET dynamics.  We define the
{\it triple partition\/} to be the common refinement
\begin{equation}
{\cal TP\/}={\cal P\/}_{-1} \bowtie {\cal P\/}_0 \bowtie {\cal P\/}_1.
\end{equation}
We will prove that every vertex of every polytope
in ${\cal TP\/}$ has a coordinates which are
divisors of $60$.  That is, if
$Q$ is a polytope of ${\cal TP\/}$, then
$60Q$ is an integer polytope.

Each polytope in $\cal TP$ has a $6$ letter label.  We
simply concatenate the labels for each of the $3$
polytopes involved in the intersection, starting
with the label for the polytope in ${\cal P\/}_{-1}$
and ending with the label for the polytope in
${\cal P\/}_1$.  This $6$ letter label has the
following meaning.  We look at the label of
the polytope containing $\Phi_A(c)$ and the
label determines the shape of the length $3$ arc
of the oriented plaid component that is centered
at $c$.  Figure 6.3 illustrates this principle
with two examples.  In Figure 6.3 we show the
label and the corresponding path.  We hope that
these examples suffice to convey the general
idea.

\begin{center}
\resizebox{!}{2in}{\includegraphics{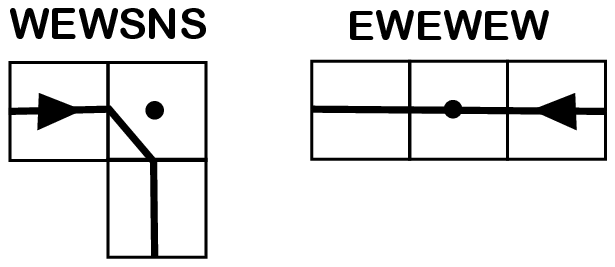}}
\newline
{\bf Figure 6.3:\/} The meaning of the $6$ letter labels
\end{center}

We define the {\it reduced triple partition\/} to be
partition
\begin{equation}
{\cal RTP\/}={\cal TP\/} \bowtie ([-1,1]^3 \times [0,1]).
\end{equation}
Essentially we are taking the pieces of
${\cal TP\/}$ which lie inside our favorite
fundamental domain for $\Lambda_1$, but in case
any of the pieces slop over the boundary of
this fundamental domain, we chop them off.
Once again, for every $Q \in {\cal RTP\/}$, the
scaled polytope $60Q$ is integral.

We think of the polytopes in ${\cal RTP\/}$ as being
labeled by the same $6$-letter labels, modulo
reversal of the labeling.  With this interpretation,
we can interpret that map
$\Phi_A G \to \widehat X$ as a classifying map
for unoriented arcs of the plaid model having
combinatorial length $3$. We just use the action
of $\Lambda_1$ to move $\Phi_A(c)$ into
${\cal RTP\/}$ and then we read off the label.

There is one minor issue that we need to address.
It follows from Theorem \ref{iso2} that
$\Phi_A(c)$ always lies in the interior of
a polytope of ${\cal TP\/}$, but it might happen
that there are several images of
$\Phi_A(c)$ on the boundary of the fundamental
domain.  This is the usual problem with
fundamental domains.  However, in this situation,
all the images will lie in polyopes having
the same labels.  So, even when there is some
ambiguity in interpreting $\Phi_A$ as a map
from $G$ into ${\cal RTP\/}$, the
ambiguity is harmless.

Even though it is somewhat less natural than 
${\cal TP\/}$, we will use
${\cal RTP\/}$ for the proof of the
Quasi-Isomorphism.  The formula for
the projective intertwiner discussed in the
introduction is simpler when we restrict
our attention to $X_1$.  The projective
intertwiner is the most important
ingredient in the proof, and so we wanted
to make it as nice as possible.

\subsection{The Fundamental Polytopes}
\label{plaid_polytopes}

Here is the list of the $10$ labeled polytopes.
The label $X$ denotes that the label is the
empty set. Also, the last polytope listed has
$12$ vertices, and the listing is spread over
$2$ lines.
\newline

\lefteqn{\left[\matrix{-1 \cr -1 \cr -1 \cr 0}\right]
\hskip 5 pt\left[\matrix{-1 \cr +1 \cr -1 \cr 0}\right]
\hskip 5 pt\left[\matrix{+1 \cr +1 \cr -1 \cr 0}\right]
\hskip 5 pt\left[\matrix{+1 \cr +1 \cr +1 \cr 0}\right]
\hskip 5 pt\left[\matrix{0 \cr 0 \cr 0 \cr +1}\right]
\hskip 20 pt (W,E)}
\vspace{4mm}

\lefteqn{\left[\matrix{+1 \cr +1 \cr -1 \cr 0}\right]
\hskip 5 pt\left[\matrix{+1 \cr +1 \cr 0 \cr +1}\right]
\hskip 5 pt\left[\matrix{+1 \cr +1 \cr -1 \cr +1}\right]
\hskip 5 pt\left[\matrix{0 \cr +1 \cr -1 \cr +1}\right]
\hskip 5 pt\left[\matrix{0 \cr 0 \cr -1 \cr +1}\right]
\hskip 20 pt (E,S)}
\vspace{4mm}

\lefteqn{\left[\matrix{+1 \cr -1 \cr -1 \cr 0}\right]
\hskip 5 pt\left[\matrix{+1 \cr 0 \cr 0 \cr +1}\right]
\hskip 5 pt\left[\matrix{+1 \cr 0 \cr -1 \cr +1}\right]
\hskip 5 pt\left[\matrix{0 \cr 0 \cr -1 \cr +1}\right]
\hskip 5 pt\left[\matrix{0 \cr -1 \cr -1 \cr +1}\right]
\hskip 20 pt (E,N)}
\vspace{4mm}

\lefteqn{\left[\matrix{-1 \cr +1 \cr -1 \cr 0}\right]
\hskip 5 pt\left[\matrix{+1 \cr +1 \cr -1 \cr 0}\right]
\hskip 5 pt\left[\matrix{+1 \cr +1 \cr +1 \cr 0}\right]
\hskip 5 pt\left[\matrix{0 \cr 0 \cr 0 \cr +1}\right]
\hskip 5 pt\left[\matrix{0 \cr +1 \cr 0 \cr +1}\right]
\hskip 5 pt\left[\matrix{+1 \cr +1 \cr 0 \cr +1}\right]
\hskip 5 pt\left[\matrix{+1 \cr +1 \cr +1 \cr +1}\right]
\hskip 20 pt (W,S)}
\vspace{4mm}

\lefteqn{\left[\matrix{-1 \cr -1 \cr -1 \cr 0}\right]
\hskip 5 pt\left[\matrix{+1 \cr -1 \cr -1 \cr 0}\right]
\hskip 5 pt\left[\matrix{+1 \cr +1 \cr -1 \cr 0}\right]
\hskip 5 pt\left[\matrix{0 \cr 0 \cr -1 \cr +1}\right]
\hskip 5 pt\left[\matrix{0 \cr 0 \cr 0 \cr +1}\right]
\hskip 5 pt\left[\matrix{+1 \cr 0 \cr 0 \cr +1}\right]
\hskip 5 pt\left[\matrix{+1 \cr +1 \cr 0 \cr +1}\right]
\hskip 20 pt (S,W)}
\vspace{4mm}

\lefteqn{\left[\matrix{-1 \cr +1 \cr +1 \cr 0}\right]
\hskip 5 pt\left[\matrix{-1 \cr 0 \cr 0 \cr +1}\right]
\hskip 5 pt\left[\matrix{-1 \cr 0 \cr +1 \cr +1}\right]
\hskip 5 pt\left[\matrix{-1 \cr +1 \cr 0 \cr +1}\right]
\hskip 5 pt\left[\matrix{-1 \cr +1 \cr +1 \cr +1}\right]
\hskip 5 pt\left[\matrix{0 \cr +1 \cr +1 \cr +1}\right]
\hskip 5 pt\left[\matrix{-2 \cr 0 \cr 0 \cr +1}\right]
\hskip 20 pt X}
\vspace{4mm}

\lefteqn{\left[\matrix{-1 \cr -1 \cr -1 \cr 0}\right]
\hskip 5 pt\left[\matrix{+1 \cr -1 \cr -1 \cr 0}\right]
\hskip 5 pt\left[\matrix{-1 \cr -1 \cr -1 \cr +1}\right]
\hskip 5 pt\left[\matrix{0 \cr 0 \cr 0 \cr +1}\right]
\hskip 5 pt\left[\matrix{0 \cr 0 \cr -1 \cr +1}\right]
\hskip 5 pt\left[\matrix{0 \cr -1 \cr 0 \cr +1}\right]
\hskip 5 pt\left[\matrix{0 \cr -1 \cr -1 \cr +1}\right]
\hskip 5 pt\left[\matrix{+1 \cr 0 \cr 0 \cr +1}\right]
\hskip 20 pt (S,N)}
\vspace{4mm}

\lefteqn{\left[\matrix{+1 \cr +1 \cr -1 \cr 0}\right]
\hskip 5 pt\left[\matrix{+1 \cr -1 \cr -1 \cr 0}\right]
\hskip 5 pt\left[\matrix{+1 \cr +1 \cr 0 \cr +1}\right]
\hskip 5 pt\left[\matrix{+1 \cr 0 \cr 0 \cr +1}\right]
\hskip 5 pt\left[\matrix{+1 \cr +1 \cr -1 \cr +1}\right]
\hskip 5 pt\left[\matrix{+1 \cr 0 \cr -1 \cr +1}\right]
\hskip 5 pt\left[\matrix{0 \cr 0 \cr -1 \cr +1}\right]
\hskip 5 pt\left[\matrix{2 \cr +1 \cr 0 \cr +1}\right]
\hskip 20 pt (E,W)}
\vspace{4mm}

\lefteqn{\left[\matrix{-1 \cr -1 \cr +1 \cr 0}\right]
\hskip 5 pt\left[\matrix{-1 \cr -1 \cr 0 \cr +1}\right]
\hskip 5 pt\left[\matrix{0 \cr -1 \cr 0 \cr +1}\right]
\hskip 5 pt\left[\matrix{0 \cr 0 \cr 0 \cr +1}\right]
\hskip 5 pt\left[\matrix{0 \cr -1 \cr +1 \cr +1}\right]
\hskip 5 pt\left[\matrix{0 \cr 0 \cr +1 \cr +1}\right]
\hskip 5 pt\left[\matrix{+1 \cr 0 \cr +1 \cr +1}\right]
\hskip 5 pt\left[\matrix{+1 \cr -1 \cr +1 \cr 0}\right]
\hskip 20 pt X}
\vspace{4mm}

\lefteqn{\left[\matrix{-1 \cr -1 \cr -1 \cr 0}\right]
\hskip 5 pt\left[\matrix{-1 \cr +1 \cr +1 \cr 0}\right]
\hskip 5 pt\left[\matrix{-1 \cr -1 \cr +1 \cr 0}\right]
\hskip 5 pt\left[\matrix{-1 \cr +1 \cr -1 \cr 0}\right]
\hskip 5 pt\left[\matrix{+1 \cr +1 \cr +1 \cr 0}\right]
\hskip 5 pt\left[\matrix{-1 \cr -1 \cr -1 \cr +1}\right]}
\vspace{2mm}
\lefteqn{\left[\matrix{-1 \cr 0 \cr -1 \cr +1}\right]
\hskip 5 pt\left[\matrix{-1 \cr -1 \cr 0 \cr +1}\right]
\hskip 5 pt\left[\matrix{-1 \cr 0 \cr 0 \cr +1}\right]
\hskip 5 pt\left[\matrix{0 \cr 0 \cr 0 \cr +1}\right]
\hskip 5 pt\left[\matrix{-3 \cr -1 \cr -1 \cr 0}\right]
\hskip 5 pt\left[\matrix{-2 \cr -1 \cr -1 \cr +1}\right]
\hskip 20 pt X}

\newpage

\section{The Graph Master Picture Theorem}

\subsection{Heuristic Discussion}

The Master Picture Theorem from [{\bf S0\/}] 
says that the arithmetic
graph has a classifying map and space much like the plaid
model does.  Following this section, we will give an
account of the Master Picture Theorem.  Here, with a
view towards making this paper more self-contained,
I'd like to give some idea why the Master Picture Theorem
is true.

Let $K=K_A$ be the kite (implicitly)
involved in the Quasi-Isomorphism Thorem.
Given an edge $e$ of $K$, we form a strip
$\Sigma_e$ as follows:  One boundary line of $\Sigma_e$
is $L_e$, the line containing $e$.  The other boundary
line $L_e'$ is parallel to $L_e$ and such that the point
of $P$ farthest from $L_e$ lies halfway between
$L_e$ and $L_e'$.  These $4$ strips $\Sigma_1,...,\Sigma_4$
are ordered cylically, according to their slopes.
Figure 7.1 shows the picture.  The big disk in the middle
covers up the mess in the middle of the picture.

\begin{center}
\resizebox{!}{4in}{\includegraphics{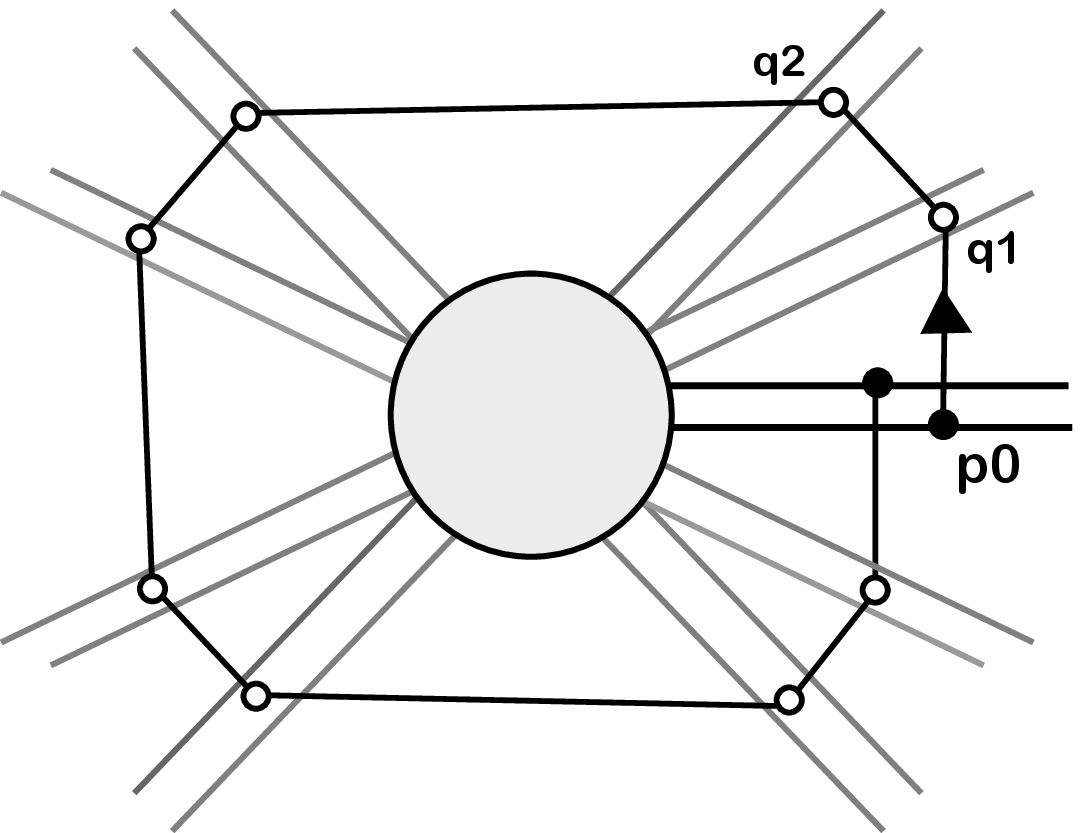}}
\newline
{\bf Figure 7.1:\/} The strips defined by the kite
\end{center}

Let $\psi$ denote the second iterate of the outer billiards
map. Let $\Xi$ denote the distinguished pair of rays
involved in the construction of the arithmetic graph.
Let $\Psi: \Xi \to \Xi$ denote the first return of
$\psi$ to $\Xi$.
A general feature of $\psi$ is the
orbit $p_0,\psi(p_0),\psi^2(p_0),...$ moves along a
straight line, in evenly spaced points, until it lands
in one of the strips.  We call the points where
the orbit lands in the strips the {\it turning points\/}.
It seems plausible (and in
fact is true) that the way the orbit returns to
$\Xi$ is determined by the relative positions of
the turning points inside the strips -- i.e., how
they sit with respect to the boundaries of the strips.

To formalize this idea, we define a map 
$$
\Phi': \Xi \to [0,1]^8,
$$
as follows.
Starting with the point $p_0$, we let
$q_1,...,q_8$ denote the successive points where
the $\psi$-orbit of $p_0$ intersects a strip.
We let $t_k \in [0,1]$ the position of $q_k$ in
the strip relative to the boundary components.
Thus, $t_k \approx 0$ if $q_k$ lies near one
boundary component of the strip,
$t_k \approx 1$ if $q_k$ lies near the other 
boundary component, and $t_k \approx 1/2$ if
$q_k$ lies near the midline of the strip.
The map is then
$$\Phi'(p_0)=(t_1,...,t_8).$$

Thanks to some relations between the slopes
of the strips, the image
$\Phi'(\Xi)$ is contained
in a $3$-dimensional slice of $[0,1]^8$.
The image $\Phi'(p_0)$ classifies
the {\it first return vector\/}
$\Psi(p_0)-p_0$ in the following sense.
There is a partition of the $3$ dimensional
slice into polyhedra such that
$p_0$ and $p_0'$ have the same first return
vector provided that $\Phi'(p_0)$ and
$\Phi'(p_0')$ lie in the same polyhedron.

In [{\bf S0\/}] we make this construction
carefully, then introduce a suitable affine
change coordinates, and then finally
interpret the result directly in terms
of the arithmetic graph. (The arithmetic
graph encodes the first return vectors.)
When all this is done, we get the
Master Picture Theorem.

It is worth pointing out that the existence
of this compactification works in much
more generality.  In my paper [{\bf S3\/}],
I give a similar result for any convex
polygon with no parallel sides.  Rather
than consider the return map to 
some distinguished pair of rays -- a construction
which doesn't even make sense in general -- I consider
a return map to one of the strips.  This
first return map has a compactification whose
dimension is comparable to the number of
sides of the polygon.  One could
probably deduce the Master Picture Theorem
from the general compactification theorem
in [{\bf S3\/}], but I haven't tried.

\subsection{The Master Picture Theorem}
\label{mpt}

We will give a simplified
presentation of what is contained in
[{\bf S0\/}, \S 6], because some of the fine
points there are irrelevant to us.  Also,
in order to line up the picture better with
the plaid model, we switch the $x$ and $z$
coordinates.

We fix $A=p/q$ as above.
As in the case of the plaid PET, we
work in the space $\widehat X=\R^3 \times [0,1]$.
This time, we use the coordinates
$(x,y,z,A)$ on $\widehat X$.

Let $\Lambda$ denote the
abelian group of generated by
the following affine transformations.
\begin{itemize}
\item $T_X(x,y,z,A)=(x+1,y-1,z-1,A)$.
\item $T_Y(x,y,z,A)=(x,y+1+A,z+1-A,A)$.
\item $T_Z(x,y,z,A)=(x,y,z+1+A,A)$.
\end{itemize}
For each parameter $A$, the rectangular solid
\begin{equation}
R_A=[0,1] \times [0,1+A] \times [0,1+A] \times \{A\}
\end{equation}
serves as a fundamental domain for the
action of $\Lambda$ on $\R^3 \times \{A\}$.
The quotient is some flat $3$-torus which
depends on $A$.  The union
\begin{equation}
\label{fd2}
R=\bigcup_{A \in [0,1]} R_A,
\end{equation}
is a convex integer polytope.

We introduce a map 
$\Phi'_A: \Z^2 \to \widehat X$, as follows.
\begin{equation}
\label{graphclass}
\phi_A(m,n)=(t,t,t,A), \hskip 30 pt
t=Am+n+\frac{1}{2q}.
\end{equation}

There are two $\Lambda$-invariant partitions
of $\widehat X$ into convex integral polytopes.
Each polytope in each partition is labeled
by a pair of integers $(i,j) \in \{-1,0,1\}$.
The labeling has the property that the
local structure of the arithmetic graph at the
point $c \in \Z^2$ is determined by the
labels of the polytopes in each partition which
contain the image $\Phi'(c)$.
Implicit in this statement is the fact that
$\Phi'(c)$ lies in the interior of a polytope
in each partition.

Here is what we mean
more precisely.  Suppose that
$\Phi'(c)$ lies in a polytope in the first
partition with label $(i_1,j_1)$ and a
polytope in the second partition with
labels $(i_2,j_2)$.  Then the arithmetic
graph has the edge connecting
$c$ to $c+(1_1,j_1)$ and the edge
connecting $c$ to $c+(i_2,j_2)$.  We
have $(i_1,j_1)=(0,0)$ if and only
if $(i_2,j_2)=(0,0)$.  In this case,
$c$ is an isolated point in the graph.

For each partition, a certain union of
$14$ polytopes forms a fundamental
domain for the action of $\Lambda$.
In both cases, the fundamental domain
is an integral translate of the
polytope $R$ from Equation \ref{fd2}.
The polytopes in the two partitions
are related as follows.  Define
\begin{equation}
I(x,y,z,A)=(1,A,2+A,A)-(x,y,z,0).
\end{equation}
Then
\begin{itemize}
\item $Q_j=I(P_j)$ for all $j=1,...,14$.
\item The label of $Q_j$ is the negative of
the label of $P_j$.
\end{itemize}
The map $I$ is an involution, so we
have $I(Q_j)=P_j$ as well.

We call these two partitions the $(+)$ graph partition
and the $(-)$ graph partition.  We denote these
partitions by ${\cal G\/}_+$ and ${\cal G\/}_-$.
The {\it double partition\/} is
${\cal G\/}_+ \bowtie {\cal G\/}_-$, the
common refinement of the two partitions.
Each polytope $Z$
in the double partition is labeled by a
quadruple $(i_+,j_+,i_-,j_i)$, where
$(i_{\pm},j_{\pm})$ is the label of the
polytope in ${\cal P\/}_{\pm}$ containing $Z$.
The polyyope of the double partition
which contains $\Phi'_A(c)$ determines the
two unoriented edges of the arithmetic
graph incident to $c$.
The polytopes in the double partition have
rational coordinates.  We did not work out
the double partition explicitly, because we
don't need to know it explicitly.  

In [{\bf S0\/}, \S 6] we give a detailed geometric
description of (translates of) the partitions
${\cal G\/}_+$ and ${\cal G\/}_-$. The geometric
description is more intricate than what we did
for the plaid PET.  We will not repeat the
description here.  The reader can get a very clear
picture of the partitions using our computer program.

In \S \ref{graphPET} we list the $14$ fundamental
polytopes of ${\cal G\/}_+$, together with their
labels.  This listing combines with the information
above to give a complete account of the statement
of the Master Picture Theorem.  My computer
program allows the user to navigate through
both partitions.  Using the program, one can
get a good visual sense of the partitions.

\subsection{Pulling Back the Maps}

Let $T$ be the canonical affine transformation
from \S 3. This map is given 
in Equation \ref{k2p} and the inverse
map is given in Equation \ref{k2p}.
We can interpret the Master Picture Theorem as a
statement about the structure of the affine
image $T(\Gamma)$, where $\Gamma$ is the 
arithmetic graph as defined in previous chapters.
After this
section, we will save words by
using the term {\it arithmetic graph\/}
to denote $T(\Gamma)$ rather than
$\Gamma$.  The vertices of
$T(\Gamma)$ are contained in the
graph grid $T(\Z^2)$.

To convert the Master Picture Theorem
into a statement about the desired
affine iamge, we simply pull back
the classifying map. We define
\begin{equation}
\Phi=\Phi' \circ T^{-1}.
\end{equation}
The domain for $\Phi$ is the
graph grid $G$.
A calculation shows that
\begin{equation}
\label{NICE}
\Phi(x,y)=(x,x,x,A).
\end{equation}
This nice equation suggests that there is
something canonical about the affine map
that appears in the Quasi-Isomorphism Theorem.

\subsection{Further Discussion}

There are three small differences between our
presentation of the Master Picture Theorem here
and the one given in [{\bf S0\/}].
For the reader who is interested in comparing
what we say here to what we say in
[{\bf S0\/}], we discuss those differences.

First, here we have one classifying map
map and two partitions, whereas in [{\bf S0\/}] we
have two classifying maps, differing from
each other by translations, and two slightly
different partitions.  The two partitions
in [{\bf S0\/}] are translates of the ones
here respectively by the vectors $(0,1,0,0)$
and $(-1,0,0,0)$.  When these changes are
made, the involution $I$ simply becomes
reflection in the midpoint of the fundamental
domain $R$, and the union of the $14$
fundamental polytopes, in either partition,
is precisely $R$.
In short, the partitions in [{\bf S0\/}] are
obtained from the ones here by translating
so that everything lies in $R$.  This 
picture is geometrically more appealing, but
it serves our purposes here to have one map
rather than two.

Second, even 
after we take into account the
translations we just discussed, the map
we have here is not quite the same as the
map in [{\bf S0\/}].  Were we to strictly
translate the map from [{\bf S0\/}] the
small term $1/2q$ in Equation
\ref{graphclass} would be replaced by an
infinitesimally small positive number $\iota=0_+$.
(If you like, you can interpret $0_+$ as an
infinitesimally small positive number in a nonstandard
field containing $\R$.)
Put another way, we set $\iota=0$, and then in those rare
cases when $p=\Phi'_A(m,n)$ lies in the boundary of
one of the polytopes, we select the polytope which
contains $p+(\epsilon,\epsilon,\epsilon,0)$ for all
sufficiently small $\epsilon$.  In fact, any choice
of $\iota$ in $(0,1/q)$ would give exactly the
same result.  Setting $\iota=0$ and then stipulating
the rules for handling boundary cases seemed at
the time to simplify the picture.  In hindsight,
the choice $\iota=1/2q$, as we make here, is more
canonical.  Also, it is better suited to our
present purposes.  See Equation \ref{NICE} in the next chapter.

Third, as we mentioned above, we have switched the
first and third coordinates.  Thus, the vertex
$(x,y,z,A)$ of a polytope here corresponds to the
vertex $(z,y,x,A)$ in [{\bf S0\/}].  Again, this
change makes the Master Picture Theorem line
up more gracefully with the Isomorphism Theorem
for the plaid model.  This switch also effects
the definition of the lattice $\Lambda$.

There is one more point we'd like to discuss.
Since the arithmetic graph comes from
following the dynamics of the first
return map, it is possible to orient
the paths in the arithmetic graph.
It is also possible to determine
when the first partition determines
the edge pointing to $c$ or the
edge pointing away to $c$.  This fine
point is not needed for the proof of
the Quasi-Isomorphism Theorem, but it
would be needed if we wanted to have
a well-defined affine PET associated
to the arithmetic graph, as we
constructed for the plaid model.

Here is the recipe for determining the
orientations from the classifying map.
We introduce the integer function
\begin{equation}
\rho_A(m,n)={\rm floor\/}(s), \hskip 30 pt
s=(1+A)m+\frac{1}{2q}.
\end{equation}
Then ${\cal G\/}_+$ determines the outward
pointing edge at $(m,n)$ if and only if
$\rho_A(m,n)$ is even.  This result is
implicit in the development given in
[{\bf S0\/}]. 

\subsection{The Fundamental Polytopes}
\label{graphPET}

Here are the $14$ fundamental polytopes for
the $(+)$ partition. Note that the listing
for the last one is spread over two lines.
\newline

\lefteqn{\left[\matrix{0 \cr -1 \cr 0 \cr 0}\right]
\hskip 1 pt\left[\matrix{0 \cr 0 \cr 0 \cr +1}\right]
\hskip 1 pt\left[\matrix{0 \cr -1 \cr +1 \cr +1}\right]
\hskip 1 pt\left[\matrix{0 \cr 0 \cr +1 \cr +1}\right]
\hskip 1 pt\left[\matrix{+1 \cr 0 \cr +1 \cr +1}\right]
\hskip 5 pt (0,+1)}
\vspace{4mm}

\lefteqn{\left[\matrix{0 \cr -1 \cr 0 \cr +1}\right]
\hskip 1 pt\left[\matrix{+1 \cr -1 \cr 0 \cr 0}\right]
\hskip 1 pt\left[\matrix{+1 \cr -1 \cr 0 \cr +1}\right]
\hskip 1 pt\left[\matrix{+1 \cr 0 \cr 0 \cr +1}\right]
\hskip 1 pt\left[\matrix{+1 \cr -1 \cr +1 \cr +1}\right]
\hskip 5 pt (0,+1)}
\vspace{4mm}

\lefteqn{\left[\matrix{0 \cr 0 \cr +1 \cr 0}\right]
\hskip 1 pt\left[\matrix{0 \cr +1 \cr +1 \cr +1}\right]
\hskip 1 pt\left[\matrix{0 \cr 0 \cr 2 \cr +1}\right]
\hskip 1 pt\left[\matrix{0 \cr +1 \cr 2 \cr +1}\right]
\hskip 1 pt\left[\matrix{+1 \cr +1 \cr 2 \cr +1}\right]
\hskip 5 pt (-1,0)}
\vspace{4mm}

\lefteqn{\left[\matrix{0 \cr 0 \cr 0 \cr 0}\right]
\hskip 1 pt\left[\matrix{0 \cr 0 \cr +1 \cr +1}\right]
\hskip 1 pt\left[\matrix{+1 \cr 0 \cr +1 \cr +1}\right]
\hskip 1 pt\left[\matrix{+1 \cr +1 \cr +1 \cr +1}\right]
\hskip 1 pt\left[\matrix{+1 \cr 0 \cr 2 \cr +1}\right]
\hskip 5 pt (-1,0)}
\vspace{4mm}

\lefteqn{\left[\matrix{0 \cr -1 \cr +1 \cr +1}\right]
\hskip 1 pt\left[\matrix{+1 \cr -1 \cr +1 \cr 0}\right]
\hskip 1 pt\left[\matrix{+1 \cr -1 \cr +1 \cr +1}\right]
\hskip 1 pt\left[\matrix{+1 \cr 0 \cr +1 \cr +1}\right]
\hskip 1 pt\left[\matrix{+1 \cr -1 \cr 2 \cr +1}\right]
\hskip 5 pt (+1,0)}
\vspace{4mm}

\lefteqn{\left[\matrix{0 \cr 0 \cr 0 \cr 0}\right]
\hskip 1 pt\left[\matrix{+1 \cr 0 \cr 0 \cr 0}\right]
\hskip 1 pt\left[\matrix{+1 \cr 0 \cr +1 \cr 0}\right]
\hskip 1 pt\left[\matrix{+1 \cr 0 \cr +1 \cr +1}\right]
\hskip 1 pt\left[\matrix{+1 \cr +1 \cr +1 \cr +1}\right]
\hskip 1 pt\left[\matrix{+1 \cr 0 \cr 2 \cr +1}\right]
\hskip 5 pt (-1,-1)}
\vspace{4mm}

\lefteqn{\left[\matrix{0 \cr -1 \cr 0 \cr 0}\right]
\hskip 1 pt\left[\matrix{0 \cr 0 \cr 0 \cr 0}\right]
\hskip 1 pt\left[\matrix{+1 \cr 0 \cr 0 \cr 0}\right]
\hskip 1 pt\left[\matrix{+1 \cr +1 \cr 0 \cr +1}\right]
\hskip 1 pt\left[\matrix{+1 \cr 0 \cr +1 \cr +1}\right]
\hskip 1 pt\left[\matrix{+1 \cr +1 \cr +1 \cr +1}\right]
\hskip 5 pt (-1,0)}
\vspace{4mm}

\lefteqn{\left[\matrix{0 \cr -1 \cr +1 \cr 0}\right]
\hskip 1 pt\left[\matrix{0 \cr 0 \cr +1 \cr 0}\right]
\hskip 1 pt\left[\matrix{0 \cr 0 \cr +1 \cr +1}\right]
\hskip 1 pt\left[\matrix{+1 \cr 0 \cr +1 \cr 0}\right]
\hskip 1 pt\left[\matrix{0 \cr -1 \cr 2 \cr +1}\right]
\hskip 1 pt\left[\matrix{0 \cr 0 \cr 2 \cr +1}\right]
\hskip 1 pt\left[\matrix{+1 \cr 0 \cr 2 \cr +1}\right]
\hskip 5 pt (+1,0)}
\vspace{4mm}

\lefteqn{\left[\matrix{0 \cr -1 \cr 0 \cr 0}\right]
\hskip 1 pt\left[\matrix{0 \cr 0 \cr 0 \cr 0}\right]
\hskip 1 pt\left[\matrix{0 \cr 0 \cr 0 \cr +1}\right]
\hskip 1 pt\left[\matrix{0 \cr +1 \cr 0 \cr +1}\right]
\hskip 1 pt\left[\matrix{+1 \cr +1 \cr 0 \cr +1}\right]
\hskip 1 pt\left[\matrix{0 \cr 0 \cr +1 \cr +1}\right]
\hskip 1 pt\left[\matrix{+1 \cr 0 \cr +1 \cr +1}\right]
\hskip 1 pt\left[\matrix{+1 \cr +1 \cr +1 \cr +1}\right]
\hskip 5 pt (-1,+1)}
\vspace{4mm}

\lefteqn{\left[\matrix{0 \cr 0 \cr 0 \cr 0}\right]
\hskip 1 pt\left[\matrix{0 \cr +1 \cr 0 \cr +1}\right]
\hskip 1 pt\left[\matrix{0 \cr -1 \cr +1 \cr 0}\right]
\hskip 1 pt\left[\matrix{0 \cr 0 \cr +1 \cr 0}\right]
\hskip 1 pt\left[\matrix{0 \cr 0 \cr +1 \cr +1}\right]
\hskip 1 pt\left[\matrix{+1 \cr 0 \cr +1 \cr 0}\right]
\hskip 1 pt\left[\matrix{0 \cr +1 \cr +1 \cr +1}\right]
\hskip 1 pt\left[\matrix{+1 \cr +1 \cr +1 \cr +1}\right]
\hskip 5 pt (0,+1)}
\vspace{4mm}

\lefteqn{\left[\matrix{0 \cr -1 \cr 0 \cr 0}\right]
\hskip 1 pt\left[\matrix{+1 \cr -1 \cr 0 \cr 0}\right]
\hskip 1 pt\left[\matrix{0 \cr 0 \cr 0 \cr +1}\right]
\hskip 1 pt\left[\matrix{+1 \cr 0 \cr 0 \cr 0}\right]
\hskip 1 pt\left[\matrix{+1 \cr 0 \cr 0 \cr +1}\right]
\hskip 1 pt\left[\matrix{+1 \cr +1 \cr 0 \cr +1}\right]
\hskip 1 pt\left[\matrix{+1 \cr -1 \cr +1 \cr 0}\right]
\hskip 1 pt\left[\matrix{+1 \cr 0 \cr +1 \cr +1}\right]
\hskip 5 pt (0,+1)}
\vspace{4mm}

\lefteqn{\left[\matrix{0 \cr 0 \cr +1 \cr 0}\right]
\hskip 1 pt\left[\matrix{0 \cr 0 \cr +1 \cr +1}\right]
\hskip 1 pt\left[\matrix{+1 \cr 0 \cr +1 \cr 0}\right]
\hskip 1 pt\left[\matrix{0 \cr +1 \cr +1 \cr +1}\right]
\hskip 1 pt\left[\matrix{+1 \cr +1 \cr +1 \cr +1}\right]
\hskip 1 pt\left[\matrix{0 \cr 0 \cr 2 \cr +1}\right]
\hskip 1 pt\left[\matrix{+1 \cr 0 \cr 2 \cr +1}\right]
\hskip 1 pt\left[\matrix{+1 \cr +1 \cr 2 \cr +1}\right]
\hskip 5 pt (0,0)}
\vspace{4mm}

\lefteqn{\left[\matrix{0 \cr -1 \cr 0 \cr 0}\right]
\hskip 1 pt\left[\matrix{0 \cr -1 \cr 0 \cr +1}\right]
\hskip 1 pt\left[\matrix{+1 \cr -1 \cr 0 \cr 0}\right]
\hskip 1 pt\left[\matrix{0 \cr 0 \cr 0 \cr +1}\right]
\hskip 1 pt\left[\matrix{+1 \cr 0 \cr 0 \cr +1}\right]
\hskip 1 pt\left[\matrix{0 \cr -1 \cr +1 \cr +1}\right]
\hskip 1 pt\left[\matrix{+1 \cr -1 \cr +1 \cr 0}\right]
\hskip 1 pt\left[\matrix{+1 \cr -1 \cr +1 \cr +1}\right]
\hskip 1 pt\left[\matrix{+1 \cr 0 \cr +1 \cr +1}\right]
\hskip 5 pt (+1,+1)}
\vspace{4mm}

\lefteqn{\left[\matrix{0 \cr -1 \cr 0 \cr 0}\right]
\hskip 1 pt\left[\matrix{0 \cr 0 \cr 0 \cr 0}\right]
\hskip 1 pt\left[\matrix{+1 \cr 0 \cr 0 \cr 0}\right]
\hskip 1 pt\left[\matrix{0 \cr -1 \cr +1 \cr 0}\right]
\hskip 1 pt\left[\matrix{0 \cr -1 \cr +1 \cr +1}\right]
\hskip 1 pt\left[\matrix{+1 \cr -1 \cr +1 \cr 0}\right]}
\vspace{2mm}
\lefteqn{\left[\matrix{0 \cr 0 \cr +1 \cr +1}\right]
\hskip 1 pt\left[\matrix{+1 \cr 0 \cr +1 \cr 0}\right]
\hskip 1 pt\left[\matrix{+1 \cr 0 \cr +1 \cr +1}\right]
\hskip 1 pt\left[\matrix{0 \cr -1 \cr 2 \cr +1}\right]
\hskip 1 pt\left[\matrix{+1 \cr -1 \cr 2 \cr +1}\right]
\hskip 1 pt\left[\matrix{+1 \cr 0 \cr 2 \cr +1}\right]
\hskip 5 pt (0,0)}
\vspace{4mm}

\newpage

\section{The Graph Reconstruction Formula}
\label{reconstruct0}

The purpose of this chapter is to present
a formula which we call the
{\it Reconstruction formula\/}, which
relates the geometry of the
graph grid to the graph PET classifying map.
Combining this formula with the
Master Picture Theorem, we will 
prove Statement 1 of the Pixellation Theorem.

\subsection{Main Result}

Recall that $\widehat X=\R^3 \times [0,1]$.
For each parameter $A$, the affine
group acting on $\widehat X$ in connection
with the graph PET acts as an abelian
group of translations.  Precisely, this
group is $\Lambda(\Z^3)$, where

\begin{equation}
\label{graphmatrix}
\Lambda=\left[\matrix{
 1 & 0& 0 \cr
-1 & 1+A & 0 \cr
-1 & 1-A & 1+A}\right] \Z^3
\end{equation}

Given a point $\xi \in G$, let
$[\xi] \in \R^2/\Z^2$ denote the
equivalence class of $\xi$.
We now explain how we can use the map
$\Phi$ to determine $[\xi]$.

We introduce linear maps $\Theta_1,\Theta_2: \R^3 \to \R$:
\begin{equation}
\Theta_1(x,y,z)=x, \hskip 30 pt
\Theta_2(x,y,z)=\frac{y-Ax}{1+A}.
\end{equation}

Checking on the obvious basis for $\Lambda$, we see
that $\Theta_j(\Lambda) \subset \Z$.  Therefore
$\Theta_1$ and $\Theta_2$ both give well defined
maps from $\R^3/\Lambda$, the domain
for the graph compactification at $A$, 
to the torus $\R/\Z$.  
Hence, we have a locally affine map
\begin{equation}
\Theta: \R^3/\Lambda \to \R^2/\Z^2, \hskip 30 pt
L=[(\Theta_1,\Theta_2)].
\end{equation}

The goal of this section is to establish
the following result.
\begin{lemma}[Reconstruction]
\begin{equation}
\label{reconstruct}
[\xi]=\Theta \circ \Phi(\xi).
\end{equation}
\end{lemma}
We call this equation the
Reconstruction Formula.
Equation \ref{reconstruct}
allows us to get control over how the
arithmetic graph sits with respect to
$\Z^2$. We've already proven that
no point of $G$ lies on the
boundary of an integer unit square.
So, we can interpret both sides of
Equation \ref{reconstruct} as 
referring to points in the open integer
unit square $(0,1)^2$.
We fix a parameter $A$ throughout the proof.

\begin{lemma}
\label{easy}
Equation \ref{reconstruct} holds 
at $\xi=T(0,0)$.
\end{lemma}

\startproof
A direct calculation shows that
$$[\xi]=\bigg(\frac{1}{2q},\frac{q-p}{2q(q+p)}\bigg).$$
Next, we compute
$$\Phi(\xi)=\bigg(\frac{1}{2q},\frac{1}{2q},\frac{1}{2q},A\bigg).$$
A direct calculation then shows that
$$\Theta_2(\Phi(\xi))=\frac{(p/q)(1/2q)-(1/2q)}{1+(p/q)}=\frac{q-p}{2q(q+p)}.$$
Hence, the two sides of the equation agree.
\endproof

\begin{lemma}
If Equation \ref{reconstruct} holds
at $\xi$, it also holds at
$\xi  \pm  dT(0,1)$.
\end{lemma}

\startproof
We do the $(+)$ case. The $(-)$ case has the same proof.
Let $\xi'=\xi+dT(0,1)=\xi+(1,-P)$.
Here $P=2A/(1+A)$.
We compute
$$\Theta \circ \Phi(\xi')-\Theta \circ \Psi(\xi)=
\Theta(1,1,1)=\bigg(1,\frac{1-A}{1+A}\bigg)=
(1,-P+1).$$
The last expression is congruent to
$\zeta'-\zeta=(1,-P)$ mod $\Z^2$.
\endproof

\begin{lemma}
If Equation \ref{reconstruct} holds
at $\xi$, it also holds at
$\xi \pm dT(1,1)$.
\end{lemma}

\startproof
We do the $(+)$ case. The $(-)$ case has the same proof.
We have $dT(1,1)=(1+A,1+A)$.  Using the same
notation as in the previous argument, we have
$$\Theta \circ \Phi(\xi')-\Theta \circ \Psi(\xi)=
\Theta(1+A,1+A,1+A)=(1+A,1-A)=\zeta'-\zeta.$$
This completes the proof
\endproof

These three lemmas together show that
Equation \ref{reconstruct} holds on
all of $dT(\Z^2)$, as desired.

\subsection{Eliminating Most Double Crossings}

In this section we eliminate most of the
double crossings mentioned in \S \ref{doublecross}.
This section
does not use either the Reconstruction Formula or
the Master Picture Theorem.

Let $T$ be the canonical
affine transform. See \S \ref{CAT}.
Recall that
a {\it distinguished edge\/} is an edge
having vertices $v_1,v_2 \in T(\Z)$
such that $v_1-v_2=T(\zeta)$, where
$\zeta$ is one of the $8$ shortest nonzero
vectors in $\Z^2$.  We say that
the edge belongs to the family
${\cal F\/}(i,j)$ is $\zeta=\pm(i,j)$.

A {\it double crossing\/}
is a configuration of the kind shown in
Figure 4.4 or in Figure 8.1 below.
The configuration consists
of two adjacent unit integer squares
$\Sigma_1$ and $\Sigma_2$, graph
grid points $v_1 \in \Sigma_1$ and
$v_2 \in \Sigma_2$, and 
disjoint distinguished edges
$e_1$ and $e_2$ incident to
$v_1$ and $v_2$ respectively
which both cross $\Sigma_1 \cap \Sigma_2$.

\begin{lemma}
A bad configuration must have the following
structure. \begin{itemize}
\item $\Sigma_1 \cap \Sigma_2$ is horizontal.
\item The edge connecting $v_1$ to $v_2$ belongs
${\cal F\/}(1,0)$.  
\item At least one of the two edges
$e_1$ or $e_2$ belongs to ${\cal F\/}(-1,1)$.
\end{itemize}
\end{lemma}

\startproof
Let $f=v_2-v_1$.   By
Statement 5 of the Grid Geometry Lemma,
$f$ must be a distinsuighed edge.
 
We define $S_1$ to be the
set of $8$ distinguished edges incident to
$v_1$.  We think of these edges as vectors
pointing out of $v_1$.
This set has a natural cyclic order on it.
Likewise we define $S_2$.  

Let $\Sigma_{12}=\Sigma_1 \cap \Sigma_2$.
The edges of $S_1$ which intersect
$\Sigma_{12}$ closest to $f \cap \Sigma_{12}$ are obtained
by turning $f$ one click in $S_1$,
either clockwise or counterclockwise.
Likewise, the edges of $S_2$ which intersect
$\Sigma_{12}$ closest to $f \cap \Sigma_{12}$ are obtained
by turning $(-f)$ one click in $S_2$,
counter clockwise or counterclockwise.
The two turnings must be in the same direction,
because otherwise the resulting edges would intersect.

So, the shortest possible
distance between the two intersection
points occurs when $e_1$ and $e_2$ are parallel.
Statement 6 of the Grid Geometry Lemma now says
that $\Sigma_{12}$ is horizontal and the
parallel lines are of type $(-1,1)$.  Looking at
the proof of Statement 6 of the Grid
Geometry Theorem given in \S \ref{gg0}, we
see that this forced $f \in {\cal F\/}(1,0)$.
\endproof

\subsection{Eliminating the last Double Crossing}

Now we will use the Reconstruction Formula
and the Master Picture Theorem to rule out
the last kind of double crossing.  This proves
Statement 1 of the Pixellation Lemma.
We will suppose that the arithmetic graph has
a double crossing for some parameter $A$ and
derive a contradiction.  The reader is warned
in advance that our proof, at the end, is
just a computer assisted calculation. 

We will 
consider the case when the edge $e_1$ shown 
on the left side of
Figure 8.1 is in the family ${\cal F\/}(-1,1)$.
The other case, when $e_2 \in {\cal F\/}(-1,1)$,
has a similar proof.  Indeed, the second
case follows from the first case and
from the rotational symmetry of the
arithmetic graph.

\begin{center}
\resizebox{!}{4.3in}{\includegraphics{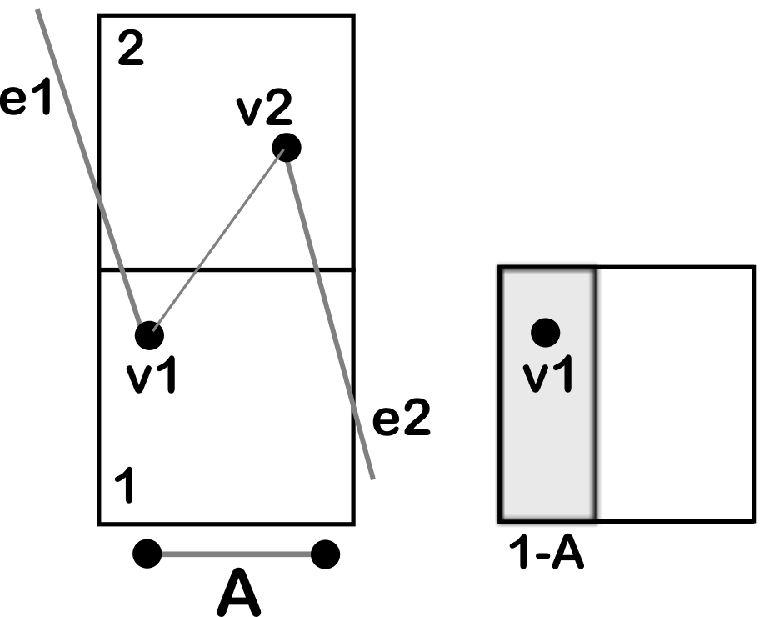}}
\newline
{\bf Figure 8.1:\/} A double crossing.
\end{center}

Let $G$ denote the graph grid and let
$\Phi: G \to \widehat X$ be the graph classifying
map.

\begin{lemma}
$\Phi(v_1)$ lies in the region of
$\widehat X \cap H$, where $H$ is the half
space given by the equation $x<1-A$.
\end{lemma}
 We know that
the edge joining $v_1$ to $v_2$ is in the family
${\cal F\/}(1,0)$.  
The first coordinate of $T(1,0)$ is $A$, so
the horizontal distance from $v_1$ to $v_2$ is
$A$. That means that $v_1$ is within $1-A$
units of the left edge of $\Sigma_1$, as
shown on the right half of Figure 8.1.

By the Reconstruction Formula, $\Phi(v_1)$ lies in
the region
\begin{equation}
R_A=X_A \cap \{x|\ x<1-A\}.
\end{equation}
$R_A$ lies in the region
\begin{equation}
H \cap \widehat X,
\end{equation}
where $H$ is the halfspace given
by the equation $x<1-A$.  
Recall that the coordinates on this space
are $(x,y,z,A)$, so $H$ is an defined
by inteer equations.  Also, $H$ is invariant
under the action of the graph
lattice $\Lambda$ defined in
\S \ref{mpt}.  
\endproof

Recall that we parametrize $\widehat X \subset \R^4$
with coordinates $(x,y,z,A)$. So, $H$ is a halfspace
defined by integer equations.  Moreover, $H$ is
invariant under the action of the graph
lattice defined in \S \ref{mpt}.

The way the graph classifying map works
is that we look at $\Phi(v_1)$ and
record the labels of the polytope
in each partition which contains this
point. 

\begin{lemma}
Assuming that a bad configuration exists, one
of the polytopes in one of the partitions
has the label $(1,-1)$ and intersects $H$.
\end{lemma}

\startproof
There are two polytopes of interest to us:
the one in each partition which contains
$\Phi(v_1)$.  These polytopes both
intersect $H$.  One of these
polytopes is responsible for the assignment
of the vector $e_1$ to $v_1$.  Call
this the {\it magic polytope\/}.
The magic polytope
Since $e_1$ is in the family
${\cal F\/}(-1,1)$, the label of the
magic polytope is either
$(-1,1)$ or $(1,-1)$.
We check that
$T(1,-1)$ is the one with positive
$y$ coordinate. Hence, the label
of the magic polytope is $(1,-1)$.
\endproof

Now for the computer assisted part of the proof.
We check that no polytope at all in the $(+)$ partition
has the label $(1,-1)$, and the only polytopes in
the $(-)$ partition intersect $H$. This is obvious
from the pictures on my computer program, and I
will describe in the last chapter the short linear
algebra computation which checks it.
This completes the proof.

\newpage

\section{The Hitset and the Intertwiner}

\subsection{The Hitset}
\label{hitset0}

We fix a paramater $A$.
Recall that a unit integer square
is {\it grid full\/} if it contains
a graph grid point, and otherwise
grid empty.  Recall that $\Phi_{\Pi}$ is
the plaid classifying map.
Let $G_{\Pi}$ denote the plaid grid - i.e.,
the centers of the unit integer squares.
Let $G^*_{\Pi} \subset G_{\Pi}$ denote
the set of centers of grid full squares.
In this section we state a result which
characterizes the image
\begin{equation}
\bigcup_{\zeta \in G^*_{\Pi}} \Phi_{\Pi}(\zeta).
\end{equation}

The domain of $\Phi_{\Pi}$ is the unit
cube 
\begin{equation}
X_{\Pi}=[-1,1]^3.
\end{equation}
We define the {\it hitset\/} to be the
subset of $X_{\Pi}$ having the form
\begin{equation}
X_{\Pi}^*=H \times [-1,1],
\end{equation}
where $H$ is the
octagon with vertices
$$
(-1,1),(-1+P,-1+P),(1-P,-1),(1,-1+P)$$
\begin{equation}
(1,1),(1-P,1-P),(-1+P,1),(-1,1-P).
\end{equation}
The vertices are listed in cyclic order.
The octagon $H$ has a kind of zig-zag shape.
Figute 9.1 shows the picture for several
parameters.

In this chapter we will prove the following
theorem.

\begin{theorem}[Hitset]
For any parameter $A$, we have
\begin{equation}
\label{HITSET}
\bigcup_{\zeta \in G_{\Pi}^*} \Phi_{\Pi}(\zeta) \subset X_{\Pi}^*.
\end{equation}
\end{theorem}

The Hitset Theorem is sharp in the following sense. 
All the objects make sense at irrational parameters
as well as rational parameters, and for irrational
parameters, the set on the left is dense in the
set on the right.  We will not prove this, because
we do not need the result, but a proof would not
be so difficult given everything else we prove in
this chapter.

In Figure 9.1 below we show
the polygon $H$ for the parameters
$k/5$ for $k=1,2,3,4$. (The picture
makes sense for all parameters, and
not just even rational ones.)
Actually, we show not just
$H$ but also the image of $H$
under the action of the
lattice $\Lambda_{\Pi}$.  We
think of this lattice as acting on
the $xy$ plane just by forgetting
about the third coordinate.
The lattice is generated by the
vectors $(2,P)$ and $(0,2)$.

Thanks to the product structure of
$X_{\Pi}$ and $X_{\Pi}^*$, the
planar pictures we show capture
all the information.
The various images of $H$ fit together
to form infinite bands which look like
zigzags.  As $A \to 0$ the zigzags come
together and fill up the plane.  

\begin{center}
\resizebox{!}{2.5in}{\includegraphics{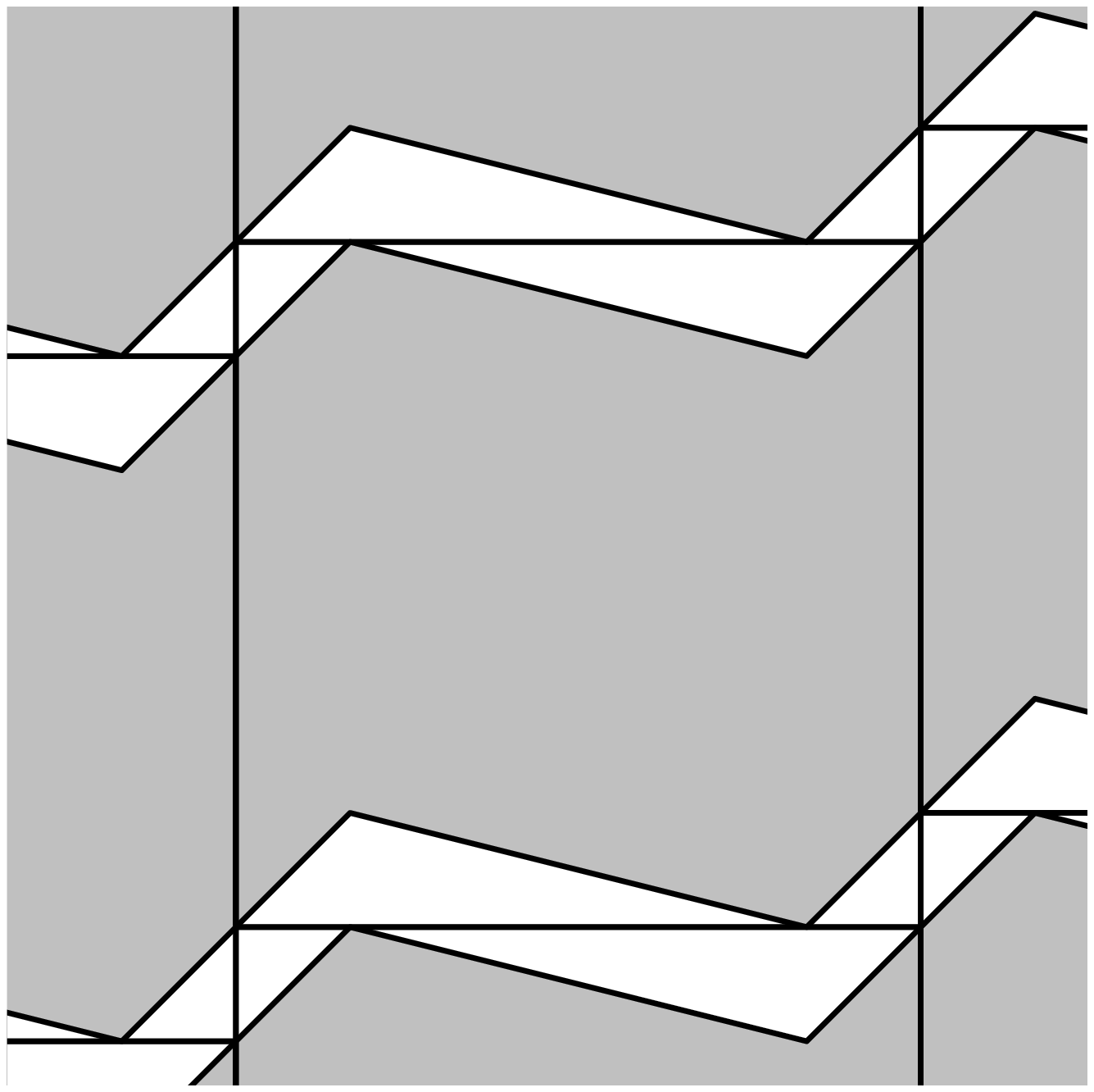}}
\hskip 10 pt
\resizebox{!}{2.5in}{\includegraphics{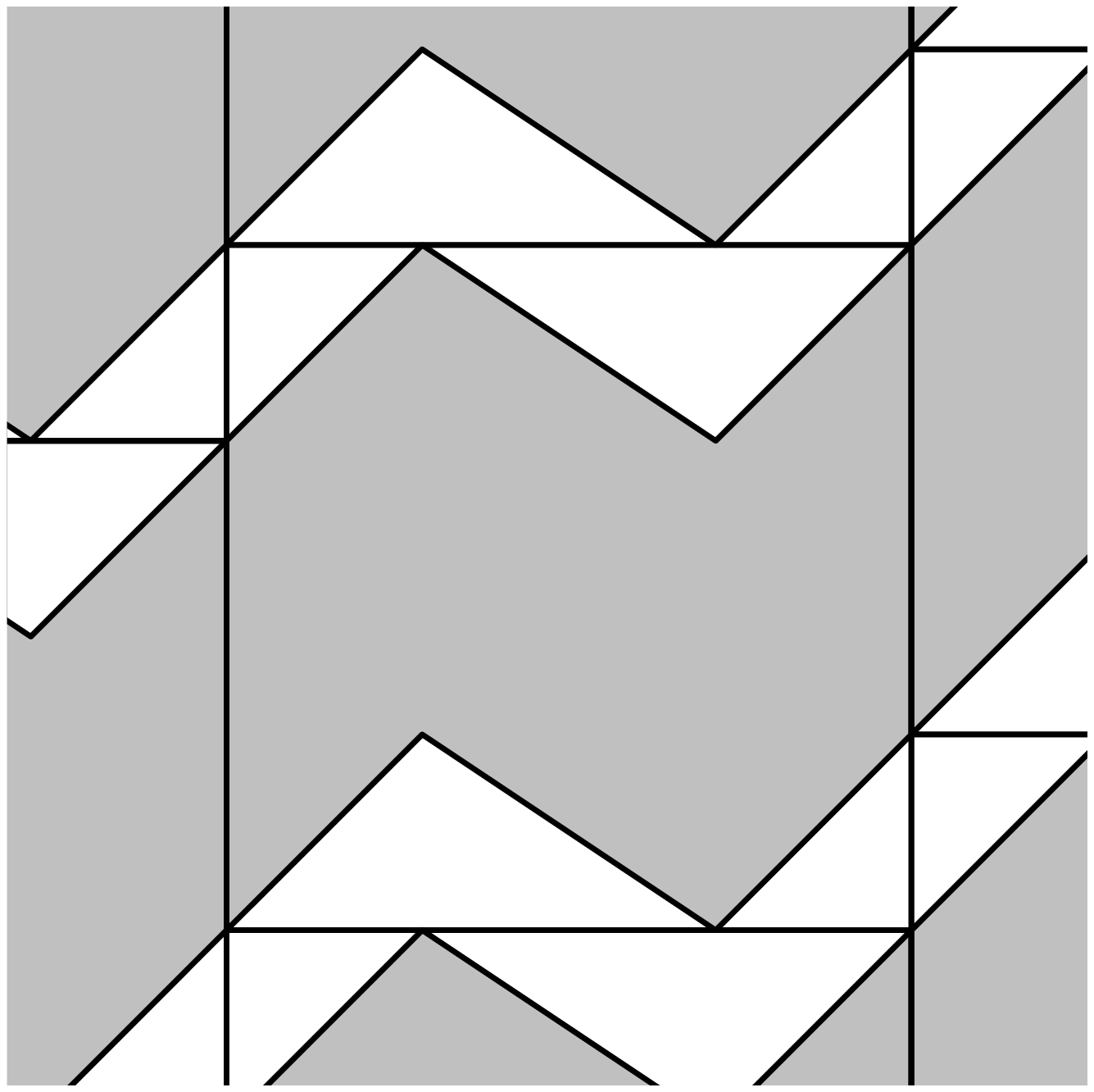}}
\vskip 10 pt
\resizebox{!}{2.5in}{\includegraphics{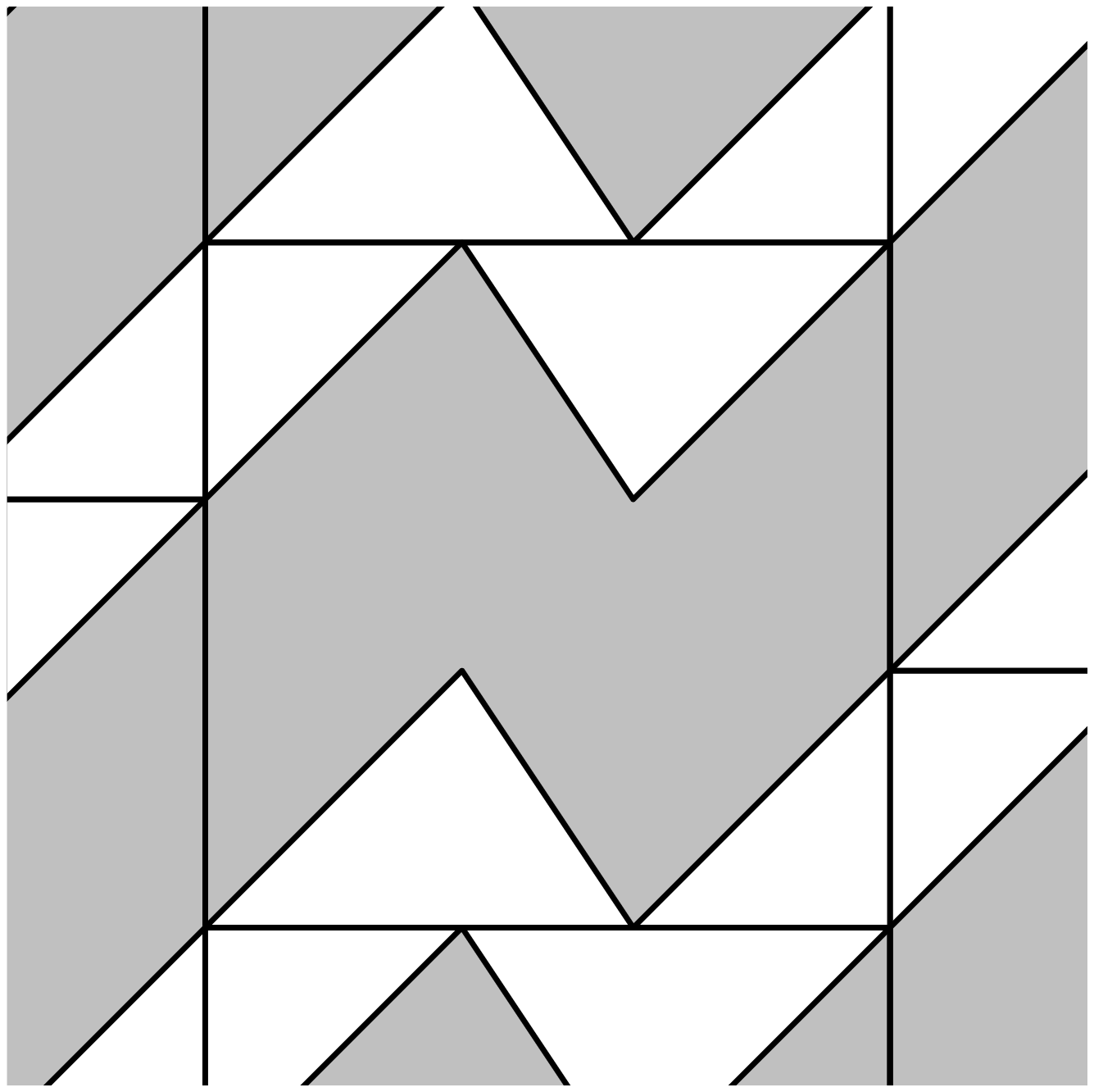}}
\hskip 10 pt
\resizebox{!}{2.5in}{\includegraphics{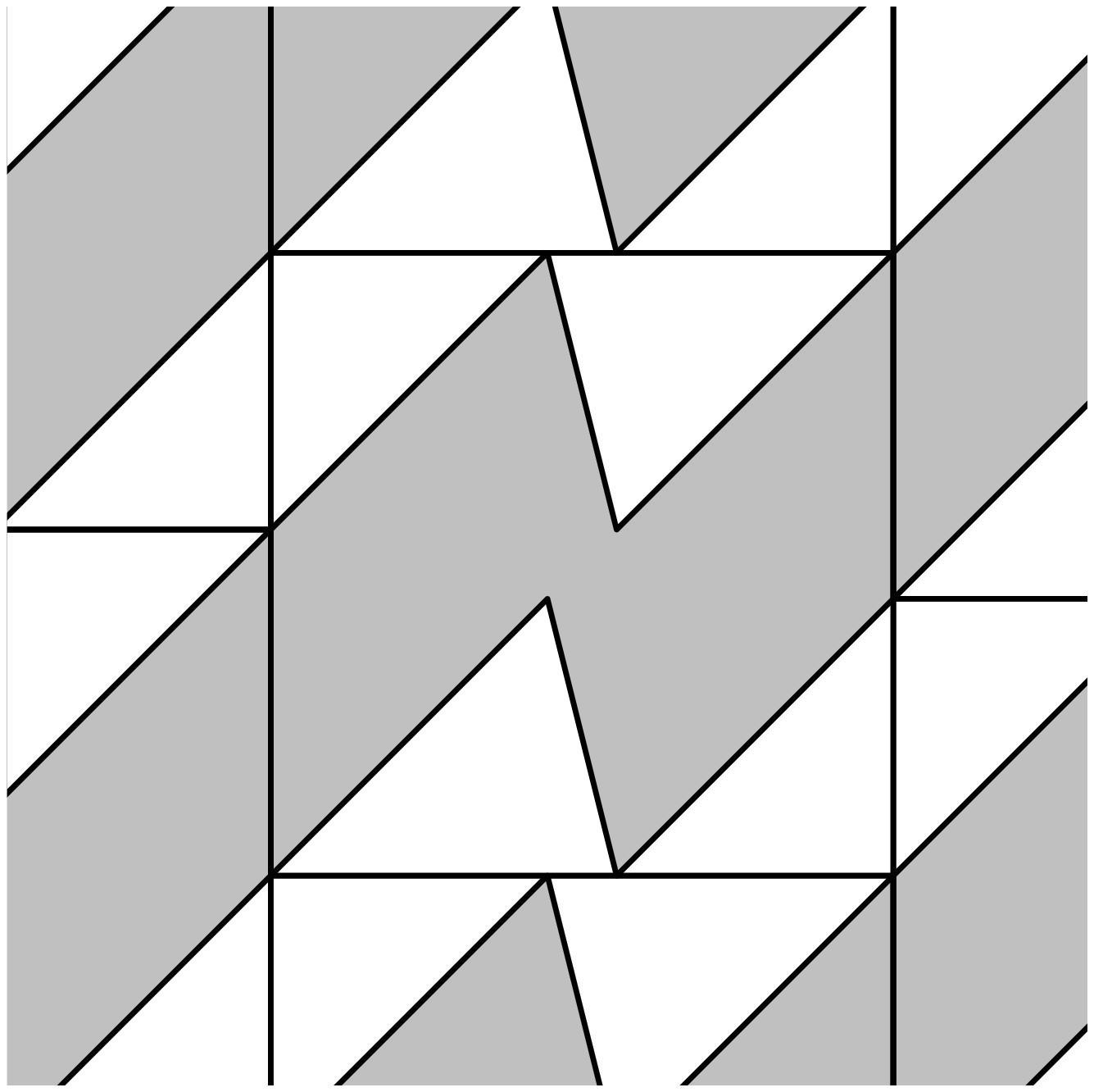}}
\newline
{\bf Figure 9.1:\/} The orbit $\Lambda_{\Pi}(H)$ for
parameters $A=k/5$, $k=1,2,3,4$.
\end{center}

\subsection{The Projective Intertwiner}

The symbol $\Psi$ will denote the map
which we call {\it the projective intertwiner\/}.
We will define
$\Psi$ after
we specify its domain and range.
Recall that our total space is
$\widehat X=\R^3 \times [0,1]$.
\newline
\newline
\noindent
{\bf Domain:\/}
The domain of $\Psi$ is $X_{\Pi}$, though
our theorem really only concerns the image
of $\Psi$ on the hitset $X_{\Pi}^*$.
The set $X_{\Pi}$ is the fundamental
domain for the lattice
$\Lambda_1$ acting on $\widehat X$.
See \S \ref{plaidspace}. 

We write
\begin{equation}
X_{\Pi}=X_{\Pi,-} \cup X_{\Pi,+},
\end{equation}
Where $X_{\Pi,+}$ consists of those points
$(x,y,z,P)$ where $x \leq y$ and
$X_{\Pi,-}$ cnsists of those points there
$y \leq x$.  This is a partition of
$X_{\Pi}$ into two isometric halves.
\newline
\newline
{\bf Range:\/}
The range of $\Psi$ is
\begin{equation}
X_{\Gamma}=\widehat X/\Lambda,
\end{equation}
where $\Lambda$
is the graph lattice
defined in \S \ref{mpt}.
\newline
\newline
{\bf The Map:\/}
Now we define the map
$$
\Psi: X_{\Pi} \to X_{\Gamma}.
$$
For $(x,y,z,P) \in X_{\Pi,\pm}$, we define
\begin{equation}
\label{MAP}
\Psi(x,y,z,P)=
\left[\frac{1}{2-P}\bigg(x-y,-y-1,z+P+1,P\bigg)
\pm (1,0,0,0)\right]_{\Lambda}
\end{equation}
That is, we add or subtract $1$ depending on which
half of the partition our point lies in, and
we take the result mod $\Lambda$.  In
the next section we check that $\Psi$ is
well defined even in boundary cases.
The map $\Psi$ is a piecewise defined
integral projective transformation.  If we
hold $P$ fixed and restrict $\Psi$ to a slice,
then $\Psi$ is an affine transformation.
We call $\Psi$ the {\it projective intertwiner\/}.
\newline
\newline
\noindent
{\bf The Main Result:\/}
Let $G_{\Gamma}$ denote the graph grid.
Let $G_{\Pi}$ denote the plaid grid.
Each point $\zeta_{\Gamma} \in G_{\Gamma}$ lies in the
interior of a unique integer square.
We let $\zeta_{\Pi} \in G_{\Pi}$ denote the
center of this square.

Let $\Phi_{\Pi}: G_{\Pi} \to X_1$ denote the
plaid classifying map. Recall that, at the
parameter $A$, the range of $\Phi_{\Pi}$ is
contained in the slice $\R^3 \times \{P\}$,
where $P=2A/(1+A)$.
Let $\Phi_{\Gamma}: G_{\Gamma} \to \widehat X/\Lambda$
be the graph classifying map.  

\begin{theorem}[Intertwining]
The following holds for every even rational paramater:
\begin{equation}
\label{INTER}
\Phi_{\Gamma}(\zeta_{\Gamma})=\Psi \circ \Phi_{\Pi}(\zeta_{\Pi})
\hskip 30 pt
\forall \zeta_{\Gamma} \in G_{\Gamma}.
\end{equation}
\end{theorem}

\noindent
{\bf Remarks:\/} 
\newline
(i) I checked Equation \ref{INTER} computationally for all relevant
points and all parameter $p/q$ with $q<30$.  This check
is not meant as a substitute for a rigorous proof, but it
is nice to know.  I didn't check the Hitset Theorem
as systematically, but my computer program plots
the left and right hand sides of Equation
\ref{HITSET}, and one can see that it always works.
\newline
(ii) Notice in Equation \ref{MAP} that the fourth
coordinate on the right hand side is $A$, because
$A=P/(2-P)$.  Thus, $\Psi$ maps the relevant
slices to each other.
\newline
(iii) It might be nicer if there were a
global projective transformation from $\widehat X$ to
itself which works in place of our piecewise
projective map $\Psi$.  However, we have
$\Psi(T_X(V)-V)=(2,0,0,0)$.
The vector on the right does not have for
them $\lambda(W)-W$ for any transformation
$\lambda \in \Lambda_1$, the plaid lattice.
This situation makes the existence of a global
projective intertwiner impossible. On the
other hand, below in \S \ref{domain}
we will modify the domain
of $\Psi$ so that $\Psi$ is projective
throughout the interior of the new domain. 
When we look at the action of
$\Psi$ on the new domain, we will see
the canonical nature of $\Psi$.  See
Lemma \ref{CANON} for instance.
\newline
(iv) When it comes time to
prove the Pixellation Theorem, we shall be interested
in the action of $\Psi$ on the polytopes of the
plaid triple partition.  At first, it looks like
we might have trouble, due to the piecewise nature
of $\Psi$.  However, it turns out that every
plaid triple polytope is contained in one of the
two pieces of the partition of $X_{\Pi}$. 
\newline
(v) Our proof will show that a suitable formulation
of the Intertwining Theorem holds for all
parameters, and not just even rational ones.
Indeed, it basically follows from continuity.
\newline
(vi) The reader who is keen to see how the
proof of the Pixellation Theorem works might want
to take the Hitset Theorem and the Intertwining
Theorem for granted on the first reading.  Our
proof of the Pixellation Theorem only uses
the truth of these statements, and not any
theory developed during their proof. \S \ref{sample0}
gives a good illustration of how we use the
Intertwining Theorem.

\subsection{Well Definedness}

Since $\Psi$ is only piecewise defined,
we have to worry about the cases when
there are two competing definitions
for $\Psi$.  In our discussion of
this matter, the symbol $(*)$ stands
for a coordinate value that we don't
care about.  We fix some even rational
parameter for the discussion.

There are three issues. One issue is that perhaps
$$\Phi_{\Pi}(\zeta)=(-1,*,*).$$
In this case, the two points
$$\phi_+=\Phi_{\Pi}(\zeta), \hskip 30 pt
\phi_-=\Phi_{\Pi}(\zeta)+(2,P,P)$$
are equally good representatives.
Here $\phi_{\pm} \in X_{\Pi,\pm}$.
An easy calculation shows that
$$
\Psi(\phi_-)-\Psi(\phi_+)=
\frac{1}{2-P}\bigg(2-P,-P,PP,-P,2-P\bigg)-(2,0,0)=$$
\begin{equation}
(-1,-A,A) \in \Lambda_1.
\end{equation}
Hence, either representative gives rise to
the same point in the range.

The second issue seems more serious, in
view of Remark (v) above.  It might happen that
$$\phi=\Phi_{\Pi}(\zeta)=(v,v,*).$$
That is, $\phi$ lies on the boundary of
both pieces of the partition of $X_{\Pi}$.
Let's check that this situation cannot
arise for the images of
points of $G_{\Pi}$.  Such points
have the form $(x,y)$ where $x$ and $y$
are half integers.  We compute that
\begin{equation}
\label{imp1}
\Phi_{\Pi}(x,y)=(2Px+2y,2Px,*)+(2m,Pm,*),
\end{equation}
for a suitable integer $m$.  We also observe
that $2y$ is odd.  Equation \ref{imp1} leads to
$$P=\frac{2m+2y}{m}.$$
This is impossible, because $P=2p/(p+q)$,
and $2m+2y$ is an odd integer. 

The third issue is that $\Phi_{\Pi}(\zeta)=(*,u,v)$
with either $u=\pm 1$ or $v=\pm 1$ or both.
The case $u=\pm 1$ is impossible for similar reasons
that we have just discussed.  When $v=\pm 1$ it has
replacing $v$ by $v \mp 2$ has no effect on the
Intertwining formula.

\subsection{Strategy of the Proof}

The rest of the chapter is devoted to proving
the Hitset Theorem and the Intertwining
Theorem.  We prove these two results at
the same time because they are closely
related to each other.   Here are the
steps of the proof. We fix a parameter $A=p/q$.

\begin{itemize}

\item We prove the Intertwining Theorem
for points of the form
\begin{equation}
\label{diag}
\zeta_n=(n+1/2)(1+A,1-A), \hskip 30 pt n=0,1,2,....
\end{equation}
These points all belong to
$G_{\Gamma}$ because $\zeta_0$ is the
anchor point and $\zeta_n-\zeta_0=ndT(1,1)$.
Here $T$ is the canonical affine transformation.
We call the points in
Equation \ref{diag} the {\it diagonal points\/}.

\item We prove the Hitset Theorem
for the points in Equation \ref{diag}.

\item We prove the following induction step:
Suppose that the Hitset Theorem is true
for some graph grid point $\zeta$. Then it
is also true for $\zeta+dT(0,1)$.
We call this {\it Hitset Induction\/}.

\item We prove the following induction step:
Suppose that the Intertwining Theorem is true
for some graph grid point $\zeta$. Then it
is also true for $\zeta+dT(0,1)$. We
call this {\it Intertwiner Induction\/}.

\end{itemize}

Let $L$ denote the lattice of symmetries 
generated by the vectors $(\omega^2,0)$ and $(0,\omega)$.
Here $\omega=p+q$, as usual.
We have already shown that both the plaid model
and the arithmetic graph are invariant under $L$.
We say that an $L$-{\it orbit\/} is an orbit
of $L$ acting on the graph grid.
By symmetry, it suffices to
prove our two theorems on a set which contains
at least one point of every $L$-orbit

The four steps above combine to
to prove that the two theorems hold true on sets
of the form
$dT(B)$ where $B$ is a ball in $\Z^n$ of arbitrarily
large radius.  For sufficiently large $B$, the
set $dT(B)$ intersects every $L$ orbit.
\newline
\newline
{\bf Remark:\/}
It would have been nice if the
set of points $\{\zeta_n\}$ intersected
every $L$-orbit. If this was true,
we would not need the induction part of the proof.
Likewise, it would
have been nice if, starting from
a single point (such as the anchor), we
could reach every $L$-orbit just
by repeatedly adding $dT(0,1)$. If
this was true, we would not need the
first two steps of our proof.
Alas, neither of these step-saving
situations is true.

\subsection{The Intertwining Theorem on the Diagonal}

In this section we prove the Intertwining Theorem
for points of the form
\begin{equation}
\zeta_n=\bigg(n+\frac{1}{2}\bigg)(1+A,1-A),
\hskip 30 pt n=0,1,2,...
\end{equation}
These points all belong to $G_{\Gamma}$, because
$\zeta_n=\zeta_0+dT(n,n)$ and $\zeta_0$ is
the anchor point.  

We fix the value of $n$.  It is convenient to 
consider the sub-intervals
\begin{equation}
R_{n,k}=\bigg(\frac{2k}{2n+1},\frac{2k+1}{2n+1}\bigg),
\hskip 15 pt k=0,...,(n-1).
\end{equation}
\begin{equation}
L_{n,k}=
\bigg(\frac{2k-1}{2n+1},\frac{2k}{2n+1}\bigg),
\hskip 15 pt k=1,...,(n-1).
\end{equation}
We ignore the endpoints of these intervals;
the boundary cases, when relevant, follow
from continuity.

Let $I=(A_0,A_1)$ be one of the intervals of interest.
We will sometimes use the following trick
to bound certain numerical quantities that
depend on $A \in I$.  When the quantity is
monotone, we get the bounds by evaluating the
expression on the boundary values 
$A=A_0$ and $A=A_1$.  We call this 
method {\it the boundary trick\/}.

When $A \in L_{n,k}$ (respectively
$R \in L_{n,k}$, the point
$\zeta_n$ lies in the left (respectively right)
half of the square
with center
\begin{equation}
\zeta_{n,k}=\bigg(n+\frac{1}{2},n+\frac{1}{2}\bigg)+(k,-k).
\end{equation}

First we consider the case when $A \in L_{n,k}$.
We have
\begin{equation}
\label{hard1}
\Phi_{\Pi}(\zeta_{n,k})\equiv
(P(2n+2k+1)+(2n-2k+1),P(2n+2k+1),2P(2n+1)).
\end{equation}
The symbol $\equiv$ means that we still need
to reduce mod $\Lambda_1$ to get a vector
in the fundamental domain. 

The boundary trick tells us that the first
coordinate in Equation \ref{hard1} lies in
\begin{equation}
\label{hard2}
\bigg(1+2k+2n-\frac{1+2n}{k+n},1+2k+2n\bigg)
\end{equation}
So, we subtract off the lattice vector
$(n+k)(2,P,P)$.  This gives
\begin{equation}
\label{hard3}
\Phi_{\Pi}(\zeta_{n,k}) \equiv
(P(2n+2k+1)-4k+1,P(1+k+n),2P(2-k+3n)).
\end{equation}
The boundary trick tells us that
the second coordinate in Equation
\ref{hard3} lies in the interval $(-1,1)+2k$. 
Subtracting $(0,2k,0)$ from Equation \ref{hard3}, we get
a point in the fundamental domain:
\begin{equation}
\label{hard4}
\Phi_{\Pi}(\zeta_{n,k})=
(P(2n+2k+1)-4k+1,P(1+k+n)-2k,P(2-k+3n)+2\beta).
\end{equation}
Here $\beta$ is some integer whose value
we don't care about.

The interval check shows that
the first coordinate in Equation
\ref{hard4} is larger than the second
coordinate. (We apply the trick to the
difference of the coordinates.)  Hence
the point in Equation \ref{hard4} lies
in $X_{\Pi,-}$.  
We also have
\begin{equation}
\Phi_{\Gamma}(\zeta) \equiv
\bigg(n+\frac{1}{2}\bigg)(1+A,1+A,1+A) \hskip 5 pt
{\rm mod\/} \hskip 5 pt \Lambda.
\end{equation}
Here $\Lambda$ is the graph lattice.
(In all these equations we are leaving off
the fourth coordinate; we know this works out 
already.)

A calculation, which we do in Mathematica,
shows that
\begin{eqnarray}
\nonumber
\Psi(\Phi_{\Pi}(\zeta_{n,k})-\Phi_{\Gamma}(\zeta)= \cr
(-1-k-n)(1,-1,-1)+\cr(-2-2n)(0,1+A,1-A)+\cr
(1-k+\beta)(0,0,1+A).
\end{eqnarray}
In other words, we have written the difference
between the quantities as an integer combination
of vectors in the graph lattice $\Lambda$.
Hence, the two quantities are equal, as desired.

When $A \in R_{n,k}$ the calculation is very similar
and we will just describe the differences.  This
time the interval trick tells us to subtract off
the vector $(n+k+1)(2,P,P)$ from the point in
Equation \ref{hard1}.  Then, as in the other case,
we subtract off $(0,-2k,0)$.  The result is
\begin{equation}
\label{hard5}
\Phi_{\Pi}(\zeta_{n,k})=
(P(2k+2n+1)-4k-1,P(k+n)-2k,P(1-k+3n)+\beta)
\end{equation}
The interval trick shows that this point lies
in $X_{\Pi,+}$.  Now we compute
\begin{eqnarray}
\nonumber
\Psi(\Phi_{\Pi}(\zeta_{n,k})-\Phi_{\Gamma}(\zeta)= \cr
(-k-n)(1,-1,-1)+\cr(-1-2n)(0,1+A,1-A)+\cr
(1-k+\beta)(0,0,1+A).
\end{eqnarray}
This gives us the same desired conclusion as in the
first case.

\subsection{The Hitset Theorem on the Diagonal}

Now we prove the Hitset Theorem for the points
we considered in Step 1.  We first consider
the case when $A \in L_{n,k}$.  The triangle
$\Delta_L$ with vertices
\begin{equation}
(-1+P,-1+P), \hskip 30 pt
(1,-1+P), \hskip 30 pt (1,1)
\end{equation}
is the convex hull of $3$ of the vertices
and is contained in the polygon defining the
hitset.

\begin{center}
\resizebox{!}{1.5in}{\includegraphics{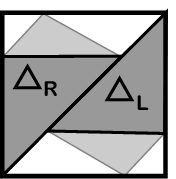}}
\newline
{\bf Figure 9.2:\/} The triangles
$\Delta_L$ and $\Delta_R$.
\end{center}

  Let
$\xi$ denote the point in the plane obtained
by taking the first two coordinates of
the point in Equation \ref{hard4}.  That is:

\begin{equation}
\label{hard6}
\xi=
(P(2n+2k+1)-4k+1,P(1+k+n)-2k).
\end{equation}

It suffices to prove that
$\xi \in \Delta_L$ for all $A \in L_{n,k}$.
One of the sides of $\Delta_L$ is the line
$x=1$, and certainly $\xi$ does not cross
this line; it corresponds to one of the
sides of the fundamental domain.
Another side of $\Delta_L$ is the line
$y=-1+P$.  The interval trick shows that
$\xi$ stays above this line.  Finally,
the other side of $\Delta_L$ is the line
$y=x$.  We already know that
$\Phi_{\Pi}(\zeta) \in X_{\Pi,-}$, and this
is the equivalent to the statement that
$\xi$ lies to the right of this line.
These three conditions together imply
that $\xi \in \Delta_L$.  My computer
program allows you to see a plot of
$\xi$ and $\Delta_L$ for all the relevant
parameters.

Now consider the case when $A \in R_{n,k}$.
This time we use the triangle $\Delta_R=-\Delta_L$,
obtained by negating all the coordinates
of the vertices of $\Delta$.  This time
we want to show that the point
\begin{equation}
\label{hard7}
(P(2k+2n+1)-4k-1,P(k+n)-2k)
\end{equation}
lies in $\Delta_R$.  The same arguments
as aboe shows that this point lies to
the right of the line $x=-1$, below the
line $y=1-P$, and to the left of the line
$y=x$. This does the job for us.

\subsection{Hitset Induction}

Here we will prove Hitset Induction
modulo what we call {\it the geometric claim\/}.
We will prove the geometric claim at the
end of the chapter.

We set
$\zeta_{\Gamma}=\zeta$ and
$\zeta_{\Gamma}'=\zeta+dT(0,1)$.
Similarly, we define
$\zeta_{\Pi}$ and $\zeta_{\Pi}'$.
Here $\zeta_{\Pi}$ is the center of
the unit integer square that contains
$\zeta_{\Gamma}$.   
Since we will be mentioning both the
plaid lattice and the graph lattice,
we let $\Lambda_{\Pi}$ denote the
plaid lattice and $\Lambda_{\Gamma}$ denote
the graph lattice.

We have
\begin{equation}
dT(0,1)=(1,-P).
\end{equation}
From this equation, we see that we have one of
two possibilities for $\zeta_{\Pi'}-\zeta_{\Pi}$.
This difference either equals $(1,0)$ or $(1,-1)$.
To be more precise, let
$G_{\Gamma}^{\rm hi\/}$ denote the union of
those choices of $\zeta_{\Gamma}=(x,y)$ such that
$x-{\rm floor\/}(x)>P$.
Define $\zeta_{\Pi}^{\rm lo\/}$ to be the
complementary set.
For even rational parameters, it never happens
that $y={\rm floor\/}(y)=P$, because then
$\zeta'_{\Gamma}$ would be lie on the boundary
of a unit integer square.

Let $H^{\rm lo\/}$ denote the parallelogram
with vertices
\begin{equation}
(1-3P,1-3P), \hskip 18 pt
(1-P,1-P), \hskip 18 pt
(-1+P,1)  \hskip 18 pt
(-1-P,1-2P).
\end{equation}
Let $H^{\rm hi\/}$ denote the polygon with vertices
\begin{equation}
(-1+P,-1+P), \hskip 18 pt
(1-P,-1), \hskip 18 pt
(3-3P,1-2P)  \hskip 18 pt
(1-P,1-P).
\end{equation}
These polygons are not subsets of $H$.  However, the
orbits $\Lambda_{\Pi}(H^{\rm hi\/})$ and
$\Lambda_{\Pi}(H^{\rm lo\/})$ together
give a partition of $\Lambda_{\Pi}(H)$.
The dark region is the {\it lo orbit\/} and
the light region is the {\it hi orbit\/}.

\begin{center}
\resizebox{!}{2in}{\includegraphics{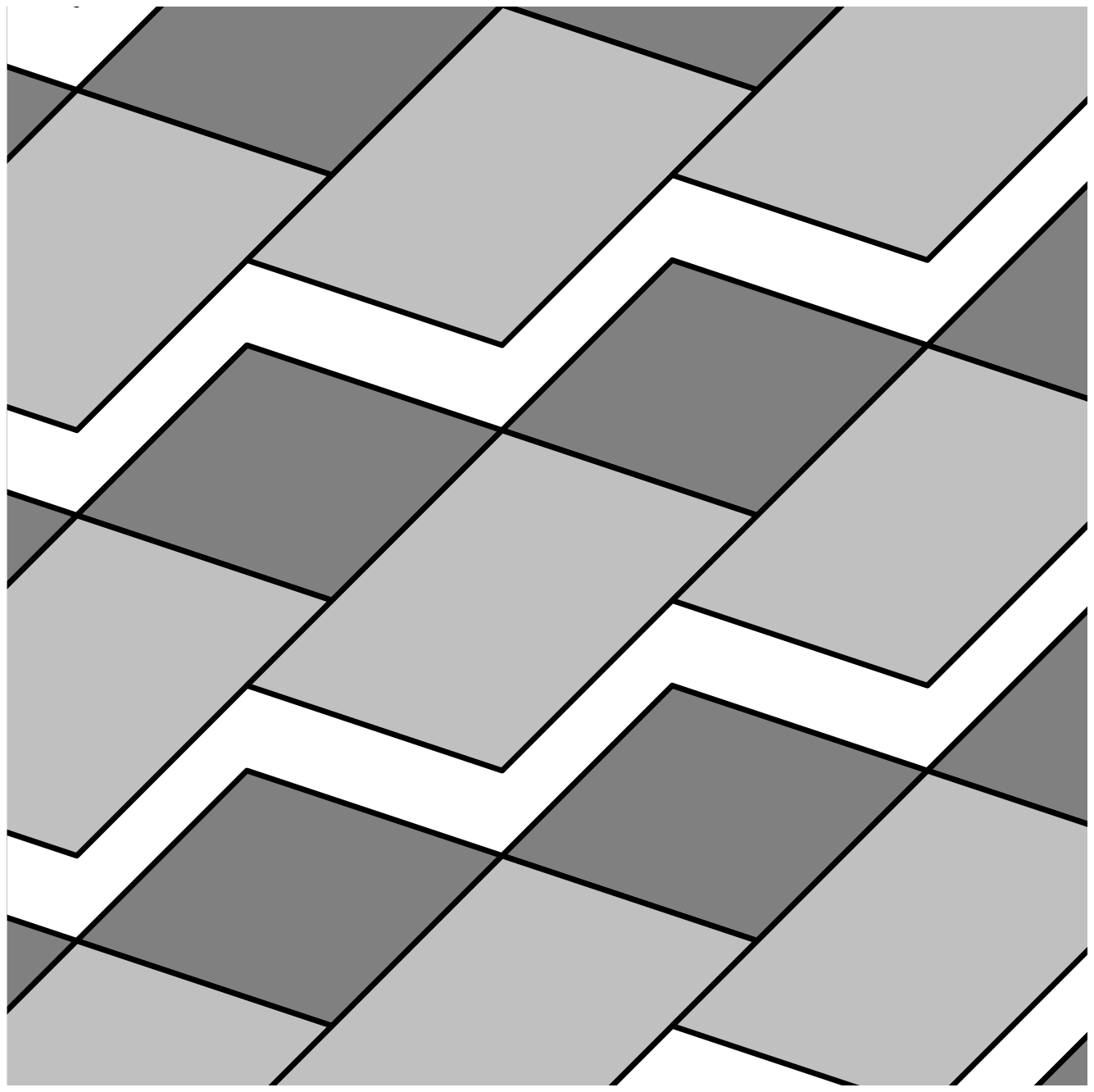}}
\hskip 50 pt
\resizebox{!}{2in}{\includegraphics{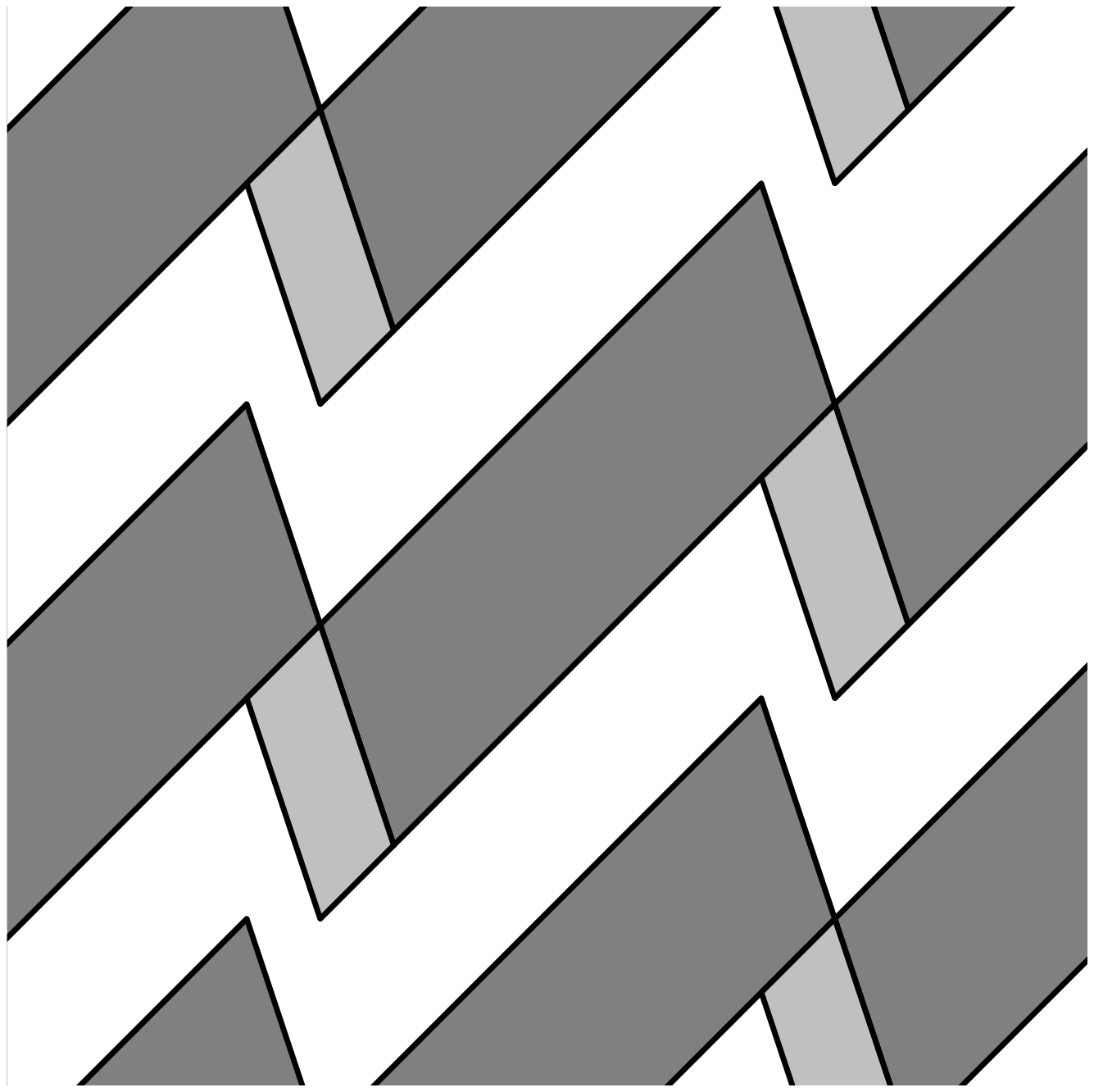}}
\newline
{\bf Figure 9.4:\/} The orbits
$\Lambda_{\Pi}H^{\rm lo\/}$ and
$\Lambda_{\Pi}H^{\rm hi\/}$ for
$A=1/4$ and $A=3/4$.
\end{center}

\noindent
{\bf Geometric Claim:\/}
\begin{itemize}
\item
$\zeta_{\Gamma} \in G_{\Gamma}^{\rm hi}$ implies 
$\Phi_{\Pi}(\zeta_{\Pi}) \in \Lambda_{\Pi}(H^{\rm hi\/} \times [-1,1])$
\item
$\zeta_{\Gamma} \in G_{\Gamma}^{\rm lo}$ implies 
$\Phi_{\Pi}(\zeta_{\Pi}) \in \Lambda_{\Pi}(H^{\rm lo\/} \times [-1,1])$
\end{itemize}
We call this claim the {\it geometric claim\/}.

We introduce the new sets
\begin{equation}
(H^{\rm hi\/})'=H^{\rm hi\/}+(2P,2P), \hskip 30 pt
(H^{\rm lo\/})'=H^{\rm lo\/}+(2P-2,2P-2).
\end{equation}
If $\zeta_{\Gamma} \in G_{\Gamma}^{\rm hi\/}$ and
the geometric claim is true, then
$$
\Phi_{\Pi}(\zeta_{\Pi}') \in \Lambda_{\Pi}((H^{\rm hi\/})' \times [-1,1]).
$$
The same goes when we replace {\it hi\/} with {\it lo\/}.
In the low case, we are using the fact that
$(2P-2,2P-2)$ and $(2P-2,2P)$ are equal up to a vector
in $\Lambda_{\Pi}$.

Now for the punchline. 
$\Lambda_{\Pi}((H^{\rm lo\/})')$ and
$\Lambda_{\Pi}((H^{\rm hi\/})')$ give
a second partition of
$\Lambda_{\Pi}(H)$. Figure 9.5 shows
the picture for the parameters $A=1/4$ and $A=3/4$.
\begin{center}
\resizebox{!}{2in}{\includegraphics{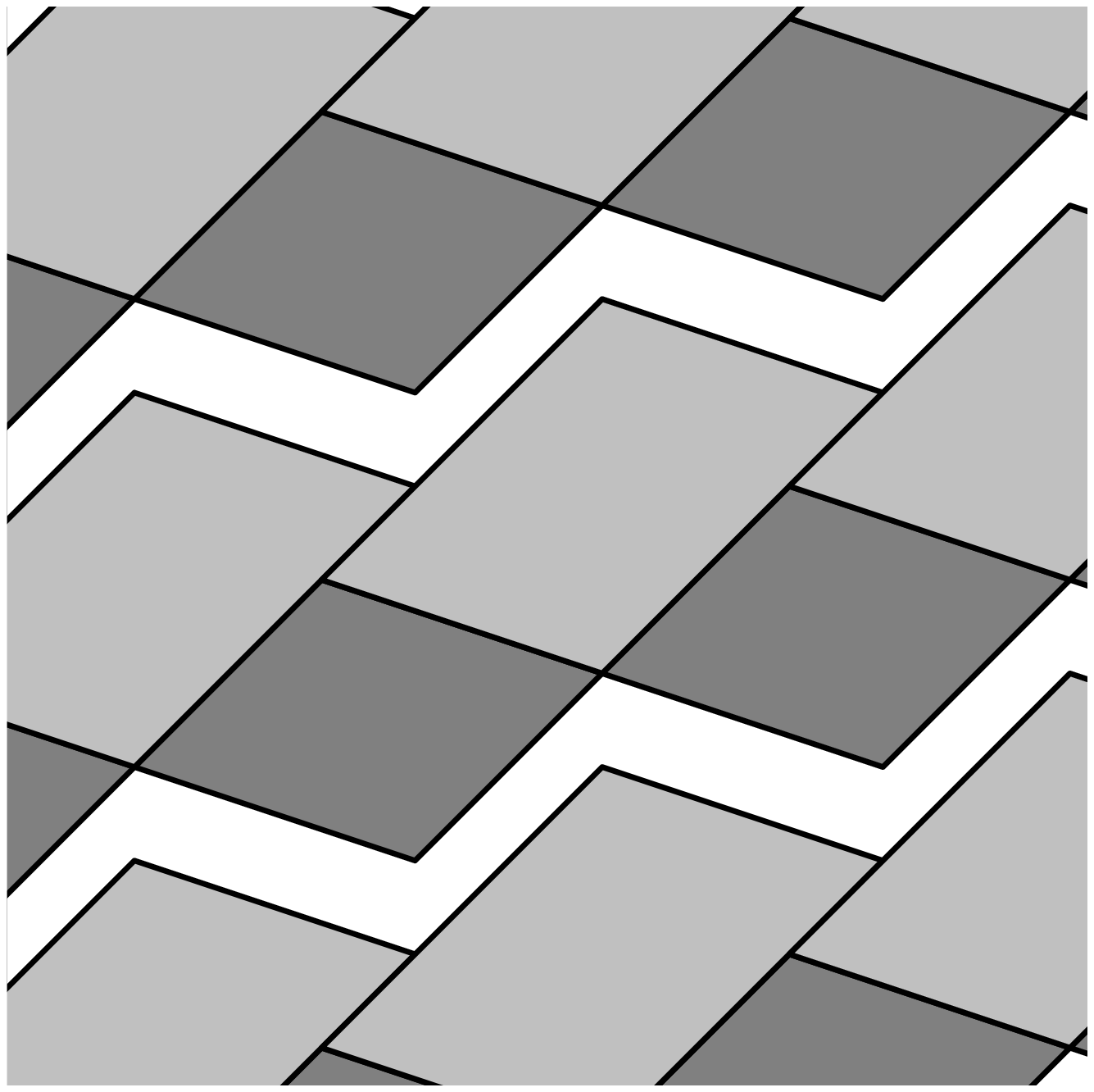}}
\hskip 50 pt
\resizebox{!}{2in}{\includegraphics{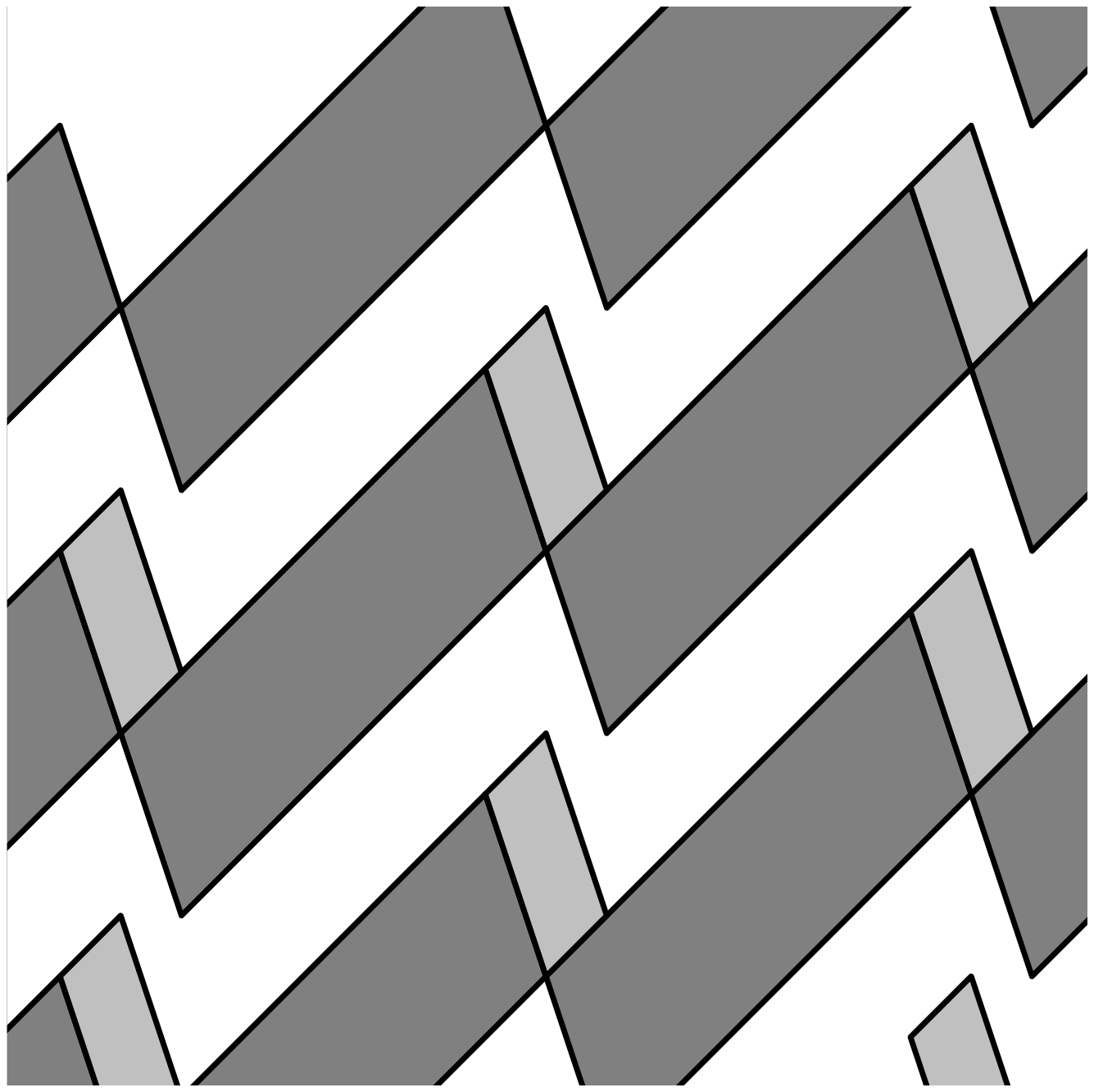}}
\newline
{\bf Figure 9.5:\/} The orbits
$\Lambda_{\Pi}(H^{\rm lo\/})'$ and
$\Lambda_{\Pi}(H^{\rm hi\/})'$ for
$A=1/4$ and $A=3/4$.
\end{center}
We get exactly the same picture
as in Figure 9.4 except that the pieces have
each been translated.   

What we are really
saying is that there is an infinite
polygon exchange transformation on 
$\Lambda_{\Pi}(H)$ which corresponds
to the operation of adding the vector
$dT(0,1)$ in $G_{\Gamma}$.  
The orbit
$\Lambda_{\Pi}(H)$ decomposes into
a countable union of parallelograms,
each partitioned into a light
parallelogram and a dark one. Our map
simply exchanges the light and dark
pieces within each component.  On
each component, our map is essentially
a rotation of a flat torus.

Figure 9.6 shows the two partitions
side by side, for the parameter $A=1/4$.
The ``components'' we are talking about are
parallelograms which are bounded on opposite
sides by lines of slope $1$.  
For later reference, we call these
parallelograms {\it dipoles\/}.
So, again, each dipole is partitioned
into a light and dark parallelogram, 
and our polygon exchange simply exchanges
the two pieces within each dipole.

\begin{center}
\resizebox{!}{2.5in}{\includegraphics{Pix/hilo1.ps}}
\hskip 10 pt
\resizebox{!}{2.5in}{\includegraphics{Pix/hilo3.ps}}
\newline
{\bf Figure 9.6:\/} The two partitions
of $\Lambda_{\Pi}(H)$ for the
parameter $A=1/4$.
\end{center}

The Geometric Claim immediately implies
Hitset Induction. Hence, the
Geometric Claim implies the Hitset Theorem.

\subsection{Changing the Fundamental Domain}
\label{domain}

Observe that the orbit
$\Lambda_{\Pi}(H)$ is simply a union of dipoles.
Rather than consider $X_{\Pi}=[-1,1]^3$ as our
fundamental domain for the
action of $\Lambda_{\Pi}$, we instead
consider the fundamental domain to be
\begin{equation}
\Upsilon \times [-1,1]
\end{equation}
where $\Upsilon$ is the dipole that intersects
$[-1,1]^2$ to the right of the diagonal line
$y=x$.  Figure 9.7 shows $\Upsilon$ and how
it sits with respect to $[-1,1]$.

\begin{center}
\resizebox{!}{3.5in}{\includegraphics{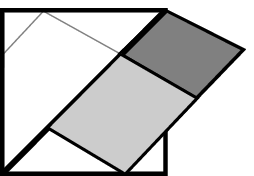}}
\newline
{\bf Figure 9.7:\/} The fundamental dipole.
\end{center}

We get an easier calculation if we
use $\Upsilon \times [-1,1]$ as the
domain for the projective intertwiner $\Psi$.
To do this, we need to define $\Psi$
on the whole fundamental domain, and
check that the new definition agrees
with the original in the appropriate sense.
\newline

The points of $\Upsilon-[-1,1]^2$ all lie
to the right of $[-1,1]^2$, and one
simply subtracts off $(-2,-P)$ to get them
back into $[-1,1]^2$.  For any
\begin{equation}
\phi \in \Bigg(\Upsilon-[-1,1]^2\bigg) \times [-1,1],
\end{equation}

We define
$\Psi(\phi)$ using the branch of $\Phi$ which is
defined for $X_{\Pi,-}$.  Let
$\phi^*=\phi-(2,P,P) \in X_{\Pi,+}$. Using the
fact that $\phi \in X_{\Pi,-}$ and
$\phi^* \in X_{\Pi,+}$, we compute
\begin{equation}
\Psi(\phi)-\Psi(\phi^*)=(-1,-A,A)=(-1,1,1)-(0,1+A,1-A) \in \Lambda_{\Gamma}.
\end{equation}
So, our redefinition does not change anything.
We can prove Induction Statement 2 using the new
fundamental domain.  Moreover, the restriction
of $\Psi$ to the new domain is projective
throughout.

\subsection{Intertwiner Induction}

Since $\Phi_{\Gamma}(x,y)=(x,x,x)$ mod $\Lambda_{\Gamma}$, we
have
\begin{equation}
\Phi_{\Gamma}(\zeta_{\Gamma}')-\Phi_{\Gamma}(\zeta_{\Gamma})=
(1,1,1) \hskip 10 pt {\rm mod\/} \hskip 5 pt
\Lambda_{\Gamma}.
\end{equation}
Let
\begin{equation}
\phi=\Phi_{\Pi}(\zeta_{\Pi}),
\hskip 40 pt
\phi'=\Phi_{\Pi}(\zeta'_{\Pi}).
\end{equation}
To establish Intertwiner Induction, we just have to prove
\begin{equation}
\Psi(\phi')-\Psi(\phi)-(1,1,1) \in \Lambda_{\Gamma}.
\end{equation}
 Using our new fundamental
domain, and also the action of the
map $\Phi_{\Pi}$, we have
one of
\begin{equation}
\phi-\phi'=(2P,2P,2P+2\beta), \hskip 30 pt
\phi'-\phi=(2P-2,2P-2,2\beta),
\end{equation}
depending on whether we are in the hi case or the lo case.
Here $\beta$ is some integer whose value does not
effect the calculation.
\newline
\newline
{\bf Hi Case:\/}
$$\Psi(\phi')-\Psi(\phi)-(1,1,1)=$$
$$\frac{1}{2-P}\bigg(0,-2P,2P+2\beta\bigg)- (1,1,1)=$$
$$(-1,-2A-1,2A-1+2\beta(1+A))=$$
$$2\beta(0,0,1+A)+(-1,1,1)+2(0,A+1,A-1) \in
\Lambda_{\Gamma}.$$
Notice that $\beta$ plays no role at all in the
final answer. To simplify the second calculation,
we assume that $\beta=0$.
\newline
\newline
\noindent
{\bf Lo Case:\/}
$$\Psi(\phi')-\Psi(\phi)-(1,1,1)=$$
$$\frac{1}{2-P}\bigg(0,-2P+2,0\bigg)- (1,1,1)=$$
$$(-1,-A,-1)=(-1,1,1) - (0,1+A,1-A) - (0,0,1+A) \in \Lambda_{\Gamma}.$$
This is what we wanted to prove.

At this point, we have reduced the Intertwining Theorem and the
Hitset Theorem to the Geometric Claim. 

\subsection{Proof of the Geometric Claim}
\label{geomproof}

We will combine the
Intertwining Formula 
with the Graph Reconstruction Formula. We
are allowed to do this for the pair
$(\zeta_{\Gamma},\zeta_{\Pi})$ by induction.

We introduce coordinates
$$
\zeta_{\Gamma}=(a,b), \hskip 18 pt
\Phi_{\Pi}(\zeta_{\Pi})=(x,y,z,P)$$.
Let $[t]=t-{\rm floor\/}(t)$.

\begin{lemma}
\begin{equation}
\label{PLAID2}
[a]=\left[\frac{x-y}{2-P}\right], \hskip 30 pt
[b]=
\left[\frac{-2-P+P^2+Px+2y-2Py}{2P-4}\right].
\end{equation}
\end{lemma}

\startproof
Define
$$
\Psi \circ \Phi_{\Pi}(\zeta_{\Pi})=(x^*,y^*,z^*,A).
$$

For convenience, we repeat Equation \ref{MAP} here:
\begin{equation}
\label{MAP2}
\Psi(x,y,z,P)=
\left[\frac{1}{2-P}\bigg(x-y,-y-1,z+P+1,P\bigg)-(1,0,0,0)\right]_{\Lambda}
\end{equation}
We always take the $(-)$ option in Equation
\ref{MAP} because we
are using the domain $\Upsilon \times [-1,1]$
described in \S \ref{domain}.
The graph reconstruction formula tells us that
$$
[a]=[x^*], \hskip 30 pt
[b]=\left[\frac{y^*-Ax^*}{1+A}\right].
$$
Combining this formula with Equation \ref{MAP2} and
doing some algebra, we get Equation \ref{PLAID2}.
\endproof

The next result says that
the formulas in Equation \ref{PLAID2}, which look
messy, are actually as nice as possible.

\begin{lemma}
\label{CANON}
Equation \ref{PLAID2} induces an affine
diffeomorphism from $\Upsilon$ to 
$[0,1]^2$.
\end{lemma}

\startproof
Forgetting about the brackets, the corresponding
map on the plane is an affine diffeomorphism $\Omega$.
The first coordinate of $\Omega$
is obviously constant along lines of slope $1$.
These lines are parallel to the diagonal sides of
$\Upsilon$.  We claim that the second coordinate
of $\Omega$ is constant along lines parallel to the other two
sides of $\Upsilon$.  To see this, we plug in
the equation
\begin{equation}
\label{lines}
y=\bigg(\frac{-P}{2-P}\bigg)\ x + c.
\end{equation}
and observe that the resulting expression
\begin{equation}
f(c)=\frac{2+P-P^2}{4-2P}+\bigg(\frac{P-1}{2-P}\bigg)\ c
\end{equation}
is independent of $x$.

The left and right sides of $\Upsilon$ are given
by the equations $y=x$ and $y=x-(2-P)$.
Hence $\Omega$ maps the left and right sides of
$\Upsilon$ respectively to the left and right sides
of the integer unit square $[0,1]^2$.
The bottom and top sides of $\Upsilon$ are given
by taking 
\begin{equation}
c_0=\frac{P-2}{2}, \hskip 30 pt
c_1=\frac{P-2}{2}+\frac{1-P}{2-P}
\end{equation}
in Equation \ref{lines}.
We check that $f(c_0)=1$ and $f(c_1)=0$.  Hence
$\Omega$ maps the top and bottom of
$\Upsilon$ respectively to the bottom and top
of $[0,1]^2$.
\endproof

The Geometric Claim follows immediately from
the analysis in the previous lemma. It says that
\begin{equation}
\Lambda_{\Gamma}(H^{\rm hi\/}) \cap \Upsilon=
\Omega^{-1}([0,1] \times [P,1]), \hskip 15 pt
\Lambda_{\Gamma}(H^{\rm lo\/}) \cap \Upsilon=
\Omega^{-1}([0,1] \times [0,P]).
\end{equation}
In light of Equation \ref{PLAID2}, this last
equation is equivalent to the Geometric Claim.

This completes the proof of the
Geometric Claim, and thereby completes
the proof of both the Hitset Theorem and
the Intertwining Theorem.

\newpage

\section{Correspondence of Polytopes}
\label{CP}

In this chapter we prove 
Statement 2  of the
Pixellation Theorem modulo
some integer linear algebra calculations.
We also set up the kind of
problem we need to solve in
order to prove Statements 3,4, and 5 of
the Pixellation Theorem.
We fix some parameter $A$ that is in the
background of the whole discussion.
Also, as in some previous chapters, we reserve
the word {\it square\/} for unit integer squares.

\subsection{Proof of Statement 2}

Statement 2 of the Pixellation Theorem
says that a grid full square is
plaid trivial if and only if it is
graph trivial.

The reduced
triple partition ${\cal RTP\/}$ consists
of $218$ polytopes in the
fundamental domain $X_{\Pi}=[-1,1]^3$.
We call these polytopes $\Pi_0,...,\Pi_{217}$.

\begin{lemma}
Suppose that $\Sigma$ is a grid full square.
If $\Sigma$ is plaid trivial then
$\Sigma$ is graph trivial.
\end{lemma}

\startproof
In our listing,
there are $6$ polytopes
in ${\cal RTP\/}$ having null labels,
namely $\Pi_0,...,\Pi_5$.
We observe that $\Psi(\Pi_j)$ is
contained in a single polytope
$$\Gamma_j=\Gamma_{j,+}=\Gamma_{j,-}$$ in
the graph partition, and that
$\Gamma_j$ is null labeled.
(This is why
$\Gamma_{j,+}=\Gamma_{j,-}$.)

Suppose that $\Sigma$ is a grid full
plaid trivial square for some parameter.
Let $\zeta_{\Pi}$ be the center of $\Sigma$.
Let $\zeta_{\Gamma}$ be the
graph grid point in $G_{\Gamma}$ which
is contained in $\Sigma$.  

By the
Plaid Master Picture Theorem,
$$\Phi_{\Pi}(\zeta_{\Pi}) \subset
\bigcup_{j=0}^5 \Pi_j.$$
By the Intertwining Theorem,
$$\Phi_{\Gamma}(\zeta_{\Gamma}) \subset
\bigcup_{j=0}^5 \Psi(\Pi_j) \subset
\bigcup_{j=0}^5 \Gamma_j.$$
By the Graph Master Picture Theorem,
the portion of the arithmetic graph
associated to $\zeta_{\Gamma}$ is
trivial. That is,
$\Sigma$ is graph trivial.
\endproof

\begin{lemma}
Suppose that $\Sigma$ is a grid full square.
If $\Sigma$ is graph trivial then
$\Sigma$ is plaid trivial.
\end{lemma}

\startproof
We will prove the contrapositive: If
$\Sigma$ is plaid nontrivial then
$\Sigma$ is graph nontrivial.
We observe that each of
the polytopes $\Pi_6,...,\Pi_{217}$ has
a nontrivial label.  We also observe the
following:
\begin{itemize}

\item $\Psi(\Pi_j)$ is contained in
a single graph polytope
$\Gamma_{j,+}$ and a single
graph polytope $\Gamma_{j,-}$ for
$j=6,...,179$.  The labels of
these graph polytopes are all nontrivial

\item $\Psi(\Pi_j)$ is contained in
a single graph polytope
$\Gamma_{j,-}$ and a union of
two graph polytope $\Gamma_{j,+,0}$ and
$\Gamma_{j,+,1}$ for $j=180,...,198$.
The labels of these graph polytopes are
all nontrivial.

\item $\Psi(\Pi_j)$ is contained in
a single graph polytope
$\Gamma_{j,+}$ and a union of
two graph polytope $\Gamma_{j,-,0}$ and
$\Gamma_{j,-,1}$ for $j=199,...,217$.
The labels of these graph polytopes are
all nontrivial.

\end{itemize}

Combining the information listed above
with the Plaid Master Picture Theorem,
the Graph Master Picture Theorem, and
the Intertwining Theorem in the same
way as the previous proof, we get the
conclusion of this lemma.
\endproof

\subsection{A Sample Result}
\label{sample0}

Let us press the method of the previous
section a bit harder, in order to see
both its power and its limitations.
Figure 10.1 shows a particular arc of the plaid
model.  Depending on the way the triple is
oriented, it corresponds to tiles in
${\cal RTP\/}$ labeled NEWEWS or SWEWEN.

\begin{center}
\resizebox{!}{.8in}{\includegraphics{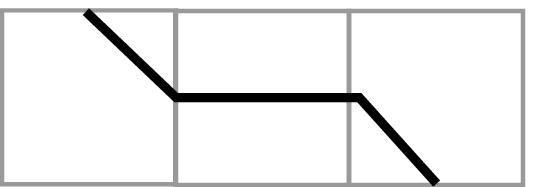}}
\newline
{\bf Figure 10.1:\/} The triple of type
NEWEWS or SWEWEN.
\end{center}

\begin{lemma}[Sample]
\label{sample}
Every time the triple of the type shown in 
Figure 10.1 appears in the arithmetic graph and
the central square $\Sigma$ 
is grid full, the two edges
incident to the grid graph vertex in $\Sigma$
are in ${\cal F\/}(0,1)$ and
${\cal F\/}(0,-1)$.  
\end{lemma}

\startproof
There are $6$ polytopes in ${\cal RTP\/}$
having the relevant labels.
The ones labeled NEWEWS are listed as
$\Pi_7,\Pi_{138},\Pi_{165}$.  The ones labeled
SWEWEN are listed as
$\Pi_{103},\Pi_{130},\Pi_{159}$.

Recall that $\Psi$ is the projective intertwining map.
We check by direct computation that, for each
$j \in \{7,103,130,138,159,165\}$ there are
polytopes $\Gamma_{j,+}$ and $\Gamma_{j,-}$ in the
$(+)$ and $(-)$ graph PET partitions respectively so that
\begin{equation}
\Psi(\Pi_j) \subset \Gamma_{j,+} \cap \Gamma_{j,-}.
\end{equation}
When we inspect the labels of these graph polytopes we
observe that
\begin{itemize}
\item The label of $\Gamma_{j,+}$ is always $(0,1)$.
\item The label of $\Gamma_{j,i}$ is always $(0,-1)$.
\end{itemize}

Now let us put this information together.
Suppose we see the triple from Figure 10.1 in
the plaid model for some parameter $A$.
The rest of the discussion implicitly refers
to the parameter $A$.
Let $\zeta_{\Gamma}$ be the vertex of the
graph grid $G_{\Gamma}$ contained in
the central square $\Sigma$ of the triple.
The same argument as in the previous section
shows that
\begin{equation}
\Phi_{\Gamma}(\zeta_{\Gamma}) \in
\bigcup_{j=7,103,130,138,159,165} \Gamma_{j,+} \cap \Gamma_{j,-}.
\end{equation}
But then, by the Master Picture Theorem,
the two edges incident to $\zeta_{\Gamma}$ are
in ${\cal F\/}(0,1)$ and ${\cal F\/}(0,-1)$.
\endproof

10.2 shows us the kind of conclusion we
can draw from Lemma \ref{sample}.

\begin{center}
\resizebox{!}{1in}{\includegraphics{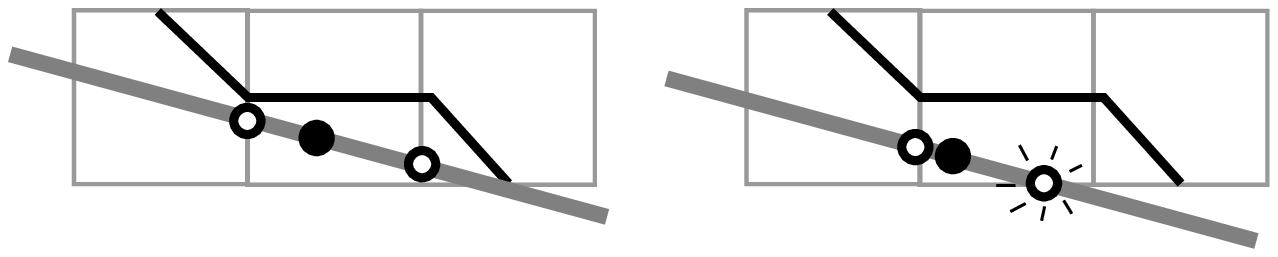}}
\newline
{\bf 10.2:\/} Some possible conclusions
from Lemma \ref{sample}.
\end{center}

Lemma \ref{sample} tells us that every time we
see the triple in Figure 10.1, we see the
arithmetic graph edges shown in (in grey) in
10.2.  However, it might happen that 
one of these grey edges crosses the top or
the bottom of the central square rather than
the sides, as it the Pixellation Theorem suggests.
Lemma \ref{sample} tells us the {\it local geometry\/}
of the arithmetic graph but not how it crosses
the edges of the central square.  So eliminate
the second option shown in 10.2 (as well
as other bad options, we need to study how our
machinery also determines the edge crossings.

In the rest of this chapter, we set up
conventions and notation for the kind of
problem we need to solve.

\subsection{Fixing Orientations}

Recall that there are both oriented and unoriented
versions of the Plaid and Graph Master
Picture Theorems.  We have stated the
Intertwining Theorem for the
unoriented versions, because the proof
is simpler.  However, this forces us
to deal with orientations in a
somewhat {\it ad hoc\/} way.

As we mentioned above, there are $218$
polytopes in ${\rm RTP\/}$.
Each of these polytopes has a $6$ letter
label which specifies an oriented
arc of combinatorial length $3$ in
the plaid model.  The catch is that
we are using the small 
domain $\widehat X/\Lambda_1$ as
the image of the plaid classifying map
$\Phi_{\Pi}$, rather than the
double cover $\widehat X/\Lambda_2$. 

Thus, if $\zeta_{\Pi} \in G_{\Pi}$
is some point, the label of the
polytope in ${\cal RTP\/}$ containing
$\Phi_{\Pi}(\zeta_{\Pi})$ might
specify the opposite orientation
on the plaid arc through $\zeta_{\Pi}$.
More precisely, the assigned orientation
is correct if and only if
$\widehat \Phi_{\Pi}(\zeta_{\Pi})$ is
congruent mod $\Lambda_2$ to a point
in the fundamental domain
$X_{\Pi}$.  Here
$\widehat \Phi_{\Pi}$ is the lift of
$\widehat \Pi$ to the double cover
$\widehat X/\Lambda_2$.
At the same time, in order to figure
out the orientations on the graph
polygons, we would need to look at
the lift $\widehat \Phi_{\Gamma}$ of
the graph classifying map
$\Phi_{\Gamma}$.

Getting these orientations right is 
a tedious business, and we have a
different way of dealing with it.
We ignore the ``true'' orientations
coming from the lifted maps and we
simply make a guess as to how
the plaid arcs correspond to the
graph arcs.  For example, in 10.2,
it is pretty clear that the grey
edge pointing left should be associated
to the left half of the plaid triple.

Now we explain the convention in more
detail, by way of example.  Figure
10.3 shows an enhanced version of
10.2, specifically for 
$\Pi_7$.  Again, the label of
$\Pi_7$ is NEWEWS, and this determines
the orientation.

\begin{center}
\resizebox{!}{1.1in}{\includegraphics{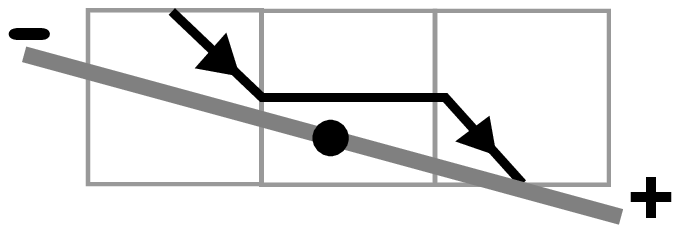}}
\newline
{\bf Figure 10.3:\/} The case of $\Pi_7$: a correspondence of type $1$
\end{center}

We observe experimentally that the $(-)$ partition
assigns the label $(0,-1)$, and this corresponds
to the leftward pointing edge. Likewise,
the $(+)$ partition assigns the label $(0,1)$,
and this corresponds to the rightward pointing
edge.  Thus, we associate the $(-)$ partition with the
tail end of the triple and the $(+)$ partition
with the head end.  For this reason, we call
$\Pi_7$ {\it type 1\/}.  We would call $\Pi_j$
{\it type 0\/} if, according to our experimental
observations, the $(+)$ edge is associated
with the tail and the $(-)$ edge is
associated with the tail.  

We guess the type for each of the $218$ polytopes
in the plaid triple partition, and these types
are stored in the computer program.  In most cases,
as for $\Pi_7$, the picture is completely obvious,
and in a few cases, like the one shown in 10.2
for $\Pi_{49}$ the picture is not quite as obvious
but still pretty clear.

\begin{center}
\resizebox{!}{1.9in}{\includegraphics{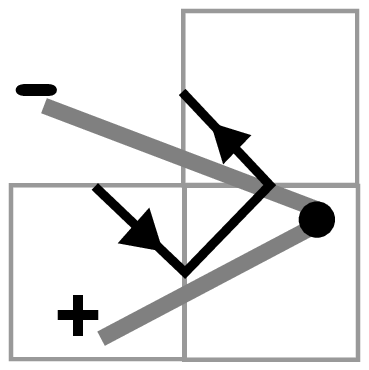}}
\newline
{\bf Figure 10.4:\/} The case of $\Pi_{48}$: a correspondence of type $0$.
\end{center}
 
Logically, it does not matter how we arrived at these
guesses, or whether they agree with the
true answer which can be cleaned by looking
at the lifts.  The point is simply that the
proof runs to completion with the guesses
made.  In hindsight, our guesses surely
agree with the true answer, but we do not
need to prove this and we will not.

\subsection{Edge Crossing Problems}
\label{ECP}

To each index $i=6,...,218$ we associate what we
call an {\it edge crossing problem\/}.  We will
first give the general definition and then we
will work out the example of $\Pi_7$.

Given $\Pi_k$, there corresponds a triple of
squares $\Sigma_{-1},\Sigma_1,\Sigma_1$
and an oriented plaid arc running through it,
as in Figures 10.3 and 10.4.  We label the
squares so that $\Sigma_1$ is the head square
in the type $1$ case and the tail square
in the type $0$ case.  In other words, the
arithmetic graph edge associated to the
$(\pm)$ partition should point generally
from $\Sigma_0$ to $\Sigma_{\pm}$ for
both types.  

There are two kinds of edge crossing problems
associated to $P_k$.  One kind is labeled
$(k,+,i,j,L)$.  Here $L \in \{N,S,E,W\}$ is an
edge of $\Sigma_0$ that the arithmetic
graph edge $dT(i,j)$ 
associated to the $(+)$ partition
could potentially cross (as in the right hand
side of Figure 10.3) but according to the
Pixellation Theorem
is not supposed to cross.   These crossing
problems really only depend on the pair
$(\Sigma,\Sigma_+)$.  The other kind of
edge crossing problem is labeled $(k,-,,i,j,L)$,
and has a similar explanation with $(-)$
in place of $(+)$.

Let's consider the case of $P_7$.  The two
crossing problems are $(7,+,0,1,N)$ and
$(7,-,0,-1,S)$.  The other two possibilities,
namely $(7,+,0,1,E)$ and $(7,-,0,-1,W)$, are the ones
predicted by the Pixellation Theorem.
They are not problems at all, but rather
{\it goals\/}.

We have already mentioned that the first
$6$ of the plaid triple polytopes correspond
to the trivial grid full squares.
There are another $174$ polytopes,
labeled $\Pi_6,...,\Pi_{179}$, which
have the property that there are unique
graph polyropes $\Gamma_{k,+}$ and
$\Gamma_{k,-}$ such that
\begin{equation}
\Psi(\Pi_k) \subset \Gamma_{k,\pm}.
\end{equation}
Each of these contributes $2$ edge
crossing problems, giving a total
of $348$.

The $19$ polytopes
$\Pi_{180},...,\Pi_{198}$ are such that
\begin{equation}
\Psi(\Pi_k) \subset \Gamma_{k,+,0} \cup \Gamma_{k,+,1}, \hskip 30 pt
\Psi(\Pi_k) \subset \Gamma_{k,-}.
\end{equation}
That is, $\Psi(\Pi_k)$ is contained in a union of two
graph polytopes from the $(+)$ partition and $1$ from
the $(-)$ partition.  In this case, there are two
possible local picture of the arithmetic graph associated 
to this plaid triple.  This does not bother us, as long
as the edge crossings come out right.  Each of these
$19$ polytopes contributes $3$ edge crossing problems.
This gives us another $57$.

The $19$ polytopes
$\Pi_{199},...,\Pi_{217}$ are such that
\begin{equation}
\Psi(\Pi_k) \subset \Gamma_{k,-,0} \cup \Gamma_{k,-,1}, \hskip 30 pt
\Psi(\Pi_k) \subset \Gamma_{k,+}.
\end{equation}
Each of these
$19$ polytopes contributes $3$ edge crossing problems,
giving yet another $57$ crossing problems.

The grand total is $462$ edge crossing problems.
In the next chapter we will introduce the machinery
needed to solve 
all these problems, so to speak. A {\it solution\/} amounts
to a proof that the given case does not actually occur.
Actually, we will be able to solve $416$ of the
problems. The remaining $46$, in the cases
where the pixellation really does fail, are
involved in the
catches for the offending edges discussed in
the Pixellation Theorem.  Once we have solved
all the problems we can solved, and classified
the exceptions, the rest of the Pixellation
Theorem just comes down to inspecting the
data generated by the program.

\newpage

\section{Edge Crossings}

\subsection{The Graph Method}

Here we explain the first method we  use for
solving the edge crossing problems.
We fix a parameter $A$ throughout the
discussion.   Let $G_{\Gamma}$ denote the
grid graph.  Let $\zeta_{\Gamma} \in G$ be some point,
contained in a square $\Sigma$.
As in previous chapters, the
word {\it square\/} always means
a unit integer square.  

We assume that $\zeta_{\Gamma}$ is a nontrivial
vertex of the arithmetic graph.
Let $e$ be one of the edges
of the arithmetic graph incident to $\zeta_{\Gamma}$.
We think of $e$ as a vector pointing
from $v$ out of $\Sigma$.  
From the Grid Geometry Lemma, we know
that $e$ crosses some edge of $\Sigma$.
We say that $e$ is of type
$(i,j,L)$ if $e=dT(i,j)$ and
$L \in \{S,W,N,E\}$ is the label of
the edge of $\Sigma$ which $e$
crosses.  
\newline
\newline
{\bf Remark:\/}
We allow the possibility
that $e$ is of two types. This
would happen if $e$ crosses
$\Sigma$ at a vertex.  We think that
this never actually happens, but
we have not ruled it out.
In any case, this eventuality
does not bother us.
\newline

Let $\Phi_{\Gamma}$ denote the graph classifying
map.  We define
\begin{equation}
(x,y,z,A)=\Phi_{\Gamma}(\zeta_{\Gamma})
\end{equation}

\begin{lemma}[Graph Avoidance]
The following is true.
\begin{enumerate}
\item If $x \in (A,1)$ then $e$ is not of type $(-1,0,W)$.
\item If $x \in (0,A)$ then $e$ is not of type $(-1,1,W)$.
\item If $y \in (A,1)$ then $e$ is not of type $(0,-1,N)$.
\item If $y \in (2A,1+A)$ then $e$ is not of type $(0,1,S)$.
\item If $x \in (1-A,1)$ then $e$ is not of type $(1,-1,E)$.
\item If $x \in (0,1-A)$ then $e$ is not of type $(1,0,E)$.
\end{enumerate}
\end{lemma}

\startproof
We will treat the first three cases.
The last three cases follow from symmetry.
More precisely, reflection in the center of the
square $[0,1]^2$ carries the sets used
to analyze Case N to the sets used to
analyze Case 7-N. Figure 11.1 illustrates
our arguments.

\begin{center}
\resizebox{!}{1in}{\includegraphics{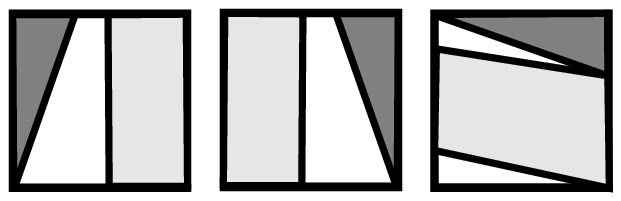}}
\newline
{\bf 11.1:\/} The first three cases
\end{center}

\noindent
{\bf Case 1:\/} We have 
$$dT(-1,0)=\bigg(-A,-\frac{1+2A-A^2}{1+A}\bigg).$$
This vector points southwest and has slope
greater than $1$.  So, in order for this edge
to cross the west edge, it must lie in the triangle
with vertices
$$(0,0), \hskip 15 pt (0,1) \hskip 15 pt \bigg(\frac{A^2+A}{1+2A-A^2}\bigg).$$
We call this triangle {\it the danger zone\/}.
The line through this third point and parallel to $e$ contains
the vertex between the south and west edges.
The dark triangles in Figure 11.1 are hand-drawn versions
of the danger zones in each case.

The Graph Reconstruction Formula tells us that
the conditions $x \in (A,1)$ correspond to the condition
that the first coordinate of $\xi_{\Gamma}$ lies in $(A,1)$.
The region of possibilities is the lightly shaded region
on the left hand side of figure 11.1.  As depicted in the
figure, the two sets we have defined are disjoint.
\newline
\newline
{\bf Case 2:\/} The argument is the same, except this time the
danger zone has vertices 
$$(1,0), \hskip 10 pt (1,1), \hskip 15 pt \bigg(\frac{4A}{4A-A^2}\bigg),$$
and the condition $x \in (0,A)$ corresponds to these same conditions
on the first coordinate of $\xi_{\Gamma}$.  Again, the two sets are disjoint.
This is shown in the middle square of Figure 11.1.
\newline
\newline
{\bf Case 3:\/} The argument is the same, except this time the
danger zone has vertices 
$$(0,1), \hskip 10 pt (1,1), \hskip 15 pt \bigg(1,\frac{1-A}{1+A}\bigg).$$
and the condition $y \in (A,1)$ corresponds, {\it via\/} the
Reconstruction Formula to $\xi_{\Gamma}$ lying in the interior
of the parallelogram with vertices
$$\bigg(0,\frac{A}{1+A}\bigg), \hskip 15 pt
\bigg(0,\frac{1}{1+A}\bigg), \hskip 15 pt
\bigg(1,\frac{1-A}{1+A}\bigg), \hskip 15 pt
\bigg(0,0\bigg).$$
Again, these sets are disjoint.  This case
is shown on the right side of Figure 11.1.
\endproof

\noindent
{\bf Remark:\/}
In the previous result, the cases $(1,1,N)$ and
$(-1,-1,S)$ are missing.  There is a similar
result for these cases, but we will prove
a result below that is more powerful and
subsumes these cases.  So, we ignore them here.

\subsection{The Sample Result Revisited}

In this section we explain how we
solve the crossing problem $(7,+,,0,1,S)$.
In other words, we are trying to 
rule out the bad crossing indicated on the
right hand side of Figure 10.2, which
we repeat here with enhanced labeling.
We really need to solve $6$ crossing

\begin{center}
\resizebox{!}{1.7in}{\includegraphics{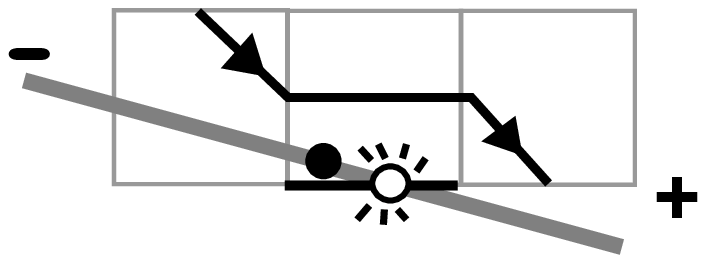}}
\newline
{\bf 11.2:\/} The crossing problem $(7,+,S)$.
\end{center}

The polytope $\Psi(\Pi_7)$ has $8$ vertices:
\begin{itemize}
\item $(60,60,0,30)/60$.
\item $(60,40,20,20)/60$.
\item $(40,40,0,20)/60$.
\item $(60,40,0,20,)/60$.
\item $(60,60,0,20)/60$.
\item $(45,45,0,15)/60$.
\item $(60,45,15,15)/60$.
\item $(60,45,0,15)/60$.
\end{itemize}
We have written the vertices this way
so as to clear denominators.  The factor
of $60$ works for every polytope in sight.
(My program has a window which allows the
user to see the vertices of any polytope
in any of the partitions.)

Recalling that we coordinatize $\widehat X$
using the variables, note that
\begin{equation}
\label{slab}
y\in [2A,1+2A].
\end{equation}
for all vertices of $\Psi(\Pi_7)$.  By
convexity, this equation holds for all
points, and we get strict inequality for
points in the interior of $\Psi(\Pi_7)$.

It follows from the Master Picture Theorems
and the Intertwining Theorem that
\begin{equation}
(x,y,z,A)=\Phi_{\Gamma}(\zeta_{\Gamma}) \in {\rm interior\/}(\Psi(\Pi_y)).
\end{equation}
Hence $y \in (2A,1+2A)$.
Case 4 of the
Graph Avoidance Lemma, which
pertains to the triple $(0,1,S)$,
solves this crossing problem.

The argument works the same way for all $6$
crossing problems associated to the
bad crossing shown in Figure 11.2.
It follows from symmetry (or from a
similar argument) that the $6$
crossings associated to the other
arithmetic graph edge in Figure 11.2
are also soluble.  Thus, every time
this pattern occurs in the plaid model for
any parameter, and the central square is
grid full, the central square is pixellated.
\newline
\newline
{\bf Remark:\/}
Sometimes 
we use the Graph Avoidance Lemma in a different way.
To illustrate our other usage, we
will solve the edge crossing
problem $(198,+,-1,1,E)$.
This time it turns out that
the criterion in the Graph Avoidance Lemma
for $(1,-1,E)$ does not hold for
$\Psi(\Pi_{198})$.  However, it does
hold for $\Gamma_{198,+,0}$, the relevant
one of the two graph polytopes in the $(+)$
graph partition whose union contains
$\Psi(\Pi_{198})$.  So, again,
we find that $\Phi_{\Gamma}(\xi_{\Gamma})$
cannot lie in the region corresponding
to a situation where the graph edge
incident to $\zeta_{\Gamma}$ crosses
$E$. 

For all the relevant crossing problems on which
we use the Graph Avoidance Lemma, we will
either apply the criteria to $\Psi(\Pi)$ or
to the relevant graph polytope $\Gamma$.
In short, we will use the one method or the other.

\subsection{The Plaid Method}

Here we discuss a second method for solving
crossing problems.
In the discussion that follows, $L \in \{N,S,E,W\}$
stands for one of the edge labels.
To $10$ of the $16$ pairs $(i,j,L)$.
Let $\langle i,j,L \rangle$ denote the
subset of $\widehat X$ which assigns the
edge $e=dT(i,j)$ to grid graph points $\zeta_{\Gamma}$.
in such a way that $e$ crosses the edge $L$
of the square containing the point $\zeta_{\Gamma}$.
The Graph Avoidance Lemma can be interpreted
as saying that certain regions in $\widehat X$
avoid $\langle i,j,L\rangle$.  For instance, the
set $x \in (A,1)$ is disjoint from $\langle -1,0,W\rangle$.

Recall that $X_{\Pi}^*$ is the hitset.
We are going to associate
a polytope $Z(i,j,L) \subset X_{\Pi}$ such that
one of two things is true:
\begin{itemize}
\item If $P \cap Z(i,j,L)=\emptyset$ then
$\Psi(P) \cap  \langle i,j,L  \rangle=\emptyset$.
\item If $P \subset Z(i,j,L)$ then
$\Psi(P \cap X_{\Pi}^*) \cap  \langle i,j,L  \rangle=\emptyset$.
\end{itemize}
In the first case, we call $P$ an {\it excluder\/} and in
the second case we call $P$ a {\it confiner\/}.
The other $6$ pairs we simply ignore.

For each fixed parameter, our polytopes all have
the form 
\begin{equation}
Z(i,j,L)=Z'(i,j,L) \times [-1,1],
\end{equation}
Where $Z'(i,j,L)$ is a polygon in the $xy$ plane.
We will list $4$ of the $10$ sets. The other $5$
are obtained from the first $5$ {\it via\/} the
following symmetry.
\begin{equation}
Z(-i,-j,L^{\rm opp\/})=-Z(i,j,L).
\end{equation}
Here $L^{\rm opp\/}$ is defined to be the edge
opposite $L$.  For instance $E^{\rm opp\/} = W$.
We list the $5$ polygons $Z'(i,j,L)$ 
as a function of the parameter
$A$, and recall that $P=2A/(1+A)$.
The first set is an excluder and the other $5$ are confiners.
Here are the sets:
\begin{itemize}
\item $Z'(1,1,N)$:
$(1-P,1-P)$, $(P-1,P-1)$, $(P-1,-1)$, $(1-P,-1)$.
\item $Z'(1,1,E)$:
$(1-P,1-P)$, $(-1,-1)$, $(1-P,-1)$.
\item $Z'(0,1,S)$:
$(-1,-1)$,$(P-1,P-1)$, $(1,P-1)$, $(1-P,-1)$.
\item $Z'(1,0,E)=Z'(1,1,E)$.
\item $Z'(0,-1,W)=Z'(0,1,S)$.
\end{itemize}
Conveniently, all these sets intersect the
hitset inside the fundamental dipole $\Upsilon$
considered in \S \ref{geomproof}.  This makes
our analysis easy.

\begin{center}
\resizebox{!}{1.4in}{\includegraphics{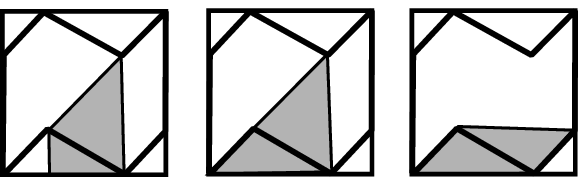}}
\newline
{\bf Figure 10.4:\/} The barrier bases
in action.
\end{center}

To establish our claims, we use the
Plaid Reconstruction Formula, Equation 
\ref{PLAID2}, to map our our
sets into the unit square, and we check
that the relevant image is disjoint
from the set of positions in $[0,1]^2$
where $\zeta_{\Gamma}$ can be placed
so that $dT(i,j)$ crosses edge $L$.
We called these sets the danger zones
in the proof of the Graph Avoidance Lemma.

\subsubsection{Case 1}

Here we consider $(1,1,N)$. Calculations like
the one done in the proof of the Graph Avoidance
Lemma show that the danger zone is the
triangle with vertices
$$(1,1), \hskip 30 pt (0,1), \hskip 30 pt(0,P).$$
Let $Z'_*(1,1,N)$ denote the intersection
of $Z'(1,1,N)$ with the planar projection of
the hitset.  The vertices of $Z'_*(1,1,N)$ are
$$(P-1,P-1), \hskip 30 pt
(1-P,-1), \hskip 30 pt (1-P,1-P).$$
Beautifully, the image of this triangle
under the affine diffeomorphism from
Equation \ref{PLAID2} is exactly
the danger zone.

\subsubsection{Case 2}

Here we consider $(1,1,E)$.  The danger zone
is the complement
of the triangle considered in Case 1.
We have
$Z'_+(1,1,N)=Z'_+(1,1,E)$.  So, we
have already computed the relevant affine image;
it is the set $\Sigma_0(1,1)$ from Case 1.  But the
interiors of $\Sigma_0(1,1)$ and $\Sigma_1(1,1)$
are disjoint: these two sets partition the unit square.

\subsubsection{Case 3}

Here we consider $(0,1,S)$.  The danger zone
is the same as the one in Case 4 of
the Graph Avoidance Lemma. It has vertices
$$(0,0), \hskip 30 pt (1,0), \hskip 30 pt
(0,P).$$
The set $Z'_*(0,1,Z)$ is the
triangle with vertices
$$(P-1,P-1), \hskip 30 pt
(1-P,-1), \hskip 30 pt (1,P-1).$$
The image of this set under the affine diffeomorphism is
the triangle with vertices
$$(0,1), \hskip 30 pt (1,1), \hskip 30 pt (1,P).$$
This triangle is clearly disjoint from the
danger zone.

\subsubsection{Case 4}

Here we consider $(1,0,E)$.  
The danger zone has vertices
$$(1,1), \hskip 30 pt (1,0), \hskip 30 pt
\bigg(\frac{1+A-2A^2}{1+2A-A^2}\bigg)$$
We have $Z'(1,0,E)=Z'(1,1,E)$. The analysis
in Case 1 shows that the affine image of
$Z'_+(1,0,E)$ is the triangle with vertices
$$(1,1), \hskip 30 pt (0,1), \hskip 30 pt(0,P).$$
The interior of this set is clearly disjoint
from the danger zone.

\subsubsection{Case 5}

Here we consider $(0,-1,W)$.  
The danger zone has vertices
$$(0,0), \hskip 30 pt (1,0), \hskip 30 pt
(1,P).$$
We have $Z'(0,-1,W)=Z'(0,1,S)$, the set
from Case 3. As in Case 3, the vertices
of the affine image of $Z'(0,-1,W)$ is
the triangle with vertices 
$$(0,1), \hskip 30 pt (1,1), \hskip 30 pt (1,P).$$
The interior of this set is disjoint from the
danger zone.

\subsection{Another Edge Crossing Problem}

Here we solve another edge crossing problem,
to illustrate the logic behind the use of
the Plaid Method.  We use the notation from
above. One of the crossing
problems is $(50,+,-1,-1,W)$.
This time we find that
$\Pi_{50} \subset Z(-1,-1,W)$.
Therefore
Hence 
\begin{equation}
\Psi_{\Pi}(\zeta_{\Pi}) \subset Z(-1,-1,W).
\end{equation}
The existence of $\zeta_{\Gamma}$ means that
the relevant square $\Sigma$ is grid full.

\begin{equation}
\Psi_{\Pi}(\zeta_{\Pi}) \subset Z(-1,-1,W) \cap
X_{\Pi}^*.
\end{equation}
By the Intertwining Theorem 
\begin{equation}
\Phi_{\Gamma}(\zeta_{\Gamma}) 
\subset \Psi(Z(-1,-1,W) \cap X_{\Pi}^*).
\end{equation}
By the confining property of $Z(-1,-1,W)$,
\begin{equation}
\Phi_{\Gamma}(\zeta_{\Gamma}) \not \in
\langle -1,-1,W \rangle.
\end{equation}
This solves the crossing problem.

\subsection{Out of Bounds}

Sometimes none of the above methods works for a
crossing problem, but then we notice
that $\Pi_k$ lies entirely outside the
hitset $Z_{\Pi}^*$.  In this case we call
$\Pi_k$ {\it out of bounds\/}.  The
corresponding crossing problem simply does not
arise for a grid full square.

We will verify that
$\Pi_k$ is out of bounds by showing that
\begin{equation}
\Pi_k \subset
\bigcup_{i=1}^4 B_i \times [-1,1].
\end{equation}
Where $B'_1,B'_2,B'_3,B'_4$ are the following $4$
triangles, described in terms of their vertices:
\begin{enumerate}
\item $B'_1: (-1,-1),(P-1,-1),(P-1,P-1)$.
\item $B'_2: (1-P,-1),(1,-1),(1,P-1)$.
\item $B'_3=-B'_1$.
\item $B'_4=-B'_2$.
\end{enumerate}
Figure 11.5 shows how these triangles sit in
relation to the planar projection of the hitset.

\begin{center}
\resizebox{!}{2.5in}{\includegraphics{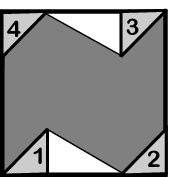}}
\newline
{\bf Figure 11,5:\/} The out-of-bounds polygons
in action.
\end{center}

Note that the union
\begin{equation}
B_i=\bigcup_{P \in [0,1]} (B'_i \times \{P\})
\end{equation}
is a convex polytope with integer coefficients.
Here $P=2A/(1+A)$ as usual.
This, showing that
$\Pi_k \subset B_i$ for some pair $(k,i)$ is just
a matter of integer linear algebra.  We
call the verification that $\Pi_k \subset B_i$
the {\it out of bounds test\/}.

\newpage

\section{Proof of the Pixellation Theorem}

Pixellation Theorem.  In this chapter, we prove
Statements 3,4,5 of the Pixellation Theorem
 modulo certain integer
computer calculations.

\subsection{Solving Most of the Crossing Problems}

Using the methods discussed in the last chapter,
namely the Graph Method, the Plaid Method, and
the Out of Bounds Test, we solve
$416$ of the $462$ crossing problems.
Each case is like one of the ones
considered in the previous chapter, but
there are too many to do by hand.
In the next chapter we will explain
the integer linear algebra tests we
use to check each case.

There are $46$ exceptional cases. Let's
call these $46$ cases {\it recalcitrant\/}.
 A few of the recalcitrant cases actually
are soluble by more delicate methods - we
would just need to look harder at the polytopes
involved - but there isn't any problem in
our proof with leaving them unsolved.

\subsection{Proof of Statement 3}

Figure 4.3 shows examples of an errant
edge when the two squares are stacked 
side by side.  We repeat the picture here.
We will consider the picture on the
right hand side of Figure 12.1 in detail.
In this picture, the square with the
dot is the central square and the top square
is not part of the triple.  So, we are only
showing two of the three squares in the triple.

The codes associated to the figure on the
right are either LLLEWS or SWELLL, depending
on the orientation. Here $L \in \{N,S,E,W\}$ is
a label we don't know.  The corresponding
errant edge must rise up at least $1$ unit. 
This leaves $dT(1,-1)$ and $dT(1,0)$ as
the only possibilities.

\begin{center}
\resizebox{!}{1.7in}{\includegraphics{Pix/errant.eps}}
\newline
{\bf Figure 12.1:\/} errant edges
\end{center}

  From the
correspondence of polytopes explained in
\S \ref{CP}, we know the complete list of
labels associated to relevant ends of the
relevant triples.  We just check, by a
quick computer search (and also by
inspection) that the errant labels never
arise.  Thus, for instance, we never
see the label $(1,0)$ associated to the
head end of LLLEWS or the tail end of
SWELLL.  We rule out the other possibilities
similarly.  Here is a list of the possible codes
and the forbidden labels.  To save space, we
only list half the possibilities.  The other
half are obtained from these by $180$ degree
rotation. Also, of the half we do list, we
only list the codes such that the corresponding
errant edge would be associated to the head
of the triple.  So, we would list LLLEWS only
in our example above. Here is the list.

\begin{itemize}
\item LLLEWS: $(1,0)$, $(1,-1)$.
\item LLLEWN: $(-1,0)$, $(-1,1)$.
\item LLLEWE: $(1,0)$, $(1,-1)$, $(-1,0)$, $(-1,1)$.
\item LLLSNE: $(-1,-1)$, $(0,-1)$.
\item LLLSNW: $(1,1)$, $(0,1)$.
\item LLLSNS: $(-1,-1)$, $(1,1)$, $(0,1)$, $(0,-1)$.
\end{itemize}

We simply check that none of these bad situations
actually arises. 
The rest of the proof of the 
Pixellation Theorem is devoted to
proving Statements 3 and 4.
This takes more work.

\subsection{Proof of Statement 4}

Statement 4 of the Pixellation
Theorem says that two arithmetic graph
edges incident to a vertex in a
square never cross the same edge of
that square.  This result follows
immediately for the $416$ cases
in which can solve the crossing problem:
These cases are all pixellated and
the edges in question cross the same
sides as the plaid model segments.
It just remains to deal with the
$46$ recalcitrant cases.

\begin{center}
\resizebox{!}{2in}{\includegraphics{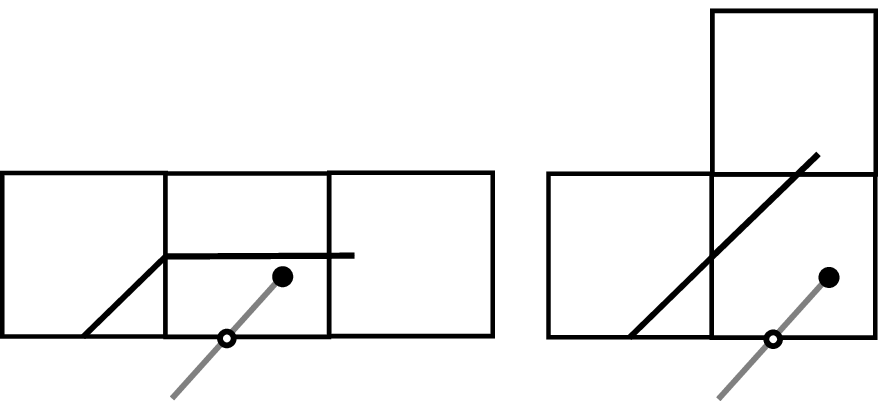}}
\newline
{\bf Figure 12.2:\/} The recalcitrant cases
\end{center}

By direct inspection, we see that
all of the recalcitrant cases are either
pixellated or equivalent to the two cases
shown in Figure 12.2. (By symmetry, we
just have to inspect $23$ of the $46$
cases.) By {\it equivalent\/}
we mean that the picture is meant to be taken
up to rotations and reflections.

Now consider the arithmetic graph edges
associated to the a recalcitrant case.
The first edge, the one shown in Figure 12.2,
is offending.  The other edge is either offending
or not.  If the other edge is not offending,
then it manifestly crosses a different edge of
the central square.  If the other edge is
offending, and crosses the same edge as
the first offending edge, then the picture
must look like Figure 12.3.

\begin{center}
\resizebox{!}{1.6in}{\includegraphics{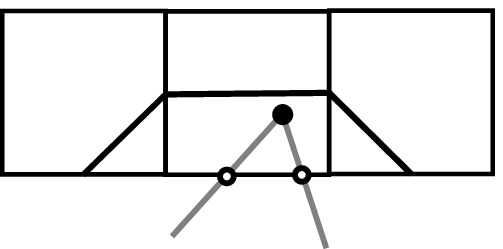}}
\newline
{\bf Figure 12.3:\/} A double edge crossing
\end{center}

There are exactly two cases like this,
corresponding to $P_{34}$ and $P_{173}$.
But, in both cases, we check that the
polytope is involved in only one 
recalcitrant crossing problem.  So, in
these two cases, only one of the two
edges is offending. 

This takes care of all the possibilities.

\subsection{Proof of Statement 5}

Each recalcitrant case corresponds to a
plaid triple.  Using the curve-following
dynamics discussed 
in \S \ref{CFD}, we check which plaid triples
could attach to the one we have.  We will
illustrate what we mean by example and
then state the general result.

\begin{center}
\resizebox{!}{1.8in}{\includegraphics{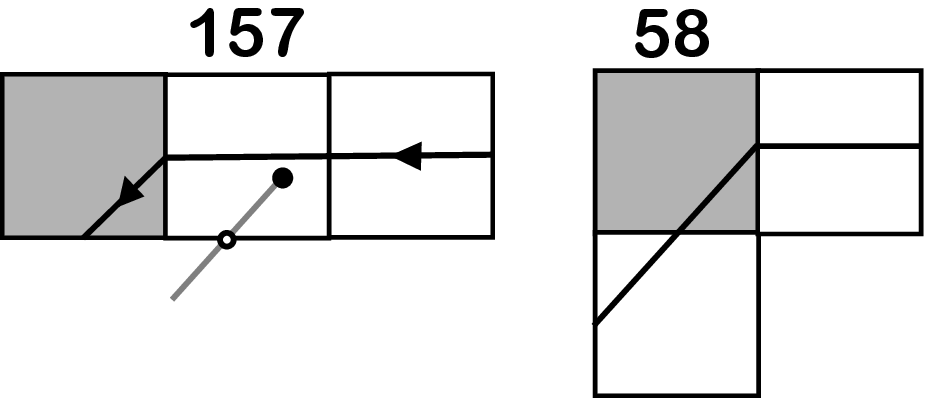}}
\newline
{\bf Figure 12.4:\/} The triples corresponding
to $\Pi_{157}$ and $\Pi_{58}$.
\end{center}

One of the recalcitrant triples
is $(157,-,-1,0,S)$.  The associated
code for $\Pi_{157}$ is EWEWES.  Figure
12.4 shows the situation.
The offending edge goes with the
orientation of the plaid arc, and
so we do the $(+)$ dynamics to figure
out what happens to points in
$\Pi_{157}$. That is, we consider
the image $F(\Pi_{157})$. (In the other
case, when the offending edge goes 
against the orientation, we would use
the map $F^{-1}$.)  We prove that
\begin{equation}
\label{follow}
F(\Pi_{157}) \subset \Lambda(\Pi_{58}).
\end{equation}
Here $\Lambda$ is the plaid lattice.

We then check that $\Pi_{58}$ is
grid empty.  This means that and
unit integer square in the plaid
model classified by $\Pi_{58}$
is grid empty.  Here we mean that
$\Phi_{\Pi}$ maps the center of
the square into $\Pi_{58}$.
But then Equation \ref{follow}
tells us that the square on the right
and side of any plaid triple associated
to $\Pi_{157}$ is grid empty.
We have shaded the trid empty squares.

Equation \ref{follow} says that
every occurance of a plaid triple
associated to $\Pi_{157}$ is conjoined,
so speak, with the a plaid triple 
associated to $\Pi_{58}$.  We have also
added in an extra square, even though
we don't know the picture in this square.
Figure
12.5 shows the situation.

\begin{center}
\resizebox{!}{2in}{\includegraphics{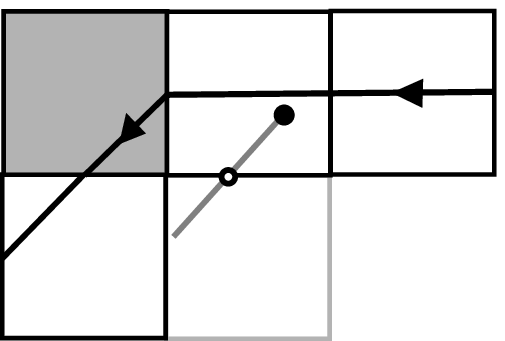}}
\newline
{\bf Figure 12.5:\/} The plaid quadruple obtained
by concatenating $\Pi_{157}$ and $\Pi_{38}$.
\end{center}

When we look at the curve following dynamics
for the other recalcitrant cases, we discover
that the same thing always happens:  The
side square associated to the offending
edge is grid empty.  Up to isometry, the
picture always looks like one of the cases
of Figure 12.6.

\begin{center}
\resizebox{!}{2.6in}{\includegraphics{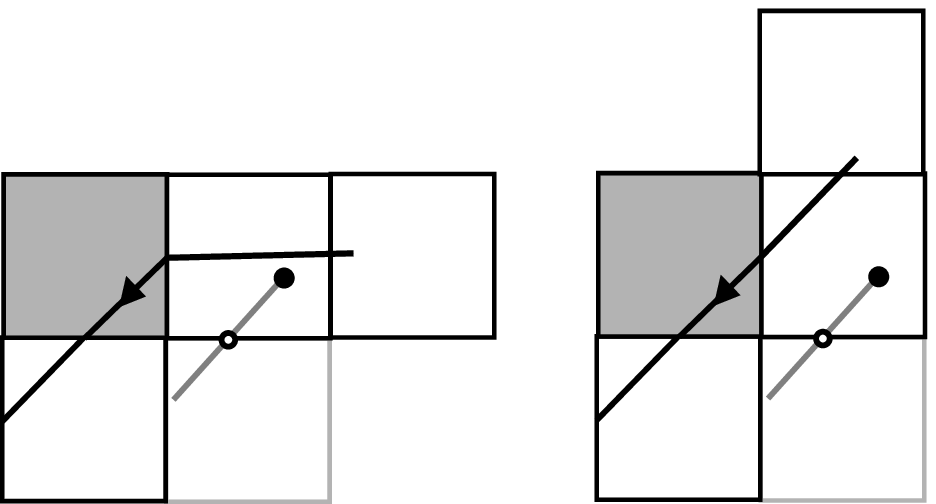}}
\newline
{\bf Figure 12.6:\/} Concatenations for the
recalcitrant triples
\end{center}

The only thing we need to do in order to finish
the proof of Statements 3 and 4 is to analyze
where the offending edge ends.  There are $8$
cases to consider, and the last $4$ cases are
rotates images of the first $4$ cases.  So,
we will just consider the first $4$ cases.

\subsubsection{Case 1}

In this case, the offending edge is
$dT(-1,0)$ and the picture is oriented as
in Figure 12.6.  There are $9$ cases like this.
We have
$$dT(0,-1)=\bigg(-A,\frac{A^2-2A-1}{1+A}\bigg).$$
For each of the $9$ cases, we check that
\begin{equation}
\label{reach}
\Psi(\Pi_k) \subset \{(x,y,z,A)| x \leq A\}.
\end{equation}
By the Reconstruction Formula, the
first coordinate of $(a,b)=\zeta_{\Gamma}$,
the graph grid point contained in the
relevant square, to satisfy
$[a]<A$.  But then the offending edge
must cross over the thick vertical
line shown in Figure 12.7 and end
in one of the two indicated squares.

\begin{center}
\resizebox{!}{3.2in}{\includegraphics{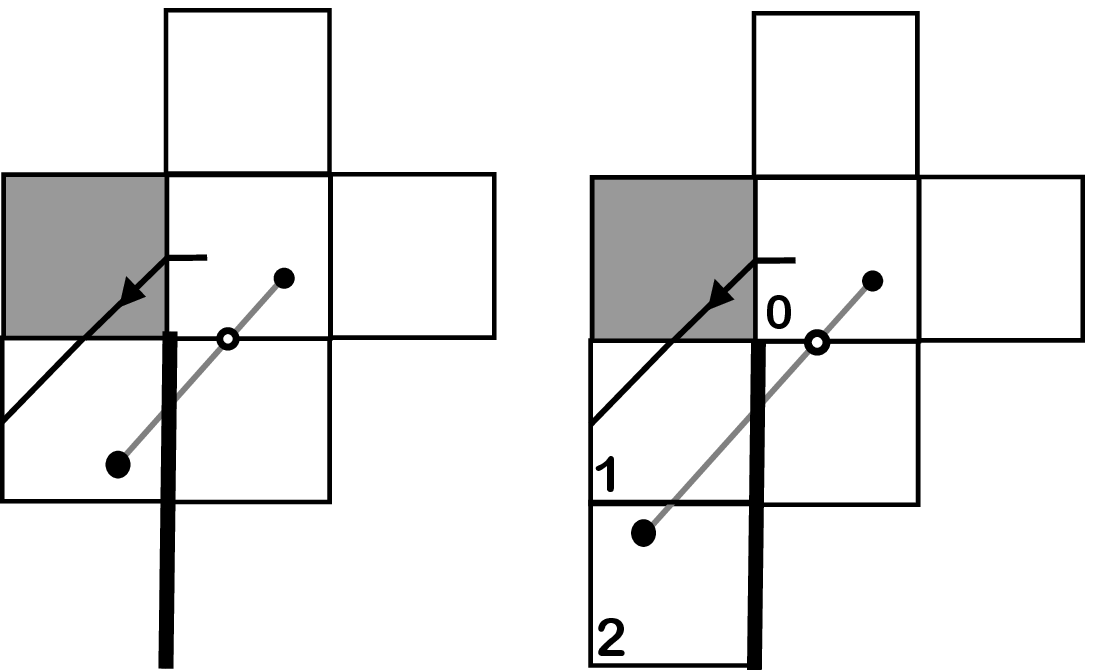}}
\newline
{\bf Figure 12.7:\/} The two possible
endings for the offending edge
\end{center}

We recognize the two cases as
reflected versions of the catches in Figure 4.2.
The proof is done in this case.

\subsubsection{Case 2}

In this case, the offending edge is
$dT(-1,0)$ and the picture is oriented as
in Figure 12.8.  There are $7$ cases like this,
and they all have the features shown in Figure 12.8.

\begin{center}
\resizebox{!}{2in}{\includegraphics{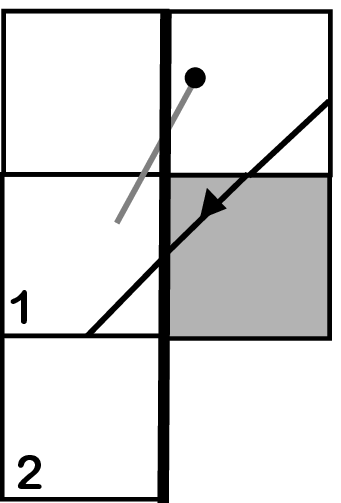}}
\newline
{\bf Figure 12.8:\/} Cast 2
\end{center}

We check that Equation \ref{reach} holds
in all $7$ cases.  Hence, the offending edge
crosses the thick vertical line and
ends in the squares marked $1$ and $2$.
Again, we have the catches shown in 
Figure 4.2.

\subsubsection{Case 3}

In this case, the offending edge is
$dT(0,-1)=(-1,P)$,
 and the picture is
oriented as in Figure 12.9. 
 There are
$3$ cases like this.

\begin{center}
\resizebox{!}{1.6in}{\includegraphics{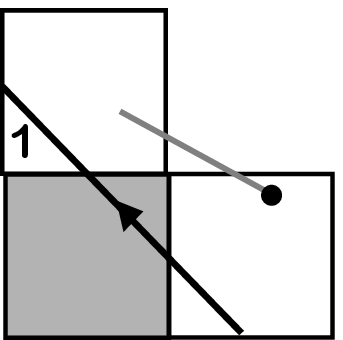}}
\newline
{\bf Figure 12.9:\/} Case 3
\end{center}

This case is easy.  We know that the offending edge
crossed the top of the square it starts in, and
given that the vector is $(-1,P)$, it must end in the
square $\Sigma_1$.  This gives us the left hand
catch in Figure 4.2, up to orientation.

\subsubsection{Case 4}

In this case, the offending edge is
$dT(0,-1)=(-1,P)$,
 and the picture is
oriented as in Figure 12.10. 
 There are
$4$ cases like this.

\begin{center}
\resizebox{!}{1.6in}{\includegraphics{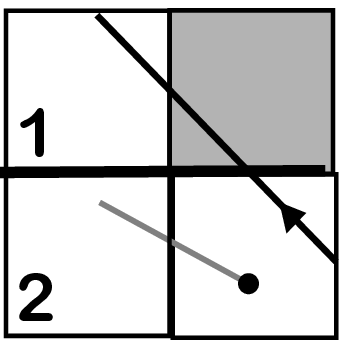}}
\newline
{\bf Figure 12.10:\/} Case 4
\end{center}

In this case, the offending edge crosses
the thick horizontal line provided that
the coordinates $(a,b)$ of $\zeta_{\Gamma}$
satisfy
\begin{equation}
\label{barrier4}
b+Pa \geq 1-P, \hskip 30 pt b \leq 1.
\end{equation}
Notice that this is a weaker condition than
$b+Pa \geq 1$, which is what we would have
needed to solve the crossing problem.
Using the Graph Reconstruction Formula,
we have
\begin{equation}
a=x, \hskip 30 pt b=\frac{y-Ax}{1+A}
\end{equation}
Plugging this equation into Equation \ref{barrier4}
and simplifying, we find that Equation \ref{barrier4}
holds provided that
\begin{equation}
\label{barrier5}
y \in [1-A-Ax,1+A+Ax]
\end{equation}
We check, for each of the cases, that
$x \geq 0$ and $y \in [1-A,1+A]$.
This does the job for us.
Thus, the offending edge ends in the square
$\Sigma_1$, and we get the catch on the
left hand side of Figure 4.2.

In all cases, each offending edge has a catch.
At the same time, each catch that appea

\newpage
\section{Computer Assisted Techniques}

\subsection{Operations on Polytopes}

\noindent
{\bf Clean Polytopes:\/}
Say that a {\it clean polytope\/} is a convex
polytope in $\R^4$ with integer vertices, such
that each vertex is the unique extreme point
of some linear functional. In ther words, a
clean polytope is the convex hull of its
(integer) vertices, and the convex hull of any proper
subset of vertices is a proper subset.
We always deal with clean polytopes.  
The polytopes in the plaid and graph
triple partitions have all the properties
mentioned above, except that their 
vertices are rational rather than integral.
We fix this problem by scaling all polytopes
in all partitions by
a factor of $60$.
\newline
\newline
{\bf Clean Polytope Test:\/}
Suppose we are given a finite number of integer
points in $\R^4$.  Here is how we test that they
are the vertices of a clean polytope.  We consider
all linear functionals of the form
\begin{equation}
\label{linear1}
L(x,y,z,A)=c_1 x + c_2 y + c_3 z + c_4 A, \hskip 30 pt
|c_i| \leq N
\end{equation}
and we wait until we have shown that each vertex is
the unique maximum for one of the functionals.
For the polytopes of interest to us, it suffices
to take $N=3$.   In general, our test halts with
success for some $N$ if and only if the polytope
is clean.
\newline
\newline
{\bf Disjointness Test:\/}
Here is how we verify that two clean polytopes
$P_1$ and $P_2$ have disjoint interiors.
We consider the same linear functionals
as listed in Equation \ref{linear1}
and we try to find some such $L$ with
the property that
\begin{equation}
\max_{v \in V(P_1)}L(v) \leq \min_{v \in V(P_2)} L(v).
\end{equation}
Here $V(P_k)$ denotes the vertex set of $P_k$.
If this happens, then we have found a hyperplane
which separates the one polytope from the other.
This time we take $N=5$.
\newline
\newline
{\bf Containment Test:\/}
Given clean polytopes $P_1$ and $P_2$, here is how
we verify that $P_1 \subset P_2$.  By convexity,
it suffices to prove that $v \in P_2$ for each
vertex of $P_1$.  So, we explain how we verify that
$P=P_2$ contains an integer point $v$.
We do not have the explicit facet structure of $P$,
though for another purposes (computing volumes)
we do find it.

Let $\{L_k\}$ denote the set of all linear
functionals determined by $4$-tuples of vertices
of $P$.  Precisely, Given $4$ vertices of
$P$, say $w_0,w_1,w_2,w_3$, and some integer point
$v$, we take the 
$4 \times 4$ matrix whose first three
rows are $w_i-w_0$ for $i=1,2,3$ and whose
last row is $v$.  Then
\begin{equation}
L_{w_0,w_1,w_2,w_3}(v)=\det(M)
\end{equation}
is the linear functional we have in mind.

If the vertices do not span a $3$-dimensional
space, then $L$ will be trivial.  This does
not bother us.  Also, some choices of vertices
will not lead to linear functions which define
a face of $P$.  This does not bother us either.
The point is that our list of linear functionals
contains all the ones which do in fact define
faces of $P$.  We take our vertices in all
orders, to make sure that we pick up every
possible relevant linear functional.  The
computer does not mind this redundancy.

It is an elementary exercise to show that
$v \not \in P$ if and only if $v$ is a
unique extreme point amongst the set
$\{v\} \cup V(P)$ for one of our
linear functionals.  Here $V(P)$ denotes
the vertex set of $P$, as above.
So, $v \in P$ is and only if $v$ is never
a unique extreme point for any of
the linear functionals on our list.
\newline
\newline
{\bf Volume:\/}
First we explain how we find the
codimension $1$ faces
of $P$. We search
for $k$-tuples of vertices which are
simultaneously in general position and
the common extreme points for one of the
linear functionals on our list.  As long
as $k \geq 4$, the list we find will be
the vertices of one of the faces of $P$.

Now we explain how we compute the
volume of $P$.  
This is a recursive problem. Let
$v_0$ be the first vertex of $P$.
Let $F_1,...,F_k$ be the codimension $1$
faces of $P$.  Let $P_j$ be the cone
of $F_j$ to $v_0$.  This is the same
as the convex hull of $F_j \cup \{v_0\}$.
Then 
\begin{equation}
{\rm vol\/}(V)=\sum_{j=1}^k {\rm vol\/}(V_j).
\end{equation}
If $v_0 \in V_j$ then the volume is $0$. These
extra trivial sums do not bother us.

To compute ${\rm vol\/}(V_j)$ we let
$w_{j0},w_{j1},w_{j2},w_{j3}$ be the first $4$ vertices
of $F_j$.  Let $L_j$ be the associated linear
functional.  Then 
\begin{equation}
4 \times {\rm vol\/}(V_j)=L_j(v_0-w_{j0}) \times {\rm vol\/}(F_j).
\end{equation}

We compute ${\rm vol\/}(F_j)$ using the same method,
one dimension down.  That is, we cone all the facets of
$F_j$ to the point $w_{j0}$.  It turns out that the
polyhedra $\{F_{ij}\}$ in the subdivision of $F$ are either
tetrahedra or pyramids with quadrilateral base.
In case $F_{ij}$ is a tetrahedron we compute
$12{\rm vol\/}(F_{ij})$ by taking the appropriate determinant
and doubling the answer.  In the other case, we
compute $6$ times the volume of each of the $4$
sub-tetrahedron of
$F_{ij}$ obtained by omitting a vertex other than $w_{0j}$
and then we add up these volumes.  This computes
$12 {\rm vol\/}(F_{ij})$ regardless of the cyclic ordering
of the vertices around the base of $F_{ij}$. 
 
When we add up all these contributions, we get
$12\ {\rm vol\/}(F_j)$.  So, our final answer
is $48\ {\rm vol\/}(V)$.  The reason we scale things
up is that we want to have entirely integer quantities.
\newline
\newline
{\bf Potential Overflow Error:\/}
With regard to the volume method,
a straightforward application of our
method can cause an overflow error. 
 We do the calculations using longs
(a $64$ bit representation of an integer) and the
total number we get by adding up $218$ smaller
numbers is (barely) too large to be reliable represented.
However, we observe that all the numbers we
compute are divisible by $480$.  So, before adding each
summand to the list, we divide out by $480$.
This puts us back into the representable range.
Recalling that we have scaled up by $60$, then
computed $48$ times the volume, then divided
by $480$, our final computation is
\begin{equation}
\Omega=6^4 10^3
\end{equation}
 times the true volume.

\subsection{The Calculations}

We work with $3$ partitions and then a few extra polytopes.
\newline
\newline
\noindent
{\bf The Plaid Partition:\/}
The plaid partition, described in
\S \ref{PLAID} has $26$ clean
polytopes modulo the action of the
plaid lattice $\Lambda_1$.
For each all $i,j,k$ with
$|i|,|j|,|k| \geq 3$, we check
that each polytope $P_k$ is
disjoint from $\lambda_{i,j,k}(P_{\ell})$,
where $\lambda_{i,j,k}$ is the
standard word in the generators of
$\Lambda_1$.  Given that all the
original polytopes intersect the
fundamental domain, and given the
sizes of the translations in
$\Lambda_1$, this check suffices
to show that the orbit
$\Lambda(\bigcup P_i)$ consists of
polytopes with pairwise disjoint
interiors and that the union
of our polytopes $P_0 \cup ... \cup P_{25}$
is contained in a fundamental domain
for $\Lambda_1$.  

At the same time, we compute that
\begin{equation}
\sum_{i=0}^{25} {\rm vol\/}(P_i)=8.
\end{equation}
Again, we compute that the scaled
volume is $8\Omega$.
This coincides with the volume of
the fundamental domain $[-1,1]^3 \times [0,1]$
for $\Lambda_1$. 

We now know the following three things:
\begin{enumerate}
\item The union of the $26$ plaid polytopes is
contained in a fundamental domain.
\item The union of the $26$ plaid polytopes
has the same volume as a fundamental domain.
\item The $\Lambda_1$ orbit of the $26$
plaid polytopes consists of polytopes
having pairwise disjoint interiors.
\end{enumerate}
We conclude from this that the $\Lambda_1$ orbit
of the plaid polytopes $P_0,...,P_{25}$ is a
partition of $\widehat X$, as desired.
\newline
\newline
{\bf Remark:\/} In view of the fact that we have
already proved the Plaid Master Picture Theorem,
these checks are really unnecessary, provided that
we have correctly interpreted the geometric
description of the plaid polytopes and copied
down the points correctly.  So, these calculations
really serve as sanity checks.
\newline
\newline
{\bf The Graph Partitions:\/}
We just deal with the $(+)$ graph partition, because
the $(-)$ graph partition is isometric to the
$(+)$ partition.  We make all the same calculations
for the $(+)$ graph partition that we made for
the plaid partition, and things come out the
same way.  Once again, the fact that things work
out is a consequence of our Master Picture Theorem
from [{\bf S0\/}], but we made several changes
to the polytopes in this paper.  We swapped
the first and third coordinates, and also
translated.  So, these calculations serve
as sanity checks.
\newline
\newline
{\bf The Reduced Plaid Triple Partition:\/}
The reduced plaid triple partition consists
of $218$ convex rational polytopes, which
we scale up by a factor of $60$. The
scaled polytopes are all clean.  Using the
above tests, we check that
the polytopes have pairwise disjoint interiors,
and are all contained in the
fundamental domain $X_{\Pi}=[-1,1]^3 \times [0,1].$
Finally, we check that the sum of the volumes is $8$.
This verifies that we really do have a partition.

There is one more important check we make.
Each plaid triple polytope has a $6$ letter
code which tells how it is obtained as
an intersection of the form
$A_{-1} \cap A_0 \cap A_1$, where
$A_k$ is a polytope in the
partition $F^{(k)}({\cal P\/})$.
Here $\cal P$ is the plaid partition
and $F$ is the curve-following
dynamics. 

We check that
each triple plaid $\Pi_k$ is contained
in the three polytopes that are supposed
to contain it, and is disjoint from all
the others in the plaid partition and
its images under the forward and
backward curve following map.
This verifies that the plaid triple
partition really is as we have defined it.
\newline
\newline
{\bf Nine Graph Doubles:\/}
We mentioned in \S \ref{CP} that sometimes
it happens that $\Psi(\Pi_k)$ is not
contained in a single graph polytope of
one of the two partitions, but rather
a union of two of them. This happens
$19$ times for the $(+)$ partition
and $19$ times for the $(-)$
partition.  We call these
unions of graph polytopes
{\it doubles\/}.

We check that each graph double 
(when scaled up by a factor of $60$) is
clean, that each graph double
indeed contains both of its
constituent graph polytopes, and that
the volume of the graph double is the sum of
the volumes of the two constituents. 
This shows that the union of the two
graph polytopes really is a clean convex
polytope.
\newline
\newline
{\bf The Rest of the Calculations}
Each edge crossing problem either involves showing
that some linear functional is positive on a
polytope, or else that two clean polytopes have
pairwise disjoint interiors, or that one
clean polytope is contained 
in another.  We simply run the tests and get
the outcomes mentioned above.  The same goes
for the several recalcitrant cases done in
connection with the proof of Statement 5 of
the Pixellation Theorem given in the last chapter.
Finally, checking that there are no errant edges
just amounts to listing out the data and checking
that there are no forbidden edge assignments.
We also survey the assignments visually and see
that there are no errant edges.

\subsection{The Computer Program}
\label{COMP}

In this last section I'll give a rough account
of some of the main features of the program.
The main purpose of this account is to show you
the kinds of things the program can do.  The
later entries in this section are rather sketchy.
They are designed to let you know what sorts of
things the program can do, but they stop short
of giving detailed instructions on how to get
the program to do it.  The program has its own
documentation - every feature is explained - and
this should help with the details.

\subsubsection{Downloading the Program}

My computer program can be downloaded from 
\newline
\newline
{\bf http://www.math.brown.edu/$\sim$res/Java/PLAID2.tar\/}
\newline
\newline
When you download the file, you get a tarred
directory. 
I untar the directory with the 
Unix command {\bf tar -xvf PLAID2.tar\/}.
(Your system might be different.)
Once this is done,  you have a new directory
called PlaidModel.  The program resides in this
directory and is spread out among many files.

\subsubsection{Running and Compiling}

The file {\bf Main.java\/} is the main file.
Assuming that the program is compiled already
(and you would know by the presence of many .class files)
compile the program with the command
{\bf javac *.java\/} and then you run the program
with the command {\bf java Main\/}.
All this assumes, of course, that your computer can run
Java programs. 
If everything works, a small and colorful
window should pop up.
This is the control panel. You can launch the other parts
of the program from this window.

\subsubsection{What to do first}

The control panel has
a smaller window which lists $10$ pop-up windows.
If you click on these buttons, additional windows
will pop up.  The first button you should press is
the {\bf Document\/} button. This will bring up a
window which has information about the program.
If you read the documentation, you will see how
to operate the program.

\subsubsection{Presets}

The program has 4 preset modes.  When you
select the {\bf preset\/} option on the
{\bf main\/} control panel, you bring
up an auxiliary control panel called
{\bf presets\/}.  This panel has
4 buttons:
\begin{itemize}
\item quasi-isomorphism theorem
\item plaid master picture theorem
\item graph master picture theorem
\item plaid-graph correspondence
\end{itemize}
If you press one of these buttons, various
windows will pop up, and they will be
automatically set up to best show the
advertised features.  Moreober, documentation
will appear which gives further instructions
and explanations.

There is one irritating feature of the program
which I should mention.  The {\bf preset\/}
features work best when all the auxiliary
pop-up windows are closed.  If you press
a {\bf preset\/} when some the pop-up windows
are already open, you run the risk of
having duplicate windows open, and this
causes problems for the program.
When you want to use one of the preset
buttons, you should close all the
auxiliary windows.

\subsubsection{Surveying the Data and Proofs}

The files starting {\bf Data\/} contain
all the data for the polytope partitions.
The file names give some idea of
what the files contain. For instance,
{\bf DataGraphPolytopes.java\/} contains
the coordinates of the graph PET polytopes.
Some of the files are harder to figure
just from the names, but the files
themselves have documentation.

The files starting {\bf Proof\/} contain
all the routines for the proof.
Again, the names indicate the
tests contained in the files. For
instance {\bf ProofVolume.java\/}
contains the volume calculations
for the partitions.

\newpage

\section{References}

[{\bf B\/}] P. Boyland, {\it Dual billiards, twist maps, and impact oscillators\/},
Nonlinearity {\bf 9\/}:1411--1438 (1996).
\newline
\newline
[{\bf DeB\/}] N. E. J. De Bruijn, {\it Algebraic theory of Penrose's nonperiodic tilings\/},
Nederl. Akad. Wentensch. Proc. {\bf 84\/}:39--66 (1981).
\newline
\newline
[{\bf D\/}] R. Douady, {\it These de 3-eme cycle\/}, Universit\'{e} de Paris 7, 1982.
\newline
\newline
[{\bf DF\/}] D. Dolyopyat and B. Fayad, {\it Unbounded orbits for semicircular
outer billiards\/}, Annales Henri Poincar\'{e}, to appear.
\newline
\newline
[{\bf G\/}] D. Genin, {\it Regular and Chaotic Dynamics of
Outer Billiards\/}, Pennsylvania State University Ph.D. thesis, State College (2005).
 \newline 
\newline
[{\bf GS\/}] E. Gutkin and N. Simanyi, {\it Dual polygonal
billiard and necklace dynamics\/}, Comm. Math. Phys.
{\bf 143\/}:431--450 (1991).
\newline
\newline
[{\bf H\/}] W. Hooper, {\it Renormalization of Polygon Exchange Transformations
arising from Corner Percolation\/}, Invent. Math. {\bf 191.2\/} (2013) pp 255-320
\newline
\newline
[{\bf Ko\/}] Kolodziej, {\it The antibilliard outside a polygon\/},
Bull. Pol. Acad Sci. Math.
{\bf 37\/}:163--168 (1994).
\newline
\newline
[{\bf M1\/}] J. Moser, {\it Is the solar system stable?\/},
Math. Intelligencer {\bf 1\/}:65--71 (1978).
\newline
\newline
[{\bf M2\/}] J. Moser, {\it Stable and random motions in dynamical systems, with
special emphasis on celestial mechanics\/},
Ann. of Math. Stud. 77, Princeton University Press, Princeton, NJ (1973).
\newline
\newline
[{\bf N\/}] B. H. Neumann, {\it Sharing ham and eggs\/},
Summary of a Manchester Mathematics Colloquium, 25 Jan 1959,
published in Iota, the Manchester University Mathematics Students' Journal.
\newline
\newline
[{\bf S0\/}] R. E. Schwartz, {\it Outer Billiard on Kites\/},
Annals of Math Studies {\bf 171\/} (2009)
\newline
\newline
[{\bf S1\/}] R. E. Schwartz, {\it Inroducing the Plaid Model\/},
preprint 2015.
\newline
\newline
[{\bf S2\/}] R. E. Schwartz, {\it Unbounded orbits for the plaid model\/},
preprint 2015.
\newline
\newline
[{\bf T1\/}] S. Tabachnikov, {\it Geometry and billiards\/},
Student Mathematical Library 30,
Amer. Math. Soc. (2005).
\newline
\newline
[{\bf T2\/}] S. Tabachnikov, {\it Billiards\/}, Soci\'{e}t\'{e} Math\'{e}matique de France, 
``Panoramas et Syntheses'' 1, 1995
\newline
\newline
[{\bf VS\/}] F. Vivaldi and A. Shaidenko, {\it Global stability of a class of discontinuous
dual billiards\/}, Comm. Math. Phys. {\bf 110\/}:625--640 (1987).

\end{document}